\documentclass[11pt]{amsart}
\usepackage[hmargin=3cm,vmargin=3cm]{geometry}
\usepackage[backend=bibtex,doi=false,isbn=false,url=false,style=alphabetic,maxnames=6]{biblatex}
\usepackage{etex}
\bibliography{../../bibliography/bib}
\AtBeginBibliography{\small}

\usepackage[latin1]{inputenc}
\usepackage{amssymb,amsfonts}
\usepackage{amsthm,amsmath,amscd,stmaryrd,latexsym,youngtab,mathtools,pinlabel,comment}
\usepackage[scaled]{helvet} % ss
\usepackage{courier} % tt

\usepackage{combelow}  %%%%% for some Romanian and Turkish characters.

\usepackage{euscript,mathrsfs}

\usepackage[nohug]{diagrams}
\usepackage{tikz, tikz-cd}
\usetikzlibrary{matrix, arrows, calc, fit}

\newcommand{\ig}[2]{\includegraphics[scale=#1]{diagrams/#2}}

\RequirePackage{color}
\definecolor{myred}{rgb}{0.75,0,0}
\definecolor{mygreen}{rgb}{0,0.5,0}
\definecolor{myblue}{rgb}{0,0.25,0.65}
\definecolor{references}{rgb}{0,0,1}
  \RequirePackage[pdftex,
   colorlinks = true,
   urlcolor = references, % \href{...}{...} external (URL)
   citecolor = references, % \cite{...}
   linkcolor = references, % \ref{...} and \pageref{...}
   ]
 {hyperref}

\newtheorem{theorem}{Theorem}[section]
\newtheorem{lemma}[theorem]{Lemma}

\newtheorem{proposition}[theorem]{Proposition}

\newtheorem{corollary}[theorem]{Corollary}

\newtheorem{conjecture}[theorem]{Conjecture}

\newtheorem{observation}[theorem]{Observation}

\theoremstyle{definition}

\newtheorem{definition}[theorem]{Definition}

\newtheorem{notation}[theorem]{Notation}

\newtheorem{example}[theorem]{Example}

\theoremstyle{remark}
\newtheorem{remark}[theorem]{Remark}

\numberwithin{equation}{section}

    \def\SM{{\mathbb{S}}}

\def\CB{{\mathbf C}}    \def\CC{{\mathcal{C}}}

\def\JB{{\mathbf J}}    
\def\KB{{\mathbf K}}    \def\KC{{\mathcal{K}}}

    \def\OC{{\mathcal{O}}}
\def\PB{{\mathbf P}}    \def\PC{{\mathcal{P}}}
\def\QB{{\mathbf Q}}

\def\XB{{\mathbf X}}

\def\CS{{\EuScript C}}

\def\a{\alpha}
\def\b{\beta}
\def\g{\gamma}

\def\d{\delta}

\def\e{\varepsilon}

\def\l{\lambda}
\def\L{\Lambda}

\let\phi=\varphi
\let\tilde=\widetilde

\usepackage{bbm}
\def\C{{\mathbbm C}}
\def\N{{\mathbbm N}}
\def\R{{\mathbbm R}}
\def\Z{{\mathbbm Z}}
\def\Q{{\mathbbm Q}}
\def\1{\mathbbm{1}}

\newcommand{\F}{\mathbb{F}}

\newcommand{\refequal}[1]{\xy {\ar@{=}^{#1}
(-1,0)*{};(1,0)*{}};
\endxy}

\newcommand{\ip}[1]{\langle #1 \rangle}

\newcommand{\inv}{^{-1}}
\renewcommand{\setminus}{\smallsetminus}
\renewcommand{\d}{\delta}
\renewcommand{\emptyset}{\varnothing}

\newcommand{\one}{\mathbbm{1}}

\newcommand{\xx}{\mathbf{x}}

\newcommand{\Hom}{\operatorname{Hom}}

\renewcommand{\deg}{\operatorname{deg}}

\newcommand{\Tot}{\operatorname{Tot}}

\newcommand{\Cone}{\operatorname{Cone}}

\newcommand{\Tr}{\operatorname{Tr}}

\newcommand{\Ext}{{\rm Ext}}

\newcommand{\Ch}{\textrm{Ch}}

\newcommand{\SBim}{\SM\textrm{Bim}}

\newcommand{\Br}{\textrm{Br}}

\renewcommand{\min}{\textrm{min}}

\newcommand{\corr}{\operatorname{corr}}
\newcommand{\tdinv}{\operatorname{dinv}}

\newcommand{\Dyck}{\operatorname{Dyck}}
\newcommand{\PF}{\operatorname{PF}}
\newcommand{\Seq}{\operatorname{Seq}}

\newcommand{\des}{\operatorname{des}}
\newcommand{\ides}{\operatorname{ides}}
\newcommand{\JO}{\operatorname{JO}}

\newcommand{\Cox}{X}
\newcommand{\otherCox}{Y}
\newcommand{\Sym}{\operatorname{Sym}}
\newcommand{\ev}{\operatorname{ev}}
\newcommand{\Hilb}{\operatorname{Hilb}}
\renewcommand{\succ}{\operatorname{succ}}
\newcommand{\FT}{\operatorname{FT}}
\newcommand{\HT}{\operatorname{HT}}

\newcommand{\dinv}{\operatorname{dinv}}
\newcommand{\Dinv}{\operatorname{Dinv}}
\newcommand{\area}{\operatorname{area}}
\newcommand{\range}{\operatorname{range}}
\newcommand{\f}{f^{(m)}}
\newcommand{\DEnd}{\operatorname{DEnd}}
\newcommand{\DHom}{\operatorname{DHom}}

\newcommand{\Coh}{\operatorname{Coh}}

\newcommand{\idemp}{\PB}

\newcommand{\HHH}{\operatorname{HHH}}
\newcommand{\HH}{\operatorname{HH}}
\newcommand{\Hecke}{\mathcal{H}}

\renewcommand{\setminus}{\smallsetminus}

\begin{document}

%\begin{abstract}
%We introduce a new method for computing triply graded link homology.  Our main application gives an exact computation for all $(n,n)$-torus links, and we verify that our answer agrees with Gorsky's magic formula, for $1\leq n\leq 4$.   We expect that our method can be adapted to many other families of links, particularly torus links. 
%\end{abstract}

\begin{abstract}
We give a simple recursion which computes the triply graded Khovanov-Rozansky homology of several infinite families of knots and links, including the $(n,nm\pm 1)$ and $(n,nm)$ torus links for $n,m\geq 1$.  We interpret our results in terms of Catalan combinatorics, proving a conjecture of Gorsky's.  Our computations agree with predictions coming from Hilbert schemes and rational DAHA, which also proves the Gorsky-Oblomkov-Rasmussen-Shende conjectures in these cases.  Additionally, our results suggest a topological interpretation of the symmetric functions which appear in the context of the $m$-shuffle conjecture of Haglund-Haiman-Loehr-Remmel-Ulyanov.
\end{abstract}

\title{Khovanov-Rozansky homology and higher Catalan sequences}

\author{Matthew Hogancamp} \address{University of Southern California}\thanks{Partially supported by NSF grant DMS-1255334.}

\maketitle

\setcounter{tocdepth}{1}
\tableofcontents

\section{Introduction}

%In \cite{DGR06} it was conjectured that there should exist a triply graded homology theory for links in $\R^3$, which categorifies the HOMFLY-PT polynomial.  In addition, it was predicted that this triply graded link homology should have a remarkable amount of structure, including certain differentials which collapse the triply graded theory to various doubly graded theories constructed earlier by Khovanov \cite{Kh00}, Khovanov-Rozansky \cite{KR08}, and others.  The authors then showed that this additional structure overdetermines the invariants of many knots, thereby resulting in remarkable ``computations'' of the (at that time nonexistent) triply graded homology of these knots.

%In \cite{KR08b}, Khovanov-Rozansky succeeded in constructing a triply graded link homology theory which categorifies the HOMFLY-PT polynomial.  J.~Rasmussen then proved \cite{Ra10} that triply graded theory can be made into the $E_1$ page of various spectral sequences which limit to the doubly graded $\sl_N$ Khovanov-Rozansky homologies.  This is somewhat weaker than the predictions in \cite{DGR06}, which amount to the collapsing of these spectral sequences at $E_2$.  The weaker form is unfortunately not enough to validate to the computational predictions in \cite{DGR06}.  

In \cite{KR08b}, Khovanov-Rozansky constructed a triply graded link homology theory which categorifies the HOMFLY-PT polynomial.   This link invariant is the subject of numerous conjectures which suggest intimate connections with Hilbert schemes, \cite{ORS12,GorNeg15,GNR16}, rational Cherednik algebras \cite{GORS12}, and mathematical physics \cite{AgShak15, OblNawata16}.  Historically the first in this circle of conjectures was formulated by Gorsky \cite{GorskyCatalan}; it is particularly important for the present work.

%This is the first instance of what is now recurring theme in the triply graded theory\footnote{We should mention that very recently Oblomkov-Rozansky constructed a potentially different triply graded link homology .  Despite numerous  conjectured formulas for the homologies of various links, the difficulty of performing actual calculations of the homology make a direct verification of these conjectures impossible in all but the simplest cases.  

%To illustrate what we mean, consider the family of torus knots.  At present time there are several beautiful conjectures expressing the triply graded homology of the $(n,m)$ torus knot $T_{n,m}$ in terms of Hilbert schemes \cite{ORS12}, representations of rational Cherednik algebras \cite{GORS12}, other Hilbert schemes \cite{GorNeg15}, and Catalan combinatorics (Gorsky's appendix to \cite{ORS12}).  

\begin{conjecture}[Gorsky's conjecture]\label{conjecture:gorsky}
The minimal $a$-degree part of the Poincar\'e series of the triply graded homology of the $(n,n+1)$ torus knot equals $C_n(q,t)$, the $q,t$ Catalan polynomial \cite{GarHaim96,GarHag-Catalan}.
\end{conjecture}

%At around the same time Oblomkov-Shende \cite{ObSh12} noticed that the HOMFLY-PT polynomials of torus links (and, more generally, algebraic links) can be computed from certain Hilbert schemes of points.  The paper \cite{ORS12} lifts \cite{ObSh12} to a family of conjectures concerning the triply graded homology of algebraic links, generalizing Gorsky's conjecture for $T_{(n,n+1)}$.  In a related but different direction, conjectures in \cite{GORS12} conjecturally extract the homologies of torus knots from representations of rational Cherednik algebras.  Very exciting recent work \cite{GorNeg15,GNR16} relates link homology to Hilbert schemes on a more categorical level, and brings arbitrary knots and links into the picture.

One major obstacle to proving any of the above conjectures is the fact that Khovanov-Rozansky homology is notoriously difficult to compute from the definitions.  That said, stable versions were proven by the author and M.~Abel \cite{Hog15,AbHog15}.  Recently, the author and Ben Elias realized \cite{ElHog16a} that the categorified Young symmetrizers from \cite{Hog15} can be used very effectively to compute the triply graded homologies of some families of links.  Our method seems particularly well adapted to the case of torus links.  In \cite{ElHog16a} we demonstrate this by giving a recursive formula for the Poincar\'e series of several links, including the $(n,n)$ torus links.

In this paper we utilize the same technique to compute the triply graded homologies of several more infinite families of links, including the $(n,nm+r)$ torus links with $r\in \{0,1,-1\}$ and $n,m\geq 1$.   Our flagship result, stated in the language of symmetric functions, is the following.

\begin{theorem}\label{thm:nablaEn}
The Poincar\'e series of the triply graded homology of the $(n,nm+1)$ torus knot equals
\[
\frac{1}{1-q} \sum_{k=0}^n \ip{\nabla^m e_n, h_ke_{n-k}}a^k.
\]
In particular, the $a^0$ coefficient is the higher $q,t$ Catalan polynomial $C_n^{(m)}(q,t)$ \cite{GarHaim96}, and the coefficient of $a^k$ is the higher $q,t$ Schr\"oder polynomial \cite{Song05}.
\end{theorem}
Here $e_i$ and $h_i$ are the elementary and complete symmetric functions and $\nabla:\L_{q,t}\rightarrow \L_{q,t}$ is the Garsia-Bergeron symmetric function operator \cite{BerGars96}, where $\L_{q,t}:=\Q(q,t)[x_1,x_2,\ldots]^{\Sym}$, and $\ip{-,-}$ is the Hall inner product, in which the Schur functions are orthonormal.

From \cite{HHLRU05} we know that $\nabla^m e_n$ is the character of the spaces of generalized diagonal coinvariant rings, which in turn are exactly the DAHA representations which conjecturally correspond to the $(n,nm+1)$ torus knots \cite{GORS12}.  Thus, our theorem proves the Gorsky-Oblomkov-Rasmussen-Shende conjecture for these knots.

We prove this theorem by introducing a special family of complexes of Soergel bimodules, and computing their homologies recursively.  We then compare our formulas with known formulas for $\nabla^m e_n$, utilizing the  $m$-shuffle conjecture \cite{HHLRU05}, recently proven by Carlsson and Mellit \cite{CarMel-pp,MellitRational-pp}.

\subsection{The recursions}
\label{subsec:mainresults}

The results of our computations follow.  Below, we use the phrase \emph{Khovanov-Rozansky (KR) series} to refer to the Poincar\'e series of the triply graded homology $\HHH(C)$ of a  complex $C$ of Soergel bimodules, and we denote it by $\PC_C$.  Given $C\in \KC^b(\SBim_n)$, the KR series $P_C(Q,A,T)$ is a actually polynomial in $A,T,T\inv,Q,Q\inv$, and $(1-Q^2)\inv$. We will almost always express the KR series $\PC_C$ in terms of the variables $q=Q^2$, $t=T^2Q^{-2}$, $a=AQ^{-2}$.  See \S \ref{subsec:soergel} for details.

Given a braid  $\b\in \Br_n$, and let $F(\b)\in \KC^b(\SBim_n)$ be the associated Rouquier complex.  The KR series of $F(\b)$ will also be referred to as the KR series of $\b$; $\PC_{F(\b)}$ depends only on the oriented link $\hat\b$ obtained by closing $\b$ (up to normalization), hence we also refer to this as the KR series of $\hat\b$ by abuse.  %Elsewhere \cite{DGR06}, the KR series of $\hat\b$ is called the \emph{superpolynomial} of $\hat\b$; we prefer to avoid this highly abused term.
 
Fix integers $n,m\geq 1$.  In \S \ref{subsec:theComplexes} we construct some special complexes of Soergel bimodules $\CB_v^{(m)}\in \KC^b(\SBim_n)$, indexed by sequences $v\in \{0,1,\ldots,m\}^n$.  In case $m=1$, we recover the complexes $D_v$ from \cite{ElHog16a}, up to conventions.  We refer to \S \ref{subsec:theComplexes} for the definition, but remark that $\CB^{(m)}_{m^n} = \FT_n^m$, where $\FT_n$ is the (Rouquier complex of the) positive full-twist braid.

\begin{theorem}\label{thm:linkCase_intro}
Fix integers $n,m\geq 1$, and let $f_{v}=f_v^{(m)}\in \Z[q,t,a,(1-q)\inv]$ denote the polynomials, indexed by $v\in\{0,1,\ldots,m\}^n$, defined by the following recursion.
\begin{enumerate}\setlength{\itemsep}{4pt}
\item[(L0)] $f_{\emptyset} = 1$
\item[(L1)] $f_{0,v} = (t^{\#\{i\:|\: v_i<m\}}+a)f_{v}$
\item[(L2)] $f_{k,v}=t^{\#\{i\:|\: v_i<k\}}f_{v,k-1}$ if $1\leq k\leq m-1$.
\item[(L3)] $f_{m,v}=f_{v,m-1}+q f_{v,m}$
\end{enumerate}
for all $v\in \{0,1,\ldots,m\}^{n-1}$.  Then $f_{m^n}$ is the KR series of the $(n,nm)$ torus link and, in  general, $f_v^{(m)}$ is the KR series of $\CB_v^{(m)}\in \KC^b(\SBim_n)$.
\end{theorem}
Note that setting $m=1$, we recover the recursions in \cite{ElHog16a} up to a difference in conventions.  The recursion below is new, even for $m=1$.

\begin{theorem}\label{thm:knotCase_intro}
Fix integers $m,n\geq 1$.  Let $g_{v}=g_v^{(m)}\in \Z[q,t,a,(1-q)\inv]$ denote the polynomials, indexed by $v\in\{0,1,\ldots,m\}^n$, defined by the following recursion.
\begin{enumerate}\setlength{\itemsep}{4pt}
\item[(K0)] $g_{\emptyset} = 1$ and $g_{0}=1+a$
\item[(K1)] $g_{0,v} = (t^{{\#\{i\:|\: v_i<m\}}}+a)g_{v}$ if $0\leq v_1\leq m-1$, and $g_{0,m,v}=t^{\#\{i\:|\: v_i<m\}}g_{v,m-1}$
\item[(K2)] $g_{k,v}=t^{\#\{i\:|\: v_i<k\}}g_{v,k-1}$ if $1\leq k\leq m-1$
\item[(K3)] $g_{m,v}=g_{v,m-1}+qg_{v,m}$
\end{enumerate}
for all sequences $v$.  Then $g_{m^n}$ is the KR series of the $(n,nm+1)$ torus knots, and in general $g_v^{(m)}$ is the KR series of $\XB\CB^{(m)}_v$, where $\XB$ is the Rouquier complex associated to $\sigma_{n-1}\cdots\sigma_2\sigma_1$ with a shift $(n-1)[1-n]$ applied.
\end{theorem}
Here, $\sigma_i\in \Br_n$ denotes the $i$-th elementary braid, which is a positive crossing between the strands $i$ and $i+1$.  In fact our work shows that the homologies $\HHH(\CB_v^{(m)})$ and $\HHH(\XB\CB_v^{(m)})$ have no $\Z$-torsion, hence the iabove polynomials determine the corresponding Khovanov-Rozansky homology groups over $\Z$, up to isomorphism.

Below we state some interesting special cases of the above.

\begin{example}
Let $X_n=\sigma_{n-1}\cdots \sigma_2\sigma_1$, and let $\FT_n= X_n^n$ denote the full twist braid.  Then the KR series of $\FT_{n-1}^{m-k}\FT_n^k$ equals $\frac{1}{1-q}f_{k,m^{n-1}}^{(m)}(q,t,a)$, and the KR series of $X_n\FT_{n-1}^{m-k}\FT_n^k$ equals $\frac{1}{1-q}g_{k,m^{n-1}}^{(m)}(q,t,a)$.
\end{example}
 In \cite{GNR16}, conjectural formulas are given for the KR series of the braids
\[
\FT_2^{m_2} \FT_3^{m_3}\cdots \FT_n^{m_n} \ \ \ \ \ \ \ \text{ and } \ \ \ \ \ \ \ \ X_n \FT_2^{m_2}\FT_3^{m_3}\cdots \FT_n^{m_n}
\]
for all $m_2,\ldots,m_n\geq 0$.  It would be very interesting to compare our formulas to these.

\begin{example}
Let $J_k = \sigma_{k-1}\cdots\sigma_2\sigma_1^2\sigma_2\cdots \sigma_{k-1}$ denote the \emph{Jucys-Murphy braid}.  Then the KR series of $\FT_n^m J_i\inv$ equals $\frac{1}{1-q}f^{(m)}_{m^{n-i},m-1,m^{i-1}}(q,t,a)$.
\end{example}

\begin{example}
The KR series of the braid
\begin{equation}\label{eq:hookBraidIntro}
 \FT_n^m \ \ \cdot\ \ \begin{minipage}{1.05in}
\labellist
\small
\pinlabel $\underbrace{\ \ \ \ \ \ \ \ \  \ \ \ \ \   }_{i-1}$ at 37 -9
\pinlabel $\underbrace{\ \ \ \ \ \ \ \ \ \ }_{n-i}$ at 81 -9
\endlabellist
\includegraphics[scale=.8]{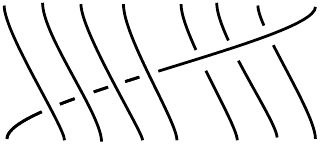}
\end{minipage}\ \
\end{equation}\vskip10pt\noindent
equals $\frac{1}{1-q}g_{m^{n-i},m-1,m^{i-1}}^{(m)}$.  The $(n,nm+1)$ and $(n,nm-1)$ torus knots correspond to $i=1$ and $i=n$, respectively.
\end{example}

\begin{example}
In general $f^{(m)}_{m^{n-i},k,m^{i-1}}(q,t,a)$ is $1-q$ times the KR series of the closure of the braid
\[
\b_{i,k,m,n}:=\sigma_{n-i}\cdots \sigma_2\sigma_1(\one_1\sqcup \FT_{n-1}^{m-k})\FT_n^k \sigma_1\sigma_2\cdots \sigma_{n-i},
\]
and $g^{(m)}_{m^{n-i},k,m^{i-1}}(q,t,a)$ is $1-q$ times the KR series of the closure of $X_n \b_{i,k,m,n}$.
\end{example}

These special cases are discussed further in \S \ref{subsec:specialBraids}.

\subsection{Combinatorial formulas}
\label{subsec:formulasIntro}
Fix integers $m,n\geq 1$, and let $\sigma\in \Z_{\geq 0}^n$ be given.  Let $v=v(\sigma)\in \{0,1,\ldots,m\}^n$ denote the sequence obtained by truncating $\sigma$:
\[
v(\sigma)_i := \begin{cases}
\sigma_i & \text{ if $\sigma_i< m$}\\ m & \text{ if $\sigma_i\geq m$}
\end{cases}.
\]
We say that $\sigma$ \emph{extends} $v$ if $v=v(\sigma)$.  Let $|\sigma|=\sigma_1+\cdots+\sigma_n$.  In \S \ref{subsec:diagrams} we introduce non-negative integers $\tdinv(\sigma)$, $e(v)$, and $d(\sigma,i)$ for each $1\leq i\leq n$.  In \S \ref{subsec:stateSum1} we prove the following.

\begin{theorem}\label{thm:combSumLink_intro}
We have
\[
f_v^{(m)} = \sum_\sigma q^{|\sigma|-|v|}t^{\tdinv(\sigma)+e(v)}\prod_{i=1}^n(1+at^{-d(\sigma,i)})
\]
for all $v\in \{0,1,\ldots,m\}^n$, where the sum is over sequences $\sigma$ which extend $v$.
\end{theorem}
Note that the above is an infinite series.  It is possible to collect terms in such a way so as to obtain a finite sum, similar to what was done in the case $m=1$ by Andy Wilson \cite{Wilson16-pp}, though we do not do so explicitly here.

Given the similarities in the recursions defining $f_v^{(m)}$ and $g_v^{(m)}$, it is perhaps not surprising that there should be a similar summation formula for $g_v^{(m)}$.  The resulting sum occurs over a restricted set of $\sigma$'s, and consequently will end up being finite unless $v=m^n$.

Let us say that $\sigma\in \Z_{\geq 0}^n$ satisfies the \emph{cyclic $m$-Dyck condition} if
\begin{enumerate}
\item[$(i)$] $\sigma_{i+1}\leq \sigma_i+m$ for $1\leq i\leq n-1$.
\item[$(ii)$] $\sigma_1-1\leq \sigma_i+m$.
\end{enumerate}
The reason for this terminology will become clear in \S \ref{sec:shuffle}.  For each index $i\in \{1,\ldots,n\}$, let $\chi_i(\sigma)\in\{0,1\}$ denote 1 if the corresponding inequality above is strict, and zero if there is equality.  In other words, for $1\leq i\leq n-1$ we have $\chi_i(\sigma)=1$ iff $\sigma_{i+1}\leq \sigma_i+m$, and for $i=n$ we have $\chi_n(\sigma)=1$ iff $\sigma_{1}-1< \sigma_n+m$.  The following is proven in \S \ref{subsec:stateSum1}.

\begin{theorem}\label{thm:combSumKnot_intro}
We have
\[
g_v^{(m)} = \sum_\sigma q^{|\sigma|-|v|}t^{\tdinv(\sigma)+e(v)}\prod_{i=1}^n(1+\chi_i(\sigma)at^{-d(\sigma;i)})
\]
for all $v\in \{0,1,\ldots,m\}^n$.  The sum is over sequences $\sigma\in \Z_{\geq 0}^n$ extending $v$, and satisfying the cyclic $m$-Dyck condition.
\end{theorem}
In case $v=(m^{n-1},m-1)$ there is a natural bijection between the terms in the sum formula for $g_{v}^{(m)}$, and $m$-Dyck paths.  Thanks to the recent proof of the shuffle conjecture and its rational version, we now have explicit combinatorial formulas for the $m$-th Catalan $C^{(m)}_n(q,t)$ (see for instance equations (97), (98), and (99) in \cite{HHLRU05}) as a weighted sum over $m$-Dyck paths, and we obtain the following.

\begin{corollary}
We have $g^{(m)}_{m^{n-1},m-1}(q,t,0) = C^{(m)}_n(q,t)$, the higher $q,t$ Catalan polynomial.  
\end{corollary}

This will be generalized and reinterpreted in \S \ref{sec:shuffle}, where we interpret $g_v^{(m)}(q,t,a)$ in terms of pieces of the $m$-shuffle conjecture in \cite{HHLRU05} (now a theorem of Carlsson-Mellit \cite{CarMel-pp} and Mellit \cite{MellitRational-pp,}).

\subsection{Auxilliary results}
\label{subsec:auxResults}
In this section we show how our main results are consequences of some relations which are more general, though a bit more technical to state.  Let $\SBim_n$ denote the category of Soergel bimodules over $\Z$, associated to $S_n$.  For each braid $\b\in \Br_n$, let $F(\b)$ denote the Rouquier complex; this is a complex of Soergel bimodules, i.e.~ an object of the homotopy category $\KC^b(\SBim_n)$.  Let $(1)$ be the grading shift of Soergel bimodules, and $[1]$ be the homological shift of complexes.  We also write $q=(-2)$ and $t=(2)[-2]$.  Let $\XB_n=t^{(1-n)/2}F(\sigma_{n-1}\cdots \sigma_2\sigma_1)$.  When the index $n$ is understood, we will omit it.

The following two relations are proven in \S \ref{subsec:rels1}:

\begin{subequations}
\begin{equation}\label{eq:Ckv_intro}
\XB  \CB_{kv}^{(m)} \XB\inv \ \ \simeq \ \ t^{\#\{i\:|\: v_i<k\}}\CB_{v(k-1)}^{(m)} \ \ \ \ \ \ \ \ \forall \ \ 1\leq k\leq m-1 \ \ \ \ \ \ \  \ \forall \ \ v\in \{0,1,\ldots,m\}^n
\end{equation}
\begin{equation}\label{eq:Cmv_intro}
\XB  \CB_{mv}^{(m)} \XB\inv \ \ \simeq \ \ \Big(\CB_{v(m-1)}^{(m)} \ \ \longrightarrow \ \ q \CB_{vm}^{(m)}\Big) \ \ \ \ \ \ \ \ \ \  \  \ \forall \ \ v\in \{0,1,\ldots,m\}^n
\end{equation}
\end{subequations}
where the notation $B\simeq (C\rightarrow A)$ indicates the existence of a distinguished triangle
\[
A\rightarrow B\rightarrow C\rightarrow A[1].
\]

Iterating these allows us to write any complex $\XB^r \CB_v^{(m)}$ as an iterated mapping cone involving complexes of the form $\XB^r\CB_{0w}^{(m)}$, as well as shifts and conjugates of these.  Thus, for the purposes of computing $\HHH(\XB^r \CB_{v}^{(m)})$ we are reduced to the case $v_1=0$ (modulo the computation of a spectral sequence, which turns out to collapse for degree reasons; see \S \ref{subsec:superpoly}).  How one proceeds from here depends on $r$.  In this paper we take care of the cases $r\in\{0,1\}$.  In \S \ref{subsec:relsII} we prove the following relations.

\begin{subequations}
\begin{equation}\label{eq:C0v_intro}
\HHH(\CB_{0v}^{(m)}) \cong (t^{\#\{i\:|\: v_i<m\}}+a) \HHH(\CB_v^{(m)}) \ \ \ \ \ \ \ \ \ \  \  \ \forall \ \ v\in \{0,1,\ldots,m\}^n
\end{equation} 
\begin{equation}\label{eq:XC0kv_intro}
\HHH(\XB\CB_{0kv}^{(m)}) \cong (t^{1+\#\{i\:|\: v_i<m\}}+a) \HHH(\CB_{kv}^{(m)})\ \ \ \ \ \ \ \ \ \  \  \ \forall \ \ v\in \{0,1,\ldots,m\}^n 
\end{equation}
\begin{equation}\label{eq:XC0mv_intro}
\HHH(\XB\CB_{0mv}^{(m)}) \cong t^{\#\{i\:|\: v_i<m\}} \HHH(\CB_{v(m-1)}^{(m)})\ \ \ \ \ \ \ \ \ \  \  \ \forall \ \ v\in \{0,1,\ldots,m\}^n.
\end{equation}
\end{subequations}
In (\ref{eq:XC0kv_intro}) we assume that $0\leq k\leq m-1$.  In the above formulas we are abusing notation by regarding a polynomial $\Gamma(q,t,a)\in \N[q,t,a]$ as the corresponding functor acting on triply graded abelian groups, which sends a triply graded abelian group $A$ to the appropriate direct sum of shifted copies of $A$.

Theorems \ref{thm:linkCase_intro} and \ref{thm:knotCase_intro} are then consequences of these.

\subsection{Relation to Hilbert schemes and symmetric functions}
\label{subsec:hilb}
Gorsky's original conjecture (Conjecture \ref{conjecture:gorsky}) can now be understood in terms of a beautiful picture that is emerging \cite{GNR16}, relating Soergel bimodules and Hilbert schemes on a categorical level.  Specifically, in \cite{GNR16}, it is conjectured that there exists a pair of adjoint functors
\[
\begin{tikzpicture}
\tikzstyle{every node}=[font=\small]
\node (a) at (0,0) {$\KC^b(\SBim_n)$};
\node (b) at (4,0) {$D^b(\Coh_{\C\times \C} (\Hilb_n(\C^2))$};
\path[->,>=stealth',shorten >=1pt,auto,node distance=1.8cm,
  thick]
([yshift=3pt] a.east) edge node[above] {$\iota_\ast$}		([yshift=3pt] b.west)
([yshift=-2pt] b.west) edge node[below] {$\iota^\ast$}		([yshift=-2pt] a.east);
\end{tikzpicture}
\]
relating $\KC^b(\SBim_n)$ and the ($\C^\times\times \C^{\times})$-equivariant) derived category of $\Hilb_n(\C^2)$ satisfying a number of conjectural properties.  Among these conjectural properties are:
\begin{enumerate}
\item Taking derived (equivariant) global sections of $\iota_\ast C$ gives a bigraded vector space which is isomorphic to $\HHH^0(C)$ (the minimal Hochshcild degree part of $\HHH(C)$).
\item Tensoring with the full twist Rouquier complex on the Soergel side corresponds to twisting by the line bundle $\OC(1)$ on the Hilbert scheme side:
\[
\iota_\ast(\FT_n\otimes C) \cong \OC(1)\otimes \iota_\ast(C).
\]
\item $\iota_\ast(C\otimes D)\cong \iota_\ast(D\otimes C)$.
\end{enumerate}

For instance if $C=X_n$ is the Rouquier complex associated to the positive braid lift of the Coxeter element $s_{n-1}\cdots s_2s_1$, then $\iota_\ast X_n$ is expected to be the trivial line bundle $\OC|_Z$, restricted to the ``punctual'' Hilbert scheme $Z\in \Hilb_n(\C^2)$, hence $\iota_\ast X_n^{nm+1}$ is expected to correspond to $\OC(m)|_Z$.  The graded dimension of the equivariant global sections of this latter object is known to recover the $q,t$ Catalan polynomials by Haiman \cite{HaimanCatalan}, hence one obtains Gorsky's conjecture as a consequence of the Gorsky-Negu{\c t}-Rasmussen (GNR) conjectures.  Thus, our work here provides significant evidence in favoe of the GNR conjecture, by proving that the homology of $\XB^{nm+1}$ agrees with the predictions from the Hilbert scheme side. 

At the moment, there seems to be no known analogue of our complexes $\XB^r\CB_v^{(m)}$ on the Hilbert scheme side, hence in general we are unable to compare the homologies which we compute with existing computations on the Hilbert scheme side.  It would be extremely interesting to find candidate complexes on the Hilbert scheme side.  It seems reasonable to expect that such complexes would play a fundamental role in the study of $\Hilb_n(\C^2)$.

\subsection{Connection to the shuffle conjecture}
\label{subsec:introShuffle}
Passing to $K$-theory gives the following combinatorial shadow of the GNR conjectures.  Associated to each complex $C\in \KC^b(\SBim_n)$ there should be a degree $n$ symmetric function $\Phi(C)=[\iota_\ast(C)]\in \L_{q,t}:=\Q(q,t)[x_1,x_2,\ldots]^{\Sym}$ such that
\begin{enumerate}
\item In good situations the Poincar\'e series of $\HHH(C)$ satisfies
\begin{equation}\label{eq:innerProductKRseries}
\PC_C(q,t,a) = \sum_{k=0}^n \ip{\Phi(C),h_ke_{n-k}}a^k
\end{equation}
where $\ip{-,-}$ is the Hall inner product on $\L_{q,t}$, in which the Schur functions are orthonormal, and $e_i, h_i$ are the elementary and complete symmetric functions.  Here, a ``good'' situation is one in which the complex $\iota_\ast C$ has no higher homology; in this case the Poincar\'e series of the  homology of $\iota_\ast C$ equals the Euler characteristic, which in turn can be computed from the class in $K$-theory.
\item  $[\Phi(\FT_n\otimes C)] = \nabla\iota_\ast (C)$.  Here, $\nabla:\L_{q,t}\rightarrow \L_{q,t}$ is the Garsia-Bergeron operator \cite{BerGars96} on symmetric functions defined in terms of the modified Macdonald basis $\tilde{H}_\mu$ by $\nabla \tilde{H}_\mu = t^{n(\mu)}q^{n(\mu^t)}\tilde{H}_\mu$, where $\mu^t$ is the transpose partition and $n(\mu)=\sum_i (i-1)\mu_i$.  
\item $\Phi(C\otimes D)\cong \Phi(D\otimes C)$.
\item $\Phi(C\sqcup D) = \Phi(C)\Phi(D)$.
\end{enumerate}

Given a bounded complex $C$ of Soergel bimodules, for instance a Rouquier complex, one can now ask what is the symmetric function corresponding $C$ under the GNR conjecture?   It is generally accepted that $X_n$ should correspond to $e_n$, hence $X_n^{nm+1}$ should correspond to $\nabla^m e_n$.  Our work provides significant evidence for this.

Before stating our last main result, we summarize the elements of the shuffle conjecture, referring to \S \ref{sec:shuffle} for the details.  Fix an integer $n\geq 1$.  An $m$-Dyck path is a path in an $n\times nm$ rectangle, from the southwest corner to the northeast corner, consisting of $n$ north steps and $nm$ east steps, which stays weakly above the diagonal.  An $m$-Dyck path $D$ determines and is determined by an associated sequence $\gamma\in \Z_{\geq 0}^n$, where $\gamma_i$ is the horizontal distance from the beginning of the $i$-th north step to the diagonal.  The set of all such sequences arising this way will be denoted  $\Dyck_m$.  To each $\gamma\in \Dyck_m$, there is an explicitly defined symmetric function $D_\gamma(\xx;q,t)\in \L_{q,t}$, which is expressed as a sum over $m$-parking functions.  The $m$-shuffle conjecture/theorem \cite{HHLRU05, CarMel-pp,MellitRational-pp}, states that
\[
\nabla^m e_n = \sum_{\gamma\in \Dyck_m} D_\gamma(\xx;q,t).
\]

We prove the following in \S \ref{subsec:connectionToShuffle}.

\begin{theorem}\label{thm:shuffleKRintro}
Let $v\in \{0,1,\ldots,m\}^n$ be given, let $r=\min\{v_1,\ldots,v_n\}$, and let $i$ be the smallest index with $v_i=r$.  Then
\begin{equation}\label{eq:shufflePoincareSeries_intro}
g_v^{(m)}(q,t,a) = q^{-\corr_q(v)}t^{-\corr_t(v)}\sum_{k=0}^n  \left\langle \sum_{\gamma}\nabla^m D_\gamma(\xx;q,t), h_ke_{n-k} \right\rangle a^k
\end{equation}
where the sum is over $m$-Dyck sequences $\gamma\in \Dyck_m$ such that
\begin{itemize}
\item if $v_j\leq m-1$ and $i+1\leq j\leq n$  then $\gamma_{j+1-i}=v_{j}-r$.
\item if $v_j\leq m-1$ and  $1\leq j\leq i-1$, then $\gamma_{j+1-i+n}=v_{j}-r-1$.
\item if $v_j=m$ and $i+1\leq j\leq n$  then $\gamma_{j+1-i}\geq m-r$.
\item if $v_j=m$  and  $1\leq j\leq i-1$, then $\gamma_{j+1-i+n}\geq m-r-1$.
\end{itemize}
Here $g_v^{(m)}(q,t,a)$ is the Poincar\'e series of $\HHH(\XB\CB_v^{(m)})$, as computed in Theorem \ref{thm:knotCase}.
\end{theorem}
See Theorem \ref{subsec:connectionToShuffle} and the remarks preceding for a definition of the \emph{corrections} $\corr_q(v)$ and $\corr_t(v)$.  In the special case $v=m^{n-1}(m-1)$, we have $\corr_q(v)=\corr_t(v)=0$, hence we obtain Theorem \ref{thm:nablaEn} as a corollary.

\subsection{Outline of the paper}
In \S \ref{sec:background} we summarize the relevant background material: \S \ref{subsec:soergel} recalls Soergel bimodules and sets up our conventions, and \S \ref{subsec:JWproj} recalls some relevant results from \cite{Hog15} and \cite{ElHog16a} regarding categorified Young symmetrizers.

Section \S \ref{sec:computations} is the heart of the paper.  In \S \ref{subsec:theComplexes} we define the complexes $\CB_v^{(m)}$.  In \S \ref{subsec:rels1} and \ref{subsec:relsII} we establish relations satisfied by the $\CB_v^{(m)}$, and in \S \ref{subsec:superpoly} we restate and prove Theorems \ref{thm:linkCase_intro} and \ref{thm:knotCase_intro} from the introduction.  In \S \ref{subsec:specialBraids} we discuss some important special cases of these computations.

In \S \ref{sec:combinatorics} we state and prove the sum formulas.  Section \S \ref{subsec:diagrams} introduces the relevant combinatorial notions.  Armed with this setup, in \S \ref{subsec:stateSum1} and \S \ref{subsec:stateSum2}, we restate and prove Theorems \ref{thm:combSumLink_intro} and \ref{thm:combSumKnot_intro}, respectively.

In \S \ref{sec:shuffle} we relate our results to the $m$-shuffle conjecture.  In \S \ref{subsec:motivation} we give some context and motivation.  In \S \ref{subsec:mDyck} and \S \ref{subsec:parkingFn} we recall the basics of $m$-Dyck paths and $m$-parking functions, respectively, and in \S \ref{subsec:quasisymmetric} we recall the definition of the symmetric functions $D_\gamma(\xx;q,t)$.  In \S \ref{subsec:connectionToShuffle} we restate and prove Theorem \ref{thm:shuffleKRintro} from the introduction, which solidifies the connection between our work and the $m$-shuffle conjecture.  In \S \ref{subsec:specialCasesShuffle} we work out some special cases of this theorem.  Finally, in \S \ref{subsec:moreconnections} we formulate some combinatorial conjectures based on our work.

\subsection{Acknowledgements}
The author would like to thank Eugene Gorsky for the enlightening conversations, and Sami Assaf for teaching him about quasi-symmetric functions.

%%%%%%%%%%%%%%%%%%%%%%%%%%%%%5
\section{Background}
\label{sec:background}
%%%%%%%%%%%%%%%%%%%%%%%%%%%%

In this section we summarize some necessary background, referring to \cite{Hog15} and \cite{ElHog16a} for the details.

%This paper builds off of \cite{ElHog16a}.  We recall the basics in the next couple of sections, but refer to \emph{loc.~cit.} for a proper introduction.

%============================
\subsection{Soergel bimodules}
\label{subsec:soergel}
%===========================

In this section we introduce some background on Soergel bimodules and set up conventions. The precise definition of Khovanov-Rozansky homology is fairly sophisticated, but we remind the reader that all of our computations will utilize only certain formal properties.  As such, the details of the definition are mostly not necessary in order to understand the results and the proofs, provided the reader is willing to accept the material in this section as a black box.

% (for instance certain long exact sequences correspond) Some material in this section gets fairly sophisticated.  We would like to assure the reader that most of these details are entirely unnecessary in order to understand our main results.  It takes a lot of work to introduce the homology theory we care about, but the technique introduced in \cite{ElHog16a} greatly simplifies the task of computing it in many cases.

If $\CC$ is an additive category then we let $\Ch^b(\CC)$ denote the category of complexes over $\CC$, with differentials of degree $+1$, and we let $\KC^b(\CC)$ denote the homotopy category of complexes, which has the same objects of $\Ch^b(\CC)$, but morphisms are chain maps modulo homotopy.  The homological shift functor is denoted $C\mapsto C[1]$, where $C[1]_k = C_{k+1}$ and $d_{C[1]}=-d_C$.

Let $R=R_n$ be the polynomial ring $\Z[x_1,\ldots,x_n]$, graded via $\deg(x_i)=2$, and let $\SBim_n$ denote the category of Soergel bimodules for $S_n$, over $\Z$.  A Soergel bimodule is, in particular, a graded $(R,R)$-bimodule, and morphisms of Soergel bimodules are just bimodule maps of degree zero.  The shift in bimodule degree is denoted $B\mapsto B(1)$, and is such that if $f\in R$ is an element of degree $2k$, then left or right multiplication by $f$ is a morphism of graded bimodules $B\rightarrow B(2k)$ for all $B$.  For each $1\leq i\leq n-1$ we have a Soergel bimodule of the form $B_i=R\otimes_{R^{s_i}}R(1)$, where $R^{s_i}\subset R$ is the subalgebra of $s_i$-invariant polynomials.

The \emph{Rouquier complexes} are complexes of Soergel bimodules associated to braids $\b\in \Br_n$.  They are defined as follows.  Let $\sigma_i$ denote the elementary braid ($1\leq i\leq n-1$), that is, a single positive crossing involving strands $i$ and $i+1$, and define
\[
F(\sigma_i) =  0 \rightarrow B_i\rightarrow R(1) \rightarrow 0  \hskip1in  F(\sigma_i\inv) =  0 \rightarrow R(-1)\rightarrow B_i \rightarrow 0,
\]
where the differential in each is the ``dot'' map, and $B_i$ sits in homological degree zero for both complexes.  In general we define
\[
F(\sigma_{i_1}^{\e_1}\cdots \sigma_{i_r}^{\e_r}):=F(\sigma_{i_1}^{\e_1})\otimes \cdots \otimes F(\sigma_{i_r}^{\e_r}),
\]
for any sequence of indices $i_j\in \{1,2,\ldots,n\}$ and signs $\e_j\in \{\pm 1\}$.  These tensor products depend only on the braid $\b=\sigma_{i_1}^{\e_1}\cdots \sigma_{i_r}^{\e_r}$, and not the expression as a product of generators, up to homotopy, hence $F(\b)$ is well-defined up to isomorphism in $\KC^b(\SBim_n)$.

Khovanov-Rozansky homology can be defined in this context as follows.  Associated to any $(R,R)$-bimodule $B$ we have the Hochschild cohomology groups $\HH^{j}(B):=\Ext_{R^e}^j(R,B)$ in the category of modules over $R^e:=R\otimes_{\Z} R$.  If $B$ is a graded bimodule (for instance a Soergel bimodule), then $\HH^j(B)$ inherits an additional grading: $\HH^j(B)=\bigoplus_{i\in \Z}\HH^{i,j}(B)$.   The functor $\HH^{i,j}$ is $\Z$-linear hence can be extended to complexes.  For any $C\in \KC(\SBim_n)$, we let $\HHH^{i,j,k}(C)$ denote the $k$-th homology group of $\HH^{i,j}(C)$.  Taking the direct sum over all $i,j,k\in \Z$ defines the triply graded abelian group $\HHH(C)$.   In \cite{Kh07}, Khovanov showed that $\HHH(F(\b))$ is isomorphic to the triply graded Khovanov-Rozansky homology \cite{KR08b} of the oriented link $\hat{\b}$ up to isomorphism and overall degree shift\footnote{For the version over the integers, see \cite{Kras10}}.

\begin{definition}
For each $C\in \KC^b(\SBim_n)$, let $\PC_C(Q,A,T)$ denote the Poincar\'e series of $\HHH(C)$:
\[
\PC_C(Q,A,T): = \sum_{i,j,k}Q^iA^jT^k \dim_\Q (\HHH^{i,j,k}(C)\otimes_\Z \Q)
\]
which we also refer to as the \emph{Khovanov-Rozansky (KR) series} of $C$.  We often express $\PC_C$ in terms of the variables $q=Q^2, a=AQ^{-2}$, and $t=T^2Q^{-2}$, in which case the KR series will be written $\PC_C(q,t,a)$. 
\end{definition}

\begin{remark}
In general $\PC_C(Q,A,T)$ is a polynomial in $T,T\inv,A$, but is a Laurent \emph{series} in $Q$.  Although it is possible to show that $\PC_C(Q,A,T)$  becomes a Laurent polynomial in $Q,A,T$ after multiplying by some number of copies of $1-Q^2$.
\end{remark}

There is a \emph{derived category of Soergel bimodules}, which we define as follows.   As above, let $R^e=R\otimes_\Q R$, and let $D^b(R^e)$ denote the bounded derived category of graded $R^e$-modules.  An object of $D(R^e)$ is naturally bigraded, since there is the bimodule degree and the homological degree.  We denote a shift in this bidegree as $B(i,j)=B(i)[j]$.  Let $\CS_n\subset D^(R^e)$ denote the smallest full subcategory containing $\SBim_n$ and closed under the shifts $(i,j)$, mapping cones, and isomorphisms.  If $A,B$ are Soergel bimodules, then the space of morphisms $A\rightarrow B(i,j)$ in $\CS_n$ is naturally identified with the degree $i$ piece of $\Ext^j_{R^e}(A,B)$.  In particular, $\SBim_n$ with its degree zero morphisms includes as a full subcategory of $\CS_n$.  This inclusion is monoidal since the derived tensor product of Soergel bimodules agrees with the usual tensor product (this is because Soergel bimodules are projective as left or right $R$-modules).

The category $\CS_n$ is $\Q$-linear, hence we may consider the category of complexes $\KC^b(\CS_n)$, which includes $\KC^b(\SBim_n)$ as a full subcategory.  The homological shift in $\KC^b(\CS_n)$ is denoted by $[1]$ (not to be confused with the homological shift in $\CS_n$, which is denoted $(0,1)$).   Thus, all said, an object of $\KC^b(\CS_n)$ is naturally triply graded.  A shift in the tridegree is denoted by $C\mapsto C(i,j)[k]$.

\begin{definition}
Given two complexes $C_1,C_2\in \KC^b(\SBim_n)$ (or more generally $\KC^b(\CS_n)$), we let $\DHom(C_1,C_2)$ denote the triply graded space of homs
\[
\DHom(C_1,C_2) = \bigoplus_{i,j,k\in\Z} \Hom_{\CS_n}(C_1,C_2(i,j)[k]).
\]
We also write $\DHom(C,C)$ as $\DEnd(C)$.
\end{definition}
With this definition in place we have $\HHH(C)=\DHom(R,C)$.

We will also regard the variables $Q,A,T$ and $q,a,t$ as grading shifts via
\vskip5pt
\begin{minipage}{2.6in}
\begin{enumerate}
\item $Q C=C(-1,0)[0]$.
\item $A C =(0,-1)[0]$.
\item $T C = (0,0)[-1]$.
\end{enumerate}
\end{minipage} 
\begin{minipage}{2.6in}
\begin{enumerate}
\item $q C=C(-2,0)[0]$.
\item $a C =(2,-1)[0]$.
\item $t C = (2,0)[-2]$.
\end{enumerate}
\end{minipage}
\vskip5pt
Similarly, if $\Gamma$ is a Laurent polynomial in these variables with positive integer coefficients, then $\Gamma C$ denotes the corresponding direct sum of shifted copies of $C$.

\begin{definition}\label{def:externalProd}
Let $\sqcup:\SBim_k\times \SBim_{n-k}\rightarrow \SBim_n$ denote the external tensor product $B\sqcup B':=B\otimes_\Q B'$, and let $\one_k=R_k$ equal the trivial bimodule.  We also denote by $\sqcup$ the induced functor on complexes $\KC^b(\CS_k)\times \KC^b(\CS_{n-k})\rightarrow \KC^b(\CS_n)$.
\end{definition}
\begin{definition}\label{def:simDef}
If $C\in \KC^b(\SBim_m)$ and $D\in \KC^b(\SBim_n)$, then we write $C\sim D$ if $\HHH(C)\cong \HHH(D)$.
\end{definition}

Now we state the equivalences corresponding to Markov moves in our conventions.

\begin{proposition}\label{prop:markov}
For any complexes $C_1,C_2 \in \KC^b(\CS_{n-1})$ we have
\begin{equation}\label{eq:markov+}
(\one_1\sqcup C_1)\otimes F(\sigma_1)\otimes (\one_1\sqcup C_2) \sim t^{\frac{1}{2}} C_1\otimes C_2
\end{equation}
and
\begin{equation}\label{eq:markov-}
(\one_1\sqcup C_1)\otimes F(\sigma_1\inv)\otimes (\one_1\sqcup C_2) \sim q^{-\frac{1}{2}}a C_1\otimes C_2
\end{equation}
\end{proposition}
In this paper we almost never have any use for (\ref{eq:markov-}).

\begin{definition}\label{def:braidExp}
Let $e$ denote the group homomorphism from braid group $\Br_n$ to the additive group of integers sending $\sigma_i^{\pm 1}\mapsto \pm 1$.  We call $e(\b)$ the \emph{braid exponent}.
\end{definition}
In other words, $e(\b)$ is the signed number of crossings, or \emph{writhe}.

\begin{remark}\label{rmk:conventionDiff}
In this paper we use a different conventions than in \cite{Hog15}.  The shift $(i,j)[k]$ here would be denoted $(-i,-j)\ip{-k}$ in \emph{loc.~cit.~}.  Furthermore, the Rouquier complex $F(\b)$ here is obtained by applying the shift $(e(\b),0)[-e(\b)]=t^{e(\b)/2}$ (the shift is in our conventions) to the Rouquier complex $F(\b)$ in \emph{loc.~cit.}.   The Markov moves above are proven in Proposition 3.10 in \emph{loc.~cit.}
\end{remark}

%================================
\subsection{Jones-Wenzl idempotents and the technique}
\label{subsec:JWproj}
%================================

In this section we summarize the technique we will use for all of our computations.  For reference, consult \cite{Hog15} and \cite{ElHog16a}.

In \cite{Hog15}, we constructed complexes $\idemp_n\in \KC^-(\SBim_n)$ which categorify the Young symmetrizers (i.e.~Young idempotents associated to one-row partitions) and proved that
\[
\DEnd(\idemp_n)\cong \Z[u_1,\ldots,u_n]\otimes_\Q \L[\xi_1,\ldots,\xi_n],
\]
where $\L$ denotes exterior algebra.  The tridegree of the even variables is $\deg(u_i)=qt^{1-i}$, and the tridegree of the odd variables is $\deg(\xi_i)=at^{1-i}$.  Let $\QB_n$ denote $t^{n-1}\Cone(u_n)$.  The complexes $\PB_n$ are idempotent, and $\PB_n\otimes \QB_n\simeq \QB_n\simeq \QB_n\otimes \PB_n$.

\begin{remark}
All of the relevant properties of $\QB_n$ are proven in \S 4.7 of \cite{Hog15}.  Remark that the complex $Q_n$ in \cite{Hog15} would be denoted by $t^{1-n}\QB_n$ in our notation.
\end{remark}

\begin{definition}\label{def:K}
 Let $\KB_n= (\QB_1\sqcup \one_{n-1})\otimes (\QB_2\sqcup \one_{n-2}) \otimes \cdots \otimes \QB_n$.
\end{definition}

\begin{remark}
We warn the reader that what we denote $\KB_n$ would have been denoted $\hat{K}_n$ in \cite{ElHog16a}, where $K_n$ denotes instead the tensor product in which we omit the factor involving  $\QB_1$.  One thinks of $\hat{K}_n$ as a \emph{reduced} version of $K_n$.  In this paper we work with the reduced version exclusively.
\end{remark}

\begin{definition}\label{def:symmetries}
The category $\SBim_n$ has the symmetries of a rectangle.  There is a functor $\tau$ which sends $B_i\leftrightarrow B_{n-i}$ and $\tau(B\otimes C)\cong \tau(B)\otimes \tau(C)$, and a functor $\omega$ which sends $B_i\mapsto B_i$ and satisfies $\omega(B\otimes C)\cong \omega(C)\otimes \omega(B)$.
\end{definition}

\begin{proposition}\label{prop:Qsymmetry}
The complexes $\QB_n$ are preserved by $\tau$ and $\omega$, up to homotopy equivalence, and
\[
(\QB_k\sqcup \one_{n-k})\otimes \QB_n \ \ \simeq \ \ (\one_{n-k}\sqcup \QB_k)\otimes \QB_n  \ \ \simeq \ \   \QB_n \otimes (\QB_k\sqcup \one_{n-k})\ \ \simeq \ \  \QB_n \otimes  (\one_{n-k}\sqcup \QB_k)
\]
\end{proposition}
This appears as Proposition 4.39 in \cite{Hog15}.  Thus, the precise form of the tensor product which defines $\KB_n$ is entirely irrelevant up to homotopy equivalence.

The complexes $\QB_n$ satisfy the following recursion:
\begin{equation} \label{eq:Qrecursion}
\begin{minipage}{.8in}
\labellist
\small
\pinlabel $\PB_{n-1}$ at 28 20
\endlabellist
\includegraphics[scale=1]{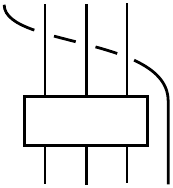}
\end{minipage}
\ \ \simeq \ \ \left(
\begin{minipage}{.8in}
\labellist
\small
\pinlabel $\QB_{n}$ at 27 20
\endlabellist
\includegraphics[scale=1]{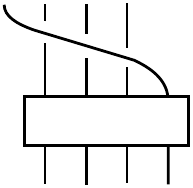}
\end{minipage}
 \ \ \longrightarrow  \ \ 
q\ \begin{minipage}{.8in}
\labellist
\small
\pinlabel $\PB_{n-1}$ at 28 20
\endlabellist
\includegraphics[scale=1]{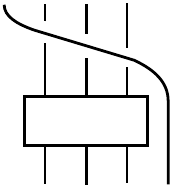}
\end{minipage}\right),
\end{equation}
where the notation $A\simeq (B\rightarrow C)$ means that there is a chain map $\d:B[-1]\rightarrow C$ such that $A\simeq  \Cone(\d)$ which, in turn, means that there is distinguished triangle
\[
C\rightarrow A\rightarrow B\rightarrow C[1].
\]
This is proven in Proposition 4.40 in \cite{Hog15}, modulo conventions.

We work exclusively with the complexes $\KB_n$, which satisfy the following recursion:
\begin{equation} \label{eq:Krecursion}
\begin{minipage}{.8in}
\labellist
\small
\pinlabel $\KB_{n-1}$ at 28 20
\endlabellist
\includegraphics[scale=1]{diagrams/symRecursion1}
\end{minipage}
\ \ \simeq \ \ \left(
\begin{minipage}{.8in}
\labellist
\small
\pinlabel $\KB_{n}$ at 27 20
\endlabellist
\includegraphics[scale=1]{diagrams/symRecursion2}
\end{minipage}
 \ \ \longrightarrow  \ \ 
q\ \begin{minipage}{.8in}
\labellist
\small
\pinlabel $\KB_{n-1}$ at 28 20
\endlabellist
\includegraphics[scale=1]{diagrams/symRecursion3}
\end{minipage}\right),
\end{equation}

\begin{proposition}\label{prop:rouquierAbsorbing}
For any braid $\b$ we have $F(\b)\otimes \KB_n \simeq t^{e(\b)/e} \KB_n \simeq \KB_n\otimes F(\b)$.
\end{proposition}
This is proven in Proposition 4.40 in \cite{Hog15}, modulo conventions.

A direct application of the The Markov move yields the following.

\begin{proposition}\label{prop:Ktrace}
For any complexes $C_1,C_2 \in \KC^b(\CS_{n-1})$ we have
\begin{equation}\label{eq:markovK}
(\one_1\sqcup C_1)\otimes (\KB_k\sqcup \one_{n-k}) \sim (t^{k-1}+a) \cdot C_1\otimes(\KB_{k-1}\sqcup \one_{n-k-1})\otimes C_2
\end{equation}
\end{proposition}
Compare with Proposition 4.4 in \cite{ElHog16a}.

\begin{proposition}\label{prop:K1superpoly}
Let $C=F(\b)$ be the Rouquier complex of a braid $\b$, and let $\hat{C}$ denote $(\KB_1\sqcup \one_{n-1})\otimes C$.  Then
\[
\HHH(C)\cong \Z[x_1]\otimes_\Z \HHH(\hat{C}).
\]
where $x_1$ has degree $q$.
In particular
\[
\PC_C(q,t,a)=\frac{1}{1-q}\PC_{\hat{C}}(q,t,a).
\]
\end{proposition}
This is Proposition 4.12 in \cite{ElHog16a}.

%%%%%%%%%%%%%%%%%%%%%%%%%%%%%%%
\section{Some special complexes of Soergel bimodules and their homology}
\label{sec:computations}
%%%%%%%%%%%%%%%%%%%%%%%%%%%%%%%
In this section we introduce the complexes of Soergel bimodules which interest us and show how to compute their homologies.

%=====================
\subsection{The complexes}
\label{subsec:theComplexes}
%=====================
We will denote objects of $\KC^b(\SBim_n)$ as diagrams drawn in a rectangle, with $n$ boundary points on the top of the rectangle, and $n$ points on the bottom, as in
\[
C  \ \ = \ \ \begin{minipage}{.4in}
\labellist
\small
\pinlabel{$C$} at 17 20
\pinlabel{$\underbrace{ \ \ \ \ \ \ \ \ \ }_n$} at 17 -10
\endlabellist
\includegraphics[scale=.8]{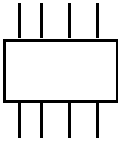}
\end{minipage}\hskip.5in (C\in \KC^b(\SBim_n)).
\]\vskip10pt
In this graphical notation, $F(\b)$ will be denoted simply by the corresponding braid diagram.  The functors $\sqcup$ and $\otimes$ will be drawn as horizontal juxtaposition and vertical concatenation, respectively.
\[
C \sqcup D \ \ = \ \ \begin{minipage}{.4in}
\labellist
\small
\pinlabel{$C$} at 17 20
%\pinlabel{$\underbrace{ \ \ \ \ \ \ \ \ \ }_n$} at 17 -10
\endlabellist
\includegraphics[scale=.8]{diagrams/box}
\end{minipage}\ \ \begin{minipage}{.4in}
\labellist
\small
\pinlabel{$D$} at 17 20
%\pinlabel{$\underbrace{ \ \ \ \ \ \ \ \ \ }_n$} at 14 -10
\endlabellist
\includegraphics[scale=.8]{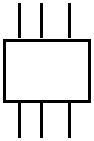}
\end{minipage}, \hskip1in
C_1 \otimes C_2 \ \ = \ \ \begin{minipage}{.4in}
\labellist
\small
\pinlabel{$C_1$} at 17 48
\pinlabel{$C_2$} at 17 20
%\pinlabel{$\underbrace{ \ \ \ \ \ \ \ \ \ }_n$} at 17 -10
\endlabellist
\includegraphics[scale=.8]{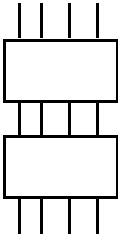}
\end{minipage}.
\]

\begin{definition}\label{def:coxAndTwists}
Let $X_n=\sigma_{n-1}\cdots\sigma_2\sigma_1$, let $\HT_n=\sigma_1(\sigma_2\sigma_1)\cdots (\sigma_{n-1}\cdots\sigma_2\sigma_1)$ denote the half-twist (positive braid lift of the longest element $w_0\in S_n$), and let $\FT_n=\HT_n^2=X_n^n$ denote the full twist.  We use the same notation to refer to the Rouquier complexes of these braids.%  Let $\XB_n$ denote $t^{(1-n)/2}X_n$.  When the index $n$ is understood, we will omit it.
\end{definition}

 We introduce a graphical short-hand for full-twists, since they appear so often in this paper.
\begin{notation}
The full twist will be denoted by a labelled ring circling some number of strands, for instance:
\[
\begin{minipage}{.55in}
\labellist
\tiny
\pinlabel $(1)$ at 50 25
\endlabellist
\includegraphics[scale=.8]{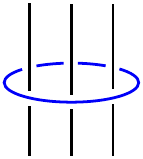}
\end{minipage} \ \ \ = \ \ \ 
\begin{minipage}{.4in}
\labellist
\small
\endlabellist
\includegraphics[scale=.8]{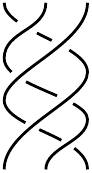}
\end{minipage}.
\]
Topologically such pictures have a precise interpretation in terms of surgery \cite{GompStip}.  The main reason we adopt this notation will be to conserve space in our diagrams.  

A strand labelled `$k$' will be used to denote $k$ parallel copies of that strand.   Then the full twists (and the above notation) can also be interpreted recursively as $\FT_1=\one_1$, the only braid on 1 strand, together with
\begin{equation}\label{eq:FTrecursion}
\begin{minipage}{.4in}
\labellist
\small
\pinlabel{$1$} at 5 -7
\pinlabel{$k$} at 22 -7
\endlabellist
\includegraphics[scale=.8]{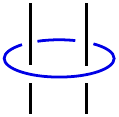}
\end{minipage}
\ \ = \ \ 
\begin{minipage}{.4in}
\labellist
\small
\pinlabel{$1$} at 5 -7
\pinlabel{$k$} at 17 -7
\endlabellist
\includegraphics[scale=.8]{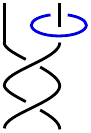}
\end{minipage}
\end{equation}\vskip7pt
\end{notation}

\begin{definition}\label{def:shuffles}
Fix an integer $m\geq 1$, and a sequence $v\in \{0,1,\ldots,m\}^n$.   The (positive) shuffle braid associated to $v$ is the braid $\g_v$ defined as follows.  Let for each $k=0,1,\ldots,m$, let $r_j$ denote the number of indices $i$ with $v_i=k$. Let $\pi_v\in S_n$ denote the unique minimal length permutation sending $(0^{r_0},1^{r_1},\ldots,m^{r_m})$ to $v$.  Note that $\pi_v$ is a minimal length representative of a coset of $S_n/S_{r_1}\times \cdots \times S_{r_m}$.

Let $\g_v$ be the positive braid lift of $\pi_v$.  Let $\omega\g_v$ be the positive braid lift of $\pi_v\inv$
\end{definition}

\begin{example}
If $v=0110201021$, then
\vskip10pt
\[
\g_v \ \ = \ \ \begin{minipage}{1.8in}
\labellist
\small
\pinlabel 0 at 0 -4
\pinlabel 0 at 12 -4
\pinlabel 0 at 24 -4
\pinlabel 0 at 36 -4
\pinlabel 1 at 48 -4
\pinlabel 1 at 62 -4
\pinlabel 1 at 74 -4
\pinlabel 1 at 86 -4
\pinlabel 2 at 100 -4
\pinlabel 2 at 112 -4
\pinlabel 0 at 0 50
\pinlabel 1 at 12  50
\pinlabel 1 at 24 50
\pinlabel 0 at 36 50
\pinlabel 2 at 48 50
\pinlabel 0 at 62 50
\pinlabel 1 at 74 50
\pinlabel 0 at 86 50
\pinlabel 2 at 100 50
\pinlabel 1 at 112 50
\endlabellist
\includegraphics[scale=1]{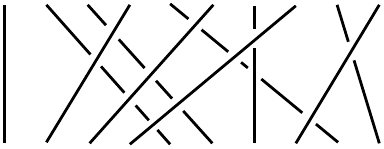}
\end{minipage} \hskip.6in 
\omega \g_v \ \  = \ \ \begin{minipage}{1.8in}
\labellist
\small
\pinlabel 0 at 0 50
\pinlabel 0 at 12 50
\pinlabel 0 at 24 50
\pinlabel 0 at 36 50
\pinlabel 1 at 48 50
\pinlabel 1 at 62 50
\pinlabel 1 at 74 50
\pinlabel 1 at 86 50
\pinlabel 2 at 100 50
\pinlabel 2 at 112 50
\pinlabel 0 at 0 -4
\pinlabel 1 at 12  -4
\pinlabel 1 at 24 -4
\pinlabel 0 at 36 -4
\pinlabel 2 at 48 -4
\pinlabel 0 at 62 -4
\pinlabel 1 at 74 -4
\pinlabel 0 at 86 -4
\pinlabel 2 at 100 -4
\pinlabel 1 at 112 -4
\endlabellist
\includegraphics[scale=1]{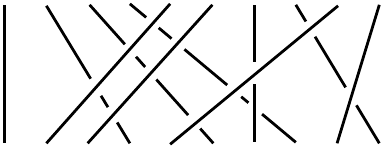}
\end{minipage}
\]\vskip7pt
\end{example}

\begin{definition}\label{def:Cv}
Fix integers $m,n\geq 1$.  For each $v\in \{0,1,\ldots,m\}^n$, we define
\[
\CB_{v}^{(m)} \ \ \ := \ \ \ \begin{minipage}{1.5in}
\labellist
\small
\pinlabel $\KB_{r_1+\cdots+r_{m-1}}$ at 45 140
\pinlabel $\gamma_v$ at 55 174
\pinlabel $\omega\gamma_v$ at 55 30
\tiny
\pinlabel $(1)$ at 125 57
\pinlabel $(1)$ at 125 75
\pinlabel $\vdots$ at 125 93
\pinlabel $(1)$ at 125 110
\pinlabel $(1)$ at 125 126
\pinlabel $r_0$ at 5 -5
\pinlabel $r_1$ at 25 -5
\pinlabel $\cdots$ at 45 -5
\pinlabel $\cdots$ at 65 -5
\pinlabel $r_{m-1}$ at 85 -5
\pinlabel $r_m$ at 105 -5
\endlabellist
\includegraphics[scale=.8]{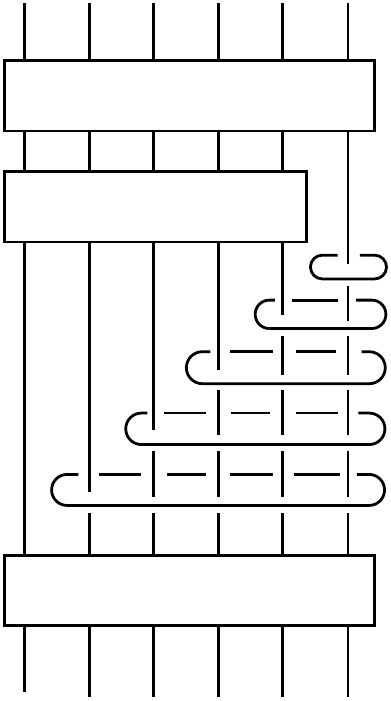}
\end{minipage},
\]\vskip7pt
\noindent
where $r_j$ is the number of indices $i$ such that $v_i=j$.   More precisely, $\CB_{v}^{(m)}$ is the following tensor product (omitting the explicit symbols $\otimes$ for reasons of space).
\[
 \g_v  (\KB_{r_0+\cdots +r_{m-1}}\sqcup \FT_{r_m}) (\one_{r_0+\cdots +r_{m-2}}\sqcup \FT_{r_{m-1}+r_{m}}) \cdots  (\one_{r_0}\sqcup \FT_{r_1+\cdots+r_m}) (\omega\g_v).
\]
\end{definition}
\begin{remark}
Let $v\in \{0,1,\ldots,m\}^n$ be given, and pick an index $i\in\{1,\ldots,n\}$.  If $v_i=k$, then we interpret $k$ as the number of full twists that the the $i$-th strand passes through in its journey from the bottom of the diagram to the top.  If $k<m$ then this strand eventually meets the symmetrizer $\KB$.  If $k=m$, then this strand passes through, unhindered, from the bottom boundary to the top boundary.
\end{remark}

Observe that
\begin{equation}\label{eq:FTtimesC}
\FT \CB^{(m)}_{v_1,\ldots,v_n} = \CB^{(m+1)}_{v_1+1,\ldots,v_n+1}
\end{equation}
for all $v\in \{0,1,\ldots,m\}^n$, and $\Cox^{nm}=\FT^m=\CB^{(m)}_{m^n}$.  More generally, $\Cox^{nm+r}=\Cox^r \CB_{m^n}^{(m)}$ is the complex which computes the triply graded homology of the $(n,nm+r)$ torus link.

It will be convenient later to work with a shifted version of $X_n$.
\begin{definition}\label{def:shiftedCox}
Let $\XB_n= t^{(1-n)/2}X_n$ denote the the Rouquier complex of $\sigma_{n-1}\cdots \sigma_2\sigma_1$, with a shift.  When the index $n$ is understood, it will be omitted.
\end{definition}

The theme of this paper (building on the main idea of \cite{ElHog16a}) is that the complexes $\XB^r \CB^{(m)}_v$ play a very imporant auxilliary role in the computation of link homology.  We also believe these these building blocks will prove to be important tools for studying other aspects of Soergel category.

%=====================
\subsection{Relations I}
\label{subsec:rels1}
%=====================
Ultimately we wish to compute $\HHH(\XB^r \CB^{(m)}_v)$ for all $r\geq 0$ and all $v\in \{0,1,\ldots,m\}^n$.  Because of (\ref{eq:FTtimesC}), we may assume that $0\leq r\leq n-1$.  Note that $\CB_{m^n}^{(m)}=\XB^{nm}$, hence $\XB^r \CB_{m^n}^{(m)}$ computes the Khovanov-Rozansky homology of the $(n,nm+r)$ torus link, up to shift.  Of course, many other knots and links occur as special cases of the $\XB^r\CB_v^{(m)}$ as well; these will investigated in the case $r=0,1$ in \S \ref{subsec:specialBraids}.

In this section we give techniques for rewriting $\XB\C_{v}^{(m)}\in \KC^b(\SBim_n)$ up to homotopy equivalence.   Most of the proofs in this section and the next will be proofs by pictures.  They can easily be converted into more rigorous (but less readable) arguments.

%In the next two sections we discuss relations satisfied by the complexes $\CB_{v}$.  The equivalences in this section involve only the grading shifts $q,t$ (and not $a$), and are valid relations in $\KC^b(\SBim_n)$, whereas the relations in \S \ref{subsec:relsII} involve the various Markov moves, hence may involve $a$ as well, and are only true upon taking $\HHH$.   We may refer to these relations of type $\simeq$ or of type $\sim$ (where $\simeq$ is homotopy equivalence, and $\sim$ is the relation from Definition \ref{def:simDef}).  Relations of type $\simeq$ are stronger since for example $C_1\simeq C_2$  implies $D\otimes C_1\simeq D\otimes C_2$, while if $C_1\sim C_2$ then it does not even make sense to form the tensor product $D\otimes C_i$ since $C_1$ and $C_2$ may lie in different categories. Even if $C_1$ and $C_2$ are both in the same $\KC^b(\CS_n)$, it is not necessarily true that $C_1\sim C_2 $ implies $D\otimes C_1\sim D\otimes C_2$.  

\begin{proposition}\label{prop:Fmv}
For $v\in \{0,1,\ldots,m\}^n$ we have
\begin{equation}\label{eq:Cmv}
\XB \CB_{mv}^{(m)}\simeq \Tot\Big(\CB_{v(m-1)}^{(m)}\XB \ \ \rightarrow \ \ q \CB_{vm}^{(m)}\XB\Big).
\end{equation}
\end{proposition}
We remind the reader that if $A,B,C$ are chain complexes, then the notation $B\simeq (C\rightarrow A)$ indicates the existence of a distinguished triangle
\[
A\rightarrow B\rightarrow C\rightarrow A[1].
\]
\begin{proof}
Observe:
\[
\XB \CB^{(m)}_{mv} \ \ = \ \ \begin{minipage}{1.4in}
\labellist
\small
\pinlabel $\KB$ at 40 105
\pinlabel $\gamma_v$ at 50 137
\pinlabel $\omega \gamma_v$ at 50 12
\tiny
\pinlabel $1$ at -3 -4
\pinlabel $r_0$ at 15 -4
\pinlabel $r_1$ at 28 -4
\pinlabel $r_2$ at 41 -4
\pinlabel $r_3$ at 53 -4
\pinlabel $r_4$ at 67 -4
\pinlabel $r_5$ at 80 -4
\pinlabel $(1)$ at 110 47
\pinlabel $(1)$ at 110 59
\pinlabel $(1)$ at 110 71
\pinlabel $(1)$ at 110 83
\pinlabel $(1)$ at 110 95
\endlabellist
\includegraphics[scale=.8]{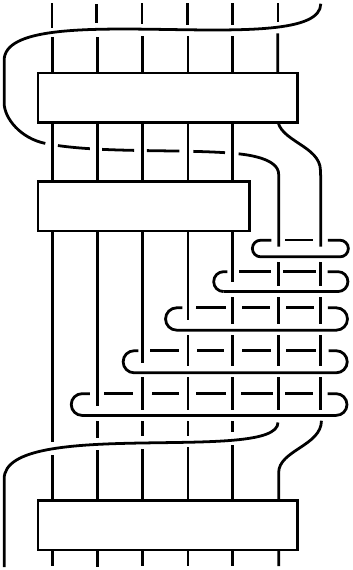}
\end{minipage}.
\]\vskip8pt
We have drawn the case $m=5$.  We remind the reader that the labels below the diagram indicate numbers of strands.  The leftmost strand on the bottom boundary corresponds to the first index of $mv$, which is $m$.  Using (\ref{eq:Krecursion}), the above is homotopy equivalent to
\begin{equation}\label{eq:Dmv_simp1}
\Tot\left(\ \ \ 
\begin{minipage}{1.4in}
\labellist
\small
\pinlabel $\KB$ at 46 111
\pinlabel $\gamma_v$ at 50 140
\pinlabel $\omega \gamma_v$ at 50 17
\tiny
\pinlabel $1$ at -3 -4
\pinlabel $r_0$ at 15 -4
\pinlabel $r_1$ at 28 -4
\pinlabel $r_2$ at 41 -4
\pinlabel $r_3$ at 53 -4
\pinlabel $r_4$ at 67 -4
\pinlabel $r_5$ at 80 -4
\pinlabel $(1)$ at 110 47
\pinlabel $(1)$ at 110 59
\pinlabel $(1)$ at 110 71
\pinlabel $(1)$ at 110 83
\pinlabel $(1)$ at 110 95
\endlabellist
\includegraphics[scale=.8]{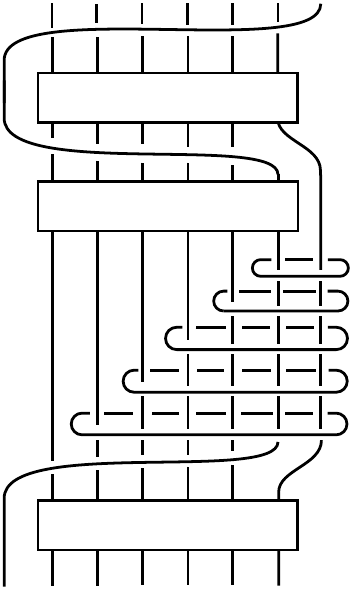}
\end{minipage}
 \ \ \rightarrow \ \ q
 \ \begin{minipage}{1.4in}
\labellist
\small
\pinlabel $\KB$ at 40 105
\pinlabel $\gamma_v$ at 50 137
\pinlabel $\omega \gamma_v$ at 50 12
\tiny
\pinlabel $1$ at -3 -4
\pinlabel $r_0$ at 15 -4
\pinlabel $r_1$ at 28 -4
\pinlabel $r_2$ at 41 -4
\pinlabel $r_3$ at 53 -4
\pinlabel $r_4$ at 67 -4
\pinlabel $r_5$ at 80 -4
\pinlabel $(1)$ at 110 47
\pinlabel $(1)$ at 110 59
\pinlabel $(1)$ at 110 71
\pinlabel $(1)$ at 110 83
\pinlabel $(1)$ at 110 95
\endlabellist
\includegraphics[scale=.8]{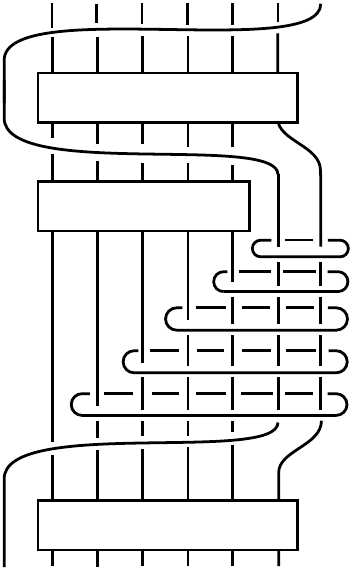}
\end{minipage}
\right)
\end{equation}\vskip8pt
The first term is homotopy equivalent to:
\begin{equation}\label{eq:overlappingBox}
\begin{minipage}{1.4in}
\labellist
\small
\pinlabel $\KB$ at 48 122
\pinlabel $\gamma_v$ at 50 154
\pinlabel $\omega \gamma_v$ at 50 19
\tiny
\pinlabel $1$ at -3 -4
\pinlabel $r_0$ at 15 -4
\pinlabel $r_1$ at 28 -4
\pinlabel $r_2$ at 41 -4
\pinlabel $r_3$ at 53 -4
\pinlabel $r_4$ at 67 -4
\pinlabel $r_5$ at 80 -4
\pinlabel $(1)$ at 110 47
\pinlabel $(1)$ at 110 59
\pinlabel $(1)$ at 110 71
\pinlabel $(1)$ at 110 83
\pinlabel $(1)$ at 110 112
\endlabellist
\includegraphics[scale=.8]{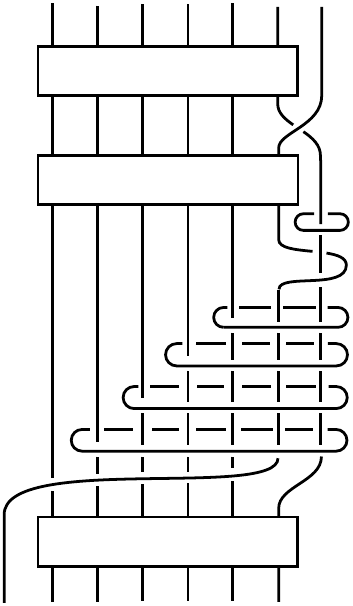}
\end{minipage}
\ \ \simeq \ \ 
\begin{minipage}{1.4in}
\labellist
\small
\pinlabel $\KB$ at 48 121
\pinlabel $\gamma_v$ at 50 152
\pinlabel $\omega \gamma_v$ at 50 17
\tiny
\pinlabel $1$ at -3 -4
\pinlabel $r_0$ at 15 -4
\pinlabel $r_1$ at 28 -4
\pinlabel $r_2$ at 41 -4
\pinlabel $r_3$ at 53 -4
\pinlabel $r_4$ at 67 -4
\pinlabel $r_5$ at 80 -4
\pinlabel $(1)$ at 110 45
\pinlabel $(1)$ at 110 57
\pinlabel $(1)$ at 110 69
\pinlabel $(1)$ at 110 81
\pinlabel $(1)$ at 110 110
\endlabellist
\includegraphics[scale=.8]{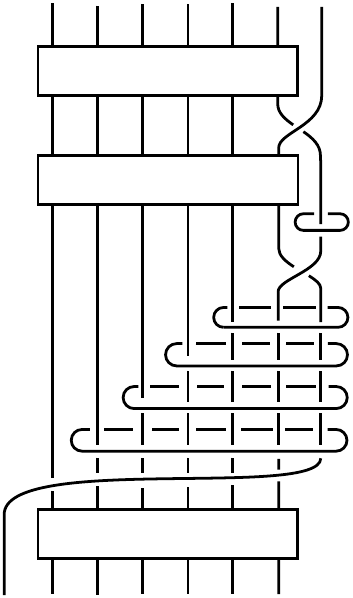}
\end{minipage}
\ \ \simeq \ \
\begin{minipage}{1.4in}
\labellist
\small
\pinlabel $\KB$ at 45 113
\pinlabel $\gamma_v$ at 50 144
\pinlabel $\omega \gamma_v$ at 50 29
\tiny
\pinlabel $1$ at -3 -4
\pinlabel $r_0$ at 15 -4
\pinlabel $r_1$ at 28 -4
\pinlabel $r_2$ at 41 -4
\pinlabel $r_3$ at 53 -4
\pinlabel $r_4$ at 67 -4
\pinlabel $r_5$ at 80 -4
\pinlabel $(1)$ at 110 52
\pinlabel $(1)$ at 110 64
\pinlabel $(1)$ at 110 76
\pinlabel $(1)$ at 110 88
\pinlabel $(1)$ at 110 100
\endlabellist
\includegraphics[scale=.8]{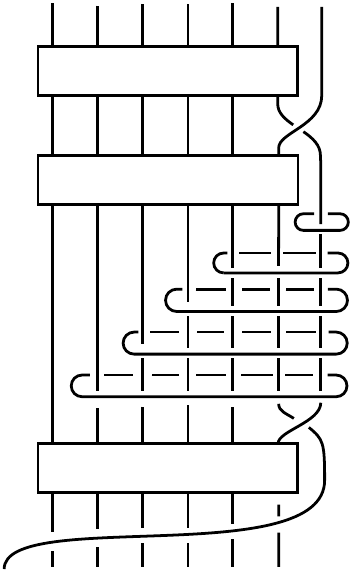}
\end{minipage},
\end{equation}\vskip8pt
which is $\CB_{v(m-1)}^{(m)}\XB$.  The equivalence relating the first term of \ref{eq:Dmv_simp1} with the first term of (\ref{eq:overlappingBox}) is an isotopy (drag the strand over $\gamma_v$) and the recursion for full twists (\ref{eq:FTrecursion}).  The remaining equivalences in (\ref{eq:overlappingBox}) are obvious isotopies. The second term in (\ref{eq:Dmv_simp1}) the right can be simplified using isotopies:
\[
\begin{minipage}{1.4in}
\labellist
\small
\pinlabel $\KB$ at 44 105
\pinlabel $\gamma_v$ at 50 131
\pinlabel $\omega \gamma_v$ at 50 13
\tiny
\pinlabel $1$ at -3 -4
\pinlabel $r_0$ at 15 -4
\pinlabel $r_1$ at 28 -4
\pinlabel $r_2$ at 41 -4
\pinlabel $r_3$ at 53 -4
\pinlabel $r_4$ at 67 -4
\pinlabel $r_5$ at 80 -4
\pinlabel $(1)$ at 110 47
\pinlabel $(1)$ at 110 59
\pinlabel $(1)$ at 110 71
\pinlabel $(1)$ at 110 83
\pinlabel $(1)$ at 110 95
\endlabellist
\includegraphics[scale=.8]{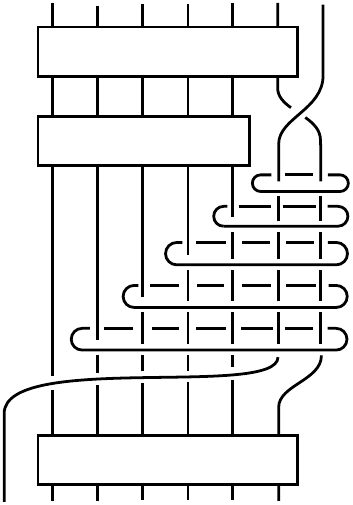}
\end{minipage}
\ \ \simeq \ \ 
\begin{minipage}{1.4in}
\labellist
\small
\pinlabel $\KB$ at 44 105
\pinlabel $\gamma_v$ at 50 131
\pinlabel $\omega \gamma_v$ at 50 13
\tiny
\pinlabel $1$ at -3 -4
\pinlabel $r_0$ at 15 -4
\pinlabel $r_1$ at 28 -4
\pinlabel $r_2$ at 41 -4
\pinlabel $r_3$ at 53 -4
\pinlabel $r_4$ at 67 -4
\pinlabel $r_5$ at 80 -4
\pinlabel $(1)$ at 110 47
\pinlabel $(1)$ at 110 59
\pinlabel $(1)$ at 110 71
\pinlabel $(1)$ at 110 83
\pinlabel $(1)$ at 110 95
\endlabellist
\includegraphics[scale=.8]{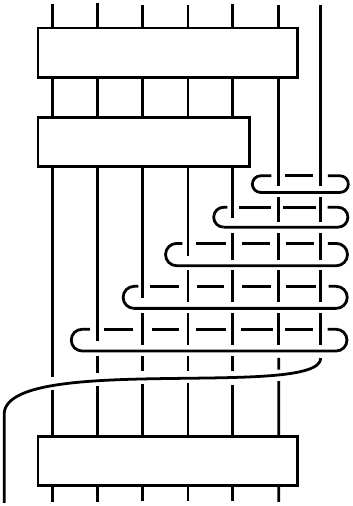}
\end{minipage}
\ \ \simeq \ \
\begin{minipage}{1.4in}
\labellist
\small
\pinlabel $\KB$ at 44 105
\pinlabel $\gamma_v$ at 50 131
\pinlabel $\omega \gamma_v$ at 50 28
\tiny
\pinlabel $1$ at -3 -4
\pinlabel $r_0$ at 15 -4
\pinlabel $r_1$ at 28 -4
\pinlabel $r_2$ at 41 -4
\pinlabel $r_3$ at 53 -4
\pinlabel $r_4$ at 67 -4
\pinlabel $r_5$ at 80 -4
\pinlabel $(1)$ at 110 47
\pinlabel $(1)$ at 110 59
\pinlabel $(1)$ at 110 71
\pinlabel $(1)$ at 110 83
\pinlabel $(1)$ at 110 95
\endlabellist
\includegraphics[scale=.8]{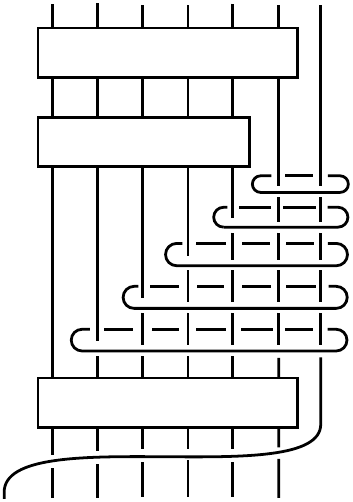}
\end{minipage},
\]\vskip7pt
which is $\CB_{vm}^{(m)} \XB$, as claimed.
\end{proof}

\begin{proposition}\label{prop:Fkv}
For $v\in \{0,1,\ldots,m\}^n$ and $1\leq k\leq m-1$ we have
\begin{equation}\label{eq:Ckv}
\XB  \CB_{k,v}^{(m)}\simeq t^{\#\{i\:|\:v_i<k\}} \CB_{v(k-1)}^{(m)} \XB.
\end{equation}
in $\KC^b(\SBim_n)$.
\end{proposition}
\begin{proof}
Let us omit the grading shifts for now to conserve space.  Observe:
\[
\begin{minipage}{1.2in}
\labellist
\small
\pinlabel $\KB$ at 44 109
\pinlabel $\gamma_v$ at 52 152
\pinlabel $\omega \gamma_v$ at 52 21
\tiny
\pinlabel $1$ at 0 -4
\pinlabel $r_0$ at 15 -4
\pinlabel $r_1$ at 30 -4
\pinlabel $r_2$ at 52 -4
\pinlabel $r_3$ at 66 -4
\pinlabel $r_4$ at 81 -4
\pinlabel $(1)$ at 100 57
\pinlabel $(1)$ at 100 71
\pinlabel $(1)$ at 100 85
\pinlabel $(1)$ at 100 99
\endlabellist
\includegraphics[scale=.80]{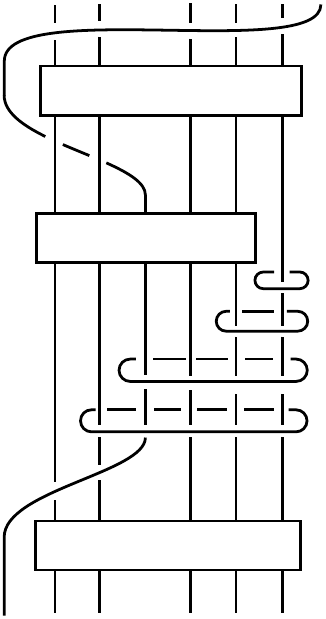}
\end{minipage}
\ \ \simeq \ \ 
\begin{minipage}{1.2in}
\labellist
\small
\pinlabel $\KB$ at 44 133
\pinlabel $\gamma_v$ at 52 175
\pinlabel $\omega \gamma_v$ at 52 21
\tiny
\pinlabel $1$ at 0 -4
\pinlabel $r_0$ at 15 -4
\pinlabel $r_1$ at 30 -4
\pinlabel $r_2$ at 52 -4
\pinlabel $r_3$ at 66 -4
\pinlabel $r_4$ at 81 -4
\pinlabel $(1)$ at 102 60
\pinlabel $(1)$ at 100 95
\pinlabel $(1)$ at 100 108
\pinlabel $(1)$ at 100 121
\endlabellist
\includegraphics[scale=.8]{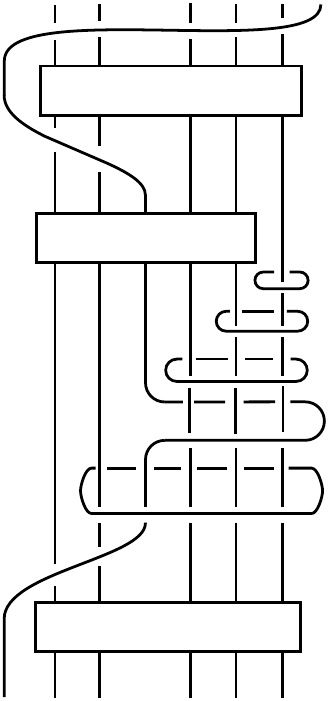}
\end{minipage}
\ \ \simeq \ \
\begin{minipage}{1.2in}
\labellist
\small
\pinlabel $\KB$ at 44 139
\pinlabel $\gamma_v$ at 55 187
\pinlabel $\omega \gamma_v$ at 47 18
\tiny
\pinlabel $1$ at 0 -4
\pinlabel $r_0$ at 15 -4
\pinlabel $r_1$ at 30 -4
\pinlabel $r_2$ at 52 -4
\pinlabel $r_3$ at 66 -4
\pinlabel $r_4$ at 81 -4
\pinlabel $(1)$ at 103 86
\pinlabel $(1)$ at 100 99
\pinlabel $(1)$ at 100 112
\pinlabel $(1)$ at 100 125
\endlabellist
\includegraphics[scale=.8]{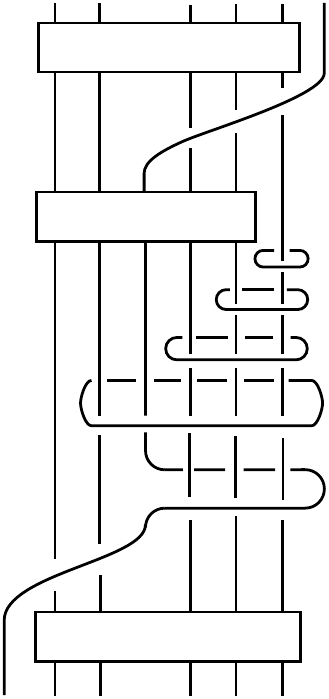}
\end{minipage}
\ \ \simeq \ \
\begin{minipage}{1.2in}
\labellist
\small
\pinlabel $\KB$ at 44 139
\pinlabel $\gamma_v$ at 52 188
\pinlabel $\omega \gamma_v$ at 52 36
\tiny
\pinlabel $1$ at 0 -4
\pinlabel $r_0$ at 15 -4
\pinlabel $r_1$ at 30 -4
\pinlabel $r_2$ at 52 -4
\pinlabel $r_3$ at 66 -4
\pinlabel $r_4$ at 81 -4
\pinlabel $(1)$ at 103 86
\pinlabel $(1)$ at 100 99
\pinlabel $(1)$ at 100 112
\pinlabel $(1)$ at 100 125
\endlabellist
\includegraphics[scale=.8]{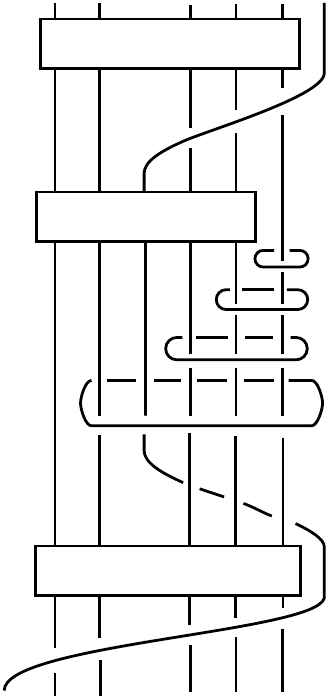}
\end{minipage}.
\]\vskip7pt
We have drawn the case $m=4$, $k=2$.  In the first equivalence we expanded the $k$-th full twist from the bottom using the full twist recursion (\ref{eq:FTrecursion}).  We also switched several negative crossings to positive crossings.  This contributes a grading shift of $t$ for each crossing, by Proposition \ref{prop:rouquierAbsorbing}.  The number of switched crossings in the above diagram is $r_0+r_1$, and in general will be $r_0+\cdots+r_{k-1}$.  This is equal to the number of indices $i$ such that $v_i<k$.  The other equivalences are all isotopies.  This last complex is $\CB_{v(k-1)}^{(m)}\XB$.  Taking into account the aforementioned grading shift proves the proposition.
\end{proof}

Iterating Propositions \ref{prop:Fkv} and \ref{prop:Fmv} we can express any complex $\CB^{(m)}_v$ as an iterated mapping cone involving shifts of conjugates of complexes $\CB^{(m)}_{0w}$ in which the first coordinate is zero.

To be more precise, let $C_i$ be chain complexes, indexed by a finite set $I$.  We say that $D\simeq \bigoplus_i C_i$ with \emph{twisted differential} if there is some partial order on $I$ with respect to which the differential of $D$ is lower triangular.  That is, a \emph{twisted differential} on $\bigoplus_{i\in I}C_i$ differs from the direct sum differential by an endomorphism which is strictly lower triangular for some partial order.

If $C$ is a chain complex, then let $[C]$ denote the set of chain complexes $D$ with $D\simeq F(\b) \otimes C\otimes F(\b)\inv$.  By abuse of notation, we write
\[
E\in \left(\bigoplus_{i\in I} [D_i], d\right)
\]
whenever there exist complexes $C_i\in [D_i]$ such that
\[
E \simeq \left(\bigoplus_{i\in I} C_i, d\right),
\]
with twisted differential.  Note that $\HHH(C_i)\cong \HHH(D_i)$; in particular, if each $\HHH(D_i)$ is even, then so is $\HHH(C)$, and $\HHH(C)\cong \bigoplus_i \HHH(D_i)$.

Now, let $v\in \{1,\ldots,m\}^{n_1}$ and $w\in \{0,1,\ldots,m\}^{n_2}$ be given.  Choose a subset $S\in v\inv(m)$, and let $v(S)\in \{0,1,\ldots,m\}^{n_1}$ be the sequence $v(S)_i=v_i=m$ if $i\in S$, and $v(S)_i=v_i -1$ otherwise.  Also, let $c(S)$ denote the number
\[
\#\{i,k\:|\: v_i\leq w_k\leq m-1 \} + \#\{i<j\:|\: v_i\leq v_j\leq m-1 \}+\#\{j<i\:|\: v_i\leq v_j\leq m-1 \},
\]
where in the above sets $i,j$ range over all indices $\{1,\ldots,n_1\}\setminus S$, and $k$ ranges over $\{1,\ldots,n_2\}$.

\begin{proposition}\label{prop:simultCycling}
If $v\in \{1,\ldots,m\}^{n_1}$ and $w\in \{0,1,\ldots,m\}^{n_2}$, we have
\begin{subequations}
\begin{equation}
\XB^r\CB^{(m)}_{mv} \simeq \Big([\XB^r \CB^{(m)}_{v(m-1)}] \rightarrow q [\XB^r\CB^{(m)}_{vm}] \Big) 
\end{equation}    
\begin{equation}
\XB^r\CB^{(m)}_{kv} \simeq t^{\#\{i\:|\: v_i<k\}} [\XB^r \CB^{(m)}_{v(k-1)}].
\end{equation}
\end{subequations}
\[
\XB^r\CB_{vw}^{(m)} \ \ \in \ \ \left(\bigoplus_S q^{|S|}t^{c(S)}[\XB^r\CB_{wv(S)}^{(m)}] \ , \ \  d\ \right)
\]
with twisted differential where the sum is over subsets $S\subset v\inv(m)$, and, $c(S)$ is defined the comments above.
\end{proposition}
\begin{proof}
This is an easy induction on $n_1$. 
\end{proof}

%=========================
\subsection{Relations II}
\label{subsec:relsII}
%==========================
Iterating Propositions \ref{prop:Fkv} and \ref{prop:Fmv} (see also Proposition \ref{prop:simultCycling}) allows us to to express $\XB^r \CB_{v}^{(m)}$ as an iterated mapping cone constructed from (shifts of conjugates of) complexes of the form $\XB^r\CB_{0w}^{(m)}$.  Thus, for the purposes of computing $\HHH(\XB^r \CB_v^{(m)})$, we are (modulo a spectral sequence argument) reduced to computing $\XB^r\HHH(\CB_{0w}^{(m)})$.

How one proceeds from here depends heavily on $r$.  We will take care of the case $r=\in \{0,1\}$, and we save the case of $r\geq 2$ for later investigations.

In case $r=0$ we have the following (compare with Proposition 4.4 in \cite{ElHog16a}).

\begin{proposition}\label{prop:CvMarkov}
For every $v\in \{0,1,\ldots,m\}^n$, we have
\begin{equation}\label{eq:C0v_rZero}
\CB_{0v}^{(m)} \sim (t^{\#\{i\:|\: v_i<m\}}+a)\CB_{v}^{(m)},
\end{equation}
where $\sim$ is the equivalence relation from Definition \ref{def:simDef}.
\end{proposition}
\begin{proof}
In this case $\CB_{0v}^{(m)}$ can be expressed as $\CB_{0v}^{(m)} = (\one_1\sqcup C_1)\otimes (\KB_{c+1}\sqcup \one_{n-c})\otimes (\one_1\sqcup C_2)$ for some complexes $C_i\in \KC^b(\SBim_n)$ with $C_1C_2\simeq \CB_{v}$.  The statement now follows from Proposition \ref{prop:Ktrace}.
\end{proof}

Our rule for simplifying $\HHH(\XB \CB_{0v}^{(m)})$ depends on whether $v_1=m$ or $v_1<m$. 

\begin{proposition}\label{prop:DvMarkov1}
For every $v\in \{0,1,\ldots,m\}^n$ and every $k\in \{0,1,\ldots,m-1\}$ we have
\begin{equation}\label{eq:C0kv_rOne}
\XB\CB_{0kv}^{(m)} \sim (t^{\#\{i\:|\: v_i<m\}+1}+a)\XB\CB_{kv}^{(m)}.
\end{equation}
\end{proposition}
\begin{proof}
Let $c:=\#\{i\:|\:v_k<m\}+1$.  Compute:
\[
 \begin{minipage}{1.55in}
\labellist
\small
\pinlabel $\gamma_v$ at 72 178
\pinlabel $\KB_c$ at 49 131
\pinlabel $\omega \gamma_v$ at 72 18
\tiny
\pinlabel $1$ at 5 -4
\pinlabel $1$ at 20 -4
\pinlabel $r_0$ at 36 -4
\pinlabel $r_1$ at 48 -4
\pinlabel $r_2$ at 71 -4
\pinlabel $r_3$ at 86 -4
\pinlabel $r_4$ at 101 -4
\pinlabel $r_5$ at 116 -4
\pinlabel $(1)$ at 130 61
\pinlabel $(1)$ at 130 75
\pinlabel $(1)$ at 130 89
\pinlabel $(1)$ at 130 103
\pinlabel $(1)$ at 130 117
\endlabellist
\includegraphics[scale=.8]{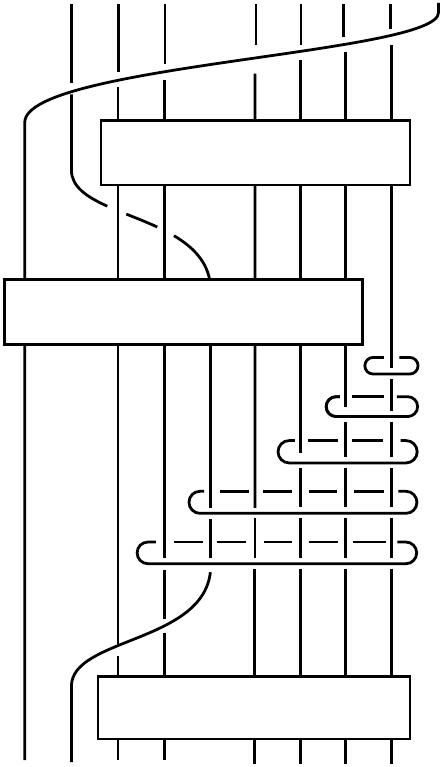}
\end{minipage}
\ \ \simeq \ \ 
t^{\frac{1}{2}}\begin{minipage}{1.55in}
\labellist
\pinlabel $\gamma_v$ at 72 178
\pinlabel $\KB_c$ at 49 131
\pinlabel $\omega \gamma_v$ at 72 18
\tiny
\pinlabel $1$ at 5 -4
\pinlabel $1$ at 20 -4
\pinlabel $r_0$ at 36 -4
\pinlabel $r_1$ at 48 -4
\pinlabel $r_2$ at 71 -4
\pinlabel $r_3$ at 86 -4
\pinlabel $r_4$ at 101 -4
\pinlabel $r_5$ at 116 -4
\pinlabel $(1)$ at 130 61
\pinlabel $(1)$ at 130 75
\pinlabel $(1)$ at 130 89
\pinlabel $(1)$ at 130 103
\pinlabel $(1)$ at 130 117
\small
\endlabellist
\includegraphics[scale=.8]{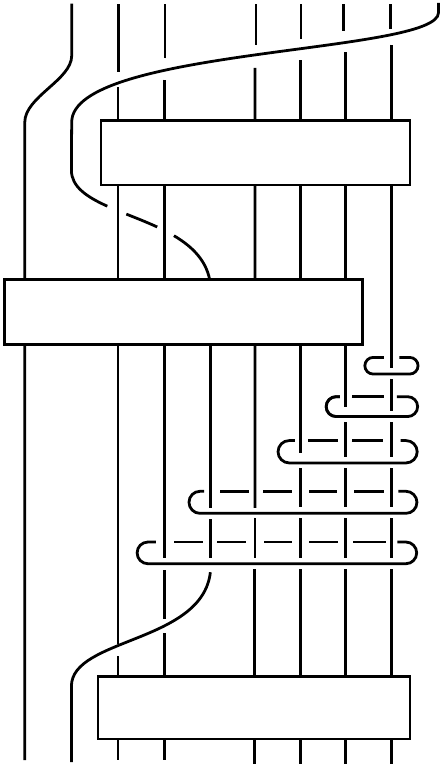}
\end{minipage}
\ \ \sim \ \ 
 \ t^{\frac{1}{2}}(t^c+a)\begin{minipage}{1.55in}
\labellist
\small
\pinlabel $\gamma_v$ at 57 178
\pinlabel $\KB_{c-1}$ at 44 131
\pinlabel $\omega \gamma_v$ at 57 18
\tiny
\pinlabel $1$ at 1 -4
\pinlabel $r_0$ at 17 -4
\pinlabel $r_1$ at 29 -4
\pinlabel $r_2$ at 54 -4
\pinlabel $r_3$ at 69 -4
\pinlabel $r_4$ at 84 -4
\pinlabel $r_5$ at 97 -4
\pinlabel $(1)$ at 115 61
\pinlabel $(1)$ at 115 75
\pinlabel $(1)$ at 115 89
\pinlabel $(1)$ at 115 103
\pinlabel $(1)$ at 115 117
\endlabellist
\includegraphics[scale=.8]{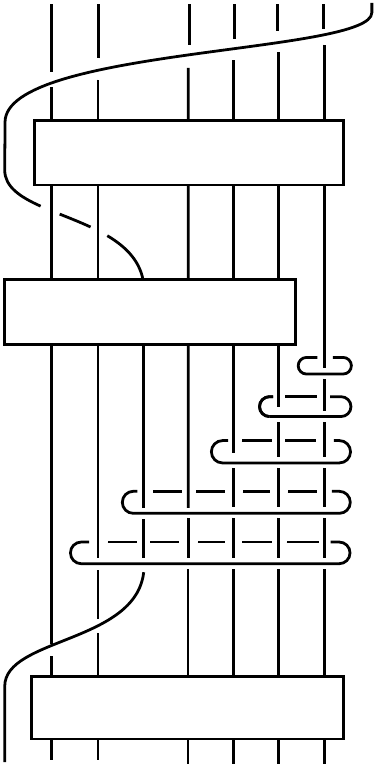}
\end{minipage}.
\]\vskip7pt
%t^{-\frac{1}{2}}(t^{c+1}+a)
We have pictured the case $m=5,k=2$.  In the first equivalence, we absorbed the leftmost crossing on the top using Proposition \ref{prop:rouquierAbsorbing}.  This contributes a grading shift of $t^{1/2}$.  This shift is cancelled by our choice of normalization for $\XB_n$, since $\XB_{n+2} = t^{(-1-n)/2}X_{n+2}$ and $\XB_{n+1}=t^{-n/2}X_{n+1}$.  In the second equivalence we used Proposition \ref{prop:Ktrace}.  This completes the proof.
\end{proof}

The final relation in this section describes how to simplify $\HHH(\CB_{0mv}^{(m)})$.

\begin{proposition}\label{prop:DvMarkov2}
For all $v\{0,1,\ldots,m\}^n$ we have
\begin{equation}\label{eq:C0mv_rOne}
\XB\CB_{0mv}^{(m)} \sim t^{\#\{i\:|\:v_i<m\}} \XB\CB_{v(m-1)}^{(m)}.
\end{equation}
\end{proposition}
\begin{proof}
Let us temporarily ignore the degree shifts.  We will discuss shifts at the end.
\begin{equation}
\XB_n\CB_{0mv}^{(m)} \ \ \ \simeq \ \ \ 
\begin{minipage}{1.4in}
\labellist
\small
\pinlabel $\gamma_v$ at 65 131
\pinlabel $\KB$ at 44 100
\pinlabel $\omega \gamma_v$ at 65 14
\tiny
\pinlabel $1$ at 8 -4
\pinlabel $1$ at 24 -4
\pinlabel $r_0$ at 38 -4
\pinlabel $r_1$ at 51 -4
\pinlabel $r_2$ at 65 -4
\pinlabel $r_3$ at 79 -4
\pinlabel $r_4$ at 102 -4
\pinlabel $(1)$ at 118 49
\pinlabel $(1)$ at 118 63
\pinlabel $(1)$ at 118 76
\pinlabel $(1)$ at 118 89
\endlabellist
\includegraphics[scale=.8]{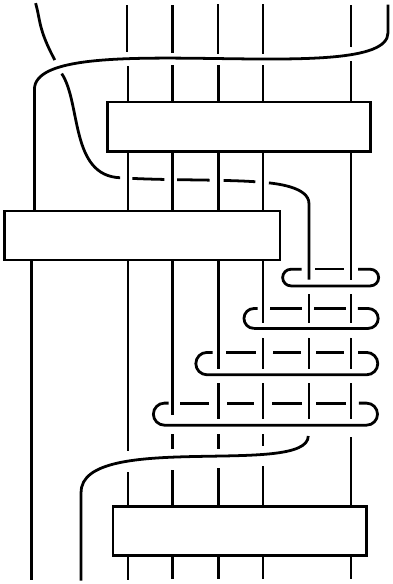}
\end{minipage}
\ \ \ \simeq \ \ \ 
\begin{minipage}{1.4in}
\labellist
\small
\pinlabel $\gamma_v$ at 70 166
\pinlabel $\KB$ at 44 101
\pinlabel $\omega \gamma_v$ at 67 16
\tiny
\pinlabel $1$ at 8 -4
\pinlabel $1$ at 24 -4
\pinlabel $r_0$ at 38 -4
\pinlabel $r_1$ at 51 -4
\pinlabel $r_2$ at 65 -4
\pinlabel $r_3$ at 79 -4
\pinlabel $r_4$ at 102 -4
\pinlabel $(1)$ at 118 49
\pinlabel $(1)$ at 118 63
\pinlabel $(1)$ at 118 76
\pinlabel $(1)$ at 118 89
\endlabellist
\includegraphics[scale=.8]{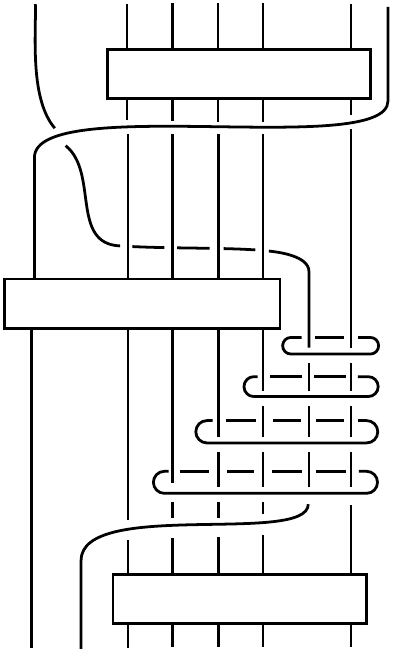}
\end{minipage}
\ \ \ \simeq \ \ \ 
\begin{minipage}{1.4in}
\labellist
\small
\pinlabel $\gamma_v$ at 60 166
\pinlabel $\KB$ at 44 127
\pinlabel $\omega \gamma_v$ at 62 15
\tiny
\pinlabel $1$ at 1 -4
\pinlabel $1$ at 17 -4
\pinlabel $r_0$ at 31 -4
\pinlabel $r_1$ at 45 -4
\pinlabel $r_2$ at 57 -4
\pinlabel $r_3$ at 71 -4
\pinlabel $r_4$ at 97 -4
\pinlabel $(1)$ at 116 49
\pinlabel $(1)$ at 116 63
\pinlabel $(1)$ at 116 76
\pinlabel $(1)$ at 116 89
\endlabellist
\includegraphics[scale=.8]{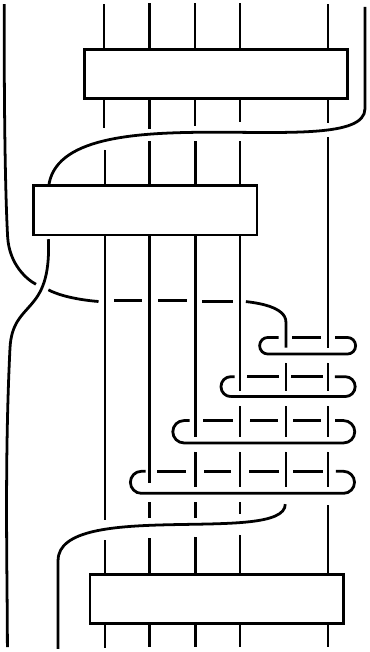}
\end{minipage}.
\end{equation}\vskip7pt
After a Markov move, this last complex is equivalent to:
\begin{equation}\label{eq:D0mvSimp2}
 \begin{minipage}{1.2in}
\labellist
\small
\pinlabel $\gamma_v$ at 56 166
\pinlabel $\KB$ at 38 127
\pinlabel $\omega \gamma_v$ at 56 15
\tiny
\pinlabel $1$ at 6 -4
\pinlabel $r_0$ at 22 -4
\pinlabel $r_1$ at 36 -4
\pinlabel $r_2$ at 49 -4
\pinlabel $r_3$ at 62 -4
\pinlabel $r_4$ at 86 -4
\pinlabel $(1)$ at 106 49
\pinlabel $(1)$ at 106 63
\pinlabel $(1)$ at 106 76
\pinlabel $(1)$ at 106 89
\endlabellist
\includegraphics[scale=.8]{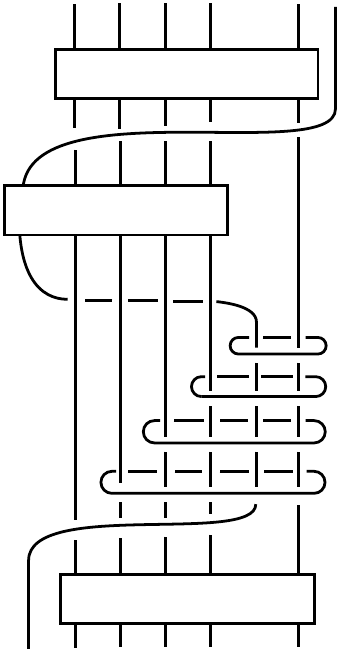}
\end{minipage}
\ \ \ \simeq \ \ \ 
\begin{minipage}{1.2in}
\labellist
\pinlabel $\gamma_v$ at 56 166
\pinlabel $\KB$ at 38 127
\pinlabel $\omega \gamma_v$ at 56 15
\tiny
\pinlabel $1$ at 1 -4
\pinlabel $r_0$ at 17 -4
\pinlabel $r_1$ at 31 -4
\pinlabel $r_2$ at 44 -4
\pinlabel $r_3$ at 57 -4
\pinlabel $r_4$ at 81 -4
\pinlabel $(1)$ at 101 49
\pinlabel $(1)$ at 101 63
\pinlabel $(1)$ at 101 76
\pinlabel $(1)$ at 101 89
\small
\endlabellist
\includegraphics[scale=.8]{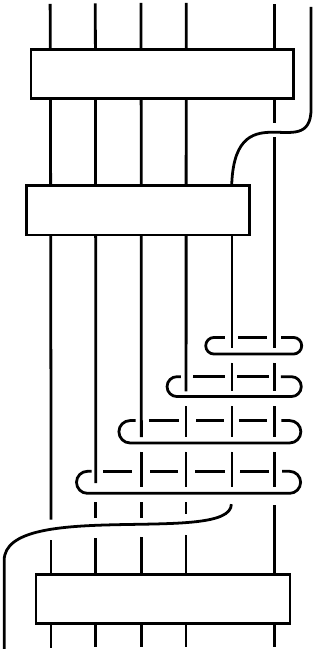}
\end{minipage}
\ \ \ \simeq \ \ \ 
\begin{minipage}{1.2in}
\labellist
\pinlabel $\gamma_v$ at 50 166
\pinlabel $\KB$ at 38 127
\pinlabel $\omega \gamma_v$ at 50 15
\tiny
\pinlabel $1$ at 1 -4
\pinlabel $r_0$ at 17 -4
\pinlabel $r_1$ at 31 -4
\pinlabel $r_2$ at 44 -4
\pinlabel $r_3$ at 57 -4
\pinlabel $r_4$ at 81 -4
\pinlabel $(1)$ at 98 49
\pinlabel $(1)$ at 98 63
\pinlabel $(1)$ at 98 76
\pinlabel $(1)$ at 98 109
\small
\endlabellist
\includegraphics[scale=.8]{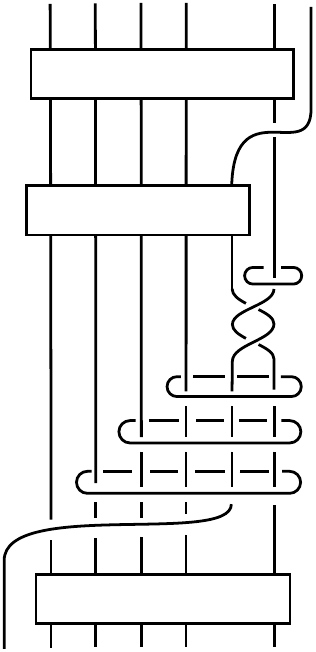}
\end{minipage} .
\end{equation}\vskip7pt
Finally, this last complex is equivalent to $\CB_{v(m-1)}^{(m)}\XB_{n-1}$.  Now we take care of the degree shifts.   The Markov move contributes a factor of $t^{1/2}$ by \ref{prop:markov}.  This shift is cancelled out by our normalization for $\XB$, as in the proof of Proposition \ref{prop:DvMarkov1}.

In the first equivalence in (\ref{eq:D0mvSimp2}) the symmetrizer $\KB$ absorbed several crossings using Proposition \ref{prop:rouquierAbsorbing}.  Each crossing contribuites a shift of $t^{1/2}$.  The number of crossings absorbed in the pictured diagram is $2(r_1+r_2+r_3)$, and in general the number of crossings absorbed will be $2(r_1+\cdots + r_{m-1})$ which is twice the number of indices $i$ such that $v_i<m$.  This completes the proof.
\end{proof}

%=========================
\subsection{The homologies}
\label{subsec:superpoly}
%========================

\begin{definition}\label{def:linkCase}
Let $f_v(q,t,a)=f^{(m)}_v(q,t,a)\in \N[q,t,a,(1-q)\inv]$ denote the polynomials defined by
\begin{enumerate}\setlength{\itemsep}{2pt}
\item[(L0)] $f_{\emptyset} = 1$.
\item[(L1)] $f_{0v} = (t^{\#\{i\:|\: v_i<m\}}+a)f_{v}$.
\item[(L2)] $f_{kv}=t^{\#\{i\:|\: v_i<k\}}f_{v(k-1)}$.
\item[(L3)] $f_{mv}=f_{v(m-1)}+q f_{vm}$.
\end{enumerate}
\end{definition}

Let us convince ourselves that the $\f_{v}$ are well-defined.   We first put an appropriate partial order on the set of sequences $v$.

\begin{definition}\label{def:partialOrder}
Given $v\in \{0,1,\ldots,m\}^n$, we write $l(v)=n$ for the length of $v$, and $|v|:=v_1+\cdots +v_n$.  Let $i(v)$ be the number of pairs of indices $i<j$ with $v_i=m$ and $v_j<m$.  Given two sequences $v,w$ (not necessarily with the same length), we write $v<w$ one of the following is true:
\begin{itemize}
\item $l(w)<l(v)$.
\item $l(w)=l(v)$ and $|w|<|v|$.
\item $l(w)=l(v)$, $|w|=|v|$, and $i(w)<i(v)$.
\end{itemize}
\end{definition}

\begin{lemma}\label{lemma:f_welldefined}
The polynomials $f_v$ are well-defined.
\end{lemma}
\begin{proof}
It is clear that the empty sequence is the unique minimum with respect to the partial order from Definition \ref{def:partialOrder}.  Further, for each $v$ there are only finitely many sequences $w$ with $\emptyset <w< v$.  Now, an application of one of the relations (L$i$) writes $f_v$ in terms of $f_w$ with $w<v$.  This is obvious except possibly in the case $v=m^n$.  In this case, relation (L3) yields
\begin{equation}\label{eq:allMs_f}
\f_{m^n}(q,t,a)=\frac{1}{1-q} \f_{m^{n-1}(m-1)}(q,t,a).
\end{equation}

Since the relations (L$i$) express $f_v$ in terms of $f_w$ with $w<v$, we conclude that the $f_v$ are unique when they exist.  For existence we must prove consistency of the relations (L$i$).  This is clear, since exactly one of the relations (L$i$) can be applied to a given sequence $v$.
\end{proof}

\begin{remark}
The denominators of $\f_{v}$ are all powers of $1-q$.  It is illustrative to normalize these denominators away by defining
\begin{equation}
\tilde{f}^{\: (m)}_{v}\ := \ (1-q)^r \f_{v},
\end{equation}
where $r$ is the number of indices $i$ such that $v_i=m$.   Then clearly the we have the following recursion for $\tilde{f}_v=\tilde{f}_v^{\:(m)}$.
\begin{enumerate}\setlength{\itemsep}{2pt}
\item[($\tilde{\text{L}}0$)] $\tilde{f}_{\emptyset} = 1$.
\item[($\tilde{\text{L}}1$)] $\tilde{f}_{0v} = (t^{\#\{i\:|\: v_i<m\}}+a)\tilde{f}_{v}$ otherwise
\item[($\tilde{\text{L}}2$)] $\tilde{f}_{kv}=t^{\#\{i\:|\: v_i<k\}}\tilde{f}_{v(k-1)}$
\item[($\tilde{\text{L}}3$)] $\tilde{f}_{mv}=(1-q)\tilde{f}_{v(r-1)}+qf_{vm}$.
\end{enumerate}
This implies that $\tilde{f}_{v}$ is a polynomial in $q,t,a$ for all $v\in \{0,1,\ldots,m\}^n$ (though typically the coefficients will be in $\Z$, not $\N$).
\end{remark}

\begin{theorem}\label{thm:linkCase}
There is no $\Z$-torsion in $\HHH(\CB_{v}^{(m)})$.  The Poincar\'e series of this homology equals $\f_{v}$ for all $v\in \{0,1,\ldots,m\}^n$.  In particular, $\f_{m^n}$ is the KR series of the $(n,nm)$ torus link.
\end{theorem}

\begin{proof}
Fix an integer $m\geq 1$.   Let $f'_{v}$ denote the KR series of $\CB_{v}$.  We prove that $f_{v}'=\f_{v}$ by induction on $v$, using an appropriate partial order.  

The base case is trivial.  For the reader's convenience, we also take care of the case $n=1$.  Observe that $\CB_{m}^{(1)}=\one_1$, while $\CB_{k}^{(1)}=\KB_1$ for $0\leq k\leq m-1$.  The Hochcshild homologies are $\Z[x_1,\xi_1]$ and $\Z[\xi_1]$, respectively, where $x_1$ has degree $q$, and $\xi_1$ is an odd variable (i.e.~$\xi^2=0$) of degree $a$.  Thus, $f_{k}'$ equals  $\frac{1+a}{1-q}$ if $k=m$ and $1+a$ otherwise. This agrees with $f_{k}^{(1)}$ in these cases.

Fix a sequence $v\in \{0,1,\ldots,m\}^n$, and assume by induction that $\HHH(\CB_w^{(m)})$ has no $\Z$-torsion, and $f_{w}'=\f_{w}$ for all $w< v$.  In particular, this implies that $\HHH(\CB_{w}^{(m)})$ is supported in even homological degrees for all $w<v$.

If $v=m^n$, then
\[
\CB_{m^{n-1}(m-1)}^{(m)} \simeq (\one_{n-1}\sqcup \KB_1)\otimes \CB_{m^n}^{(m)}.
\]
We have $\HHH(\CB_{m^n}^{(m)})\cong \Z[x_n]\otimes_\Z \HHH(\CB_{m^{n-1}(m-1)})$ by Proposition \ref{prop:K1superpoly}.  The KR series of the LHS is $\f_{m^{n-1}(m-1)}$ by induction, while the KR series of the RHS is $(1-q) f_{m^n}'$.  This shows that $f_{v}'=\f_{v}$ in this case.

So assume that $v\neq m^n$.  We have several possible cases.

Case 1. If $v=0w$, then $\HHH(\CB_{v}^{(m)})=\HHH(\CB_{0w}^{(m)})$ is isomorphic to the direct sum of two copies of $\HHH(\CB_{w}^{(m)})$ with the shifts $t^c$ and $a$, as in Proposition \ref{prop:CvMarkov}.  From this it follows easily that $f_{v}'=\f_{v}$ in this case, using relation (C1) in the definition of $\f_{v}$.

Case 2.  If $v=kw$ with $0<k<m$, then observe $f_{kw}'= t^b f_{w(k-1)}'$ by Proposition \ref{prop:Fkv}.  Comparison with (C2) gives $f_{v}'=f_{v}$ in this case.

Case 3.  If $v=mw$, then $\CB_{v}^{(m)}=\CB_{mw}^{(m)}$ fits into a distinguished triangle
\[
q\CB_{wm}^{(m)} \rightarrow \CB_{v}^{(m)}\rightarrow \CB_{w(m-1)}^{(m)}\rightarrow q\CB_{wm}^{(m)}[1].
\]
by Proposition \ref{prop:Fmv}.  The functor of Hochschild cohomology is just a linear functor extended to complexes, hence it preserves distinguished triangles (i.e.~ it commutes with mapping cones).  Thus, for each $i,j$ we have a distuinguished triangle (in the homotopy category of complexes of bigraded abelian groups)
\[
q\HH(\CB_{wm}^{(m)}) \rightarrow \HH(\CB_{v}^{(m)})\rightarrow \HH(\CB_{w(m-1)^{(m)}})\rightarrow q\HH(\CB_{wm}^{(m)})[1].
\]
Taking homology gives a long exact sequence
\[
\cdots \rightarrow q\HHH(\CB_{wm}^{(m)}) \rightarrow \HHH(\CB_{v}^{(m)})\rightarrow \HHH(\CB_{w(m-1)}^{(m)})\buildrel \d\over \rightarrow q\HHH(\CB_{wm}^{(m)})[1]\rightarrow \cdots
\]
However, the arrow labelled $\d$ is zero, since the source is supported in even homological degrees, while the target is supported in odd homological degrees.  Thus we obtain a short exact sequence of triply graded abelian groups
\[
0 \rightarrow q\HHH(\CB_{wm}^{(m)}) \rightarrow \HHH(\CB_{v}^{(m)})\rightarrow \HHH(\CB_{w(m-1)}^{(m)}) \rightarrow 0
\]
The right-most term has no $\Z$-torsion by the induction hypothesis, hence it is free.  Thus the above short exact sequence splits, and we find that
\[
\HHH(\CB_{v}^{(m)})\ \cong\ \HHH(\CB_{w(m-1)}^{(m)}) \oplus q\HHH(\CB_{wm}^{(m)})
\]
Taking KR series yields $f_{v}'=\f_{m,v}$ in this case, by comparing with the relation (C3) in the definition of $\f_{v}$.  This completes the inductive step and completes the proof.
\end{proof}

Now we do the same thing, this time with the complexes $\XB\CB_{v}$.

\begin{definition}\label{def:knotCase}
Let $g_v=g_{v}^{(m)}\in \N[q,t,a,(1-q)]$ be the polynomials defined by the following recursion.
\begin{enumerate}\setlength{\itemsep}{2pt}
\item[(K0)] $g_{\emptyset} = 1$ and $g_{0}=1+a$.
\item[(K1)] $g_{0mv}=t^{\#\{i\:|\: v_i<m\}}g_{v(m-1)}$, and $g_{0kv} = (t^{1+{\#\{i\:|\: v_i<m\}}}+a)g_{v}$ for $0\leq k\leq m-1$.
\item[(K2)] $g_{kv}=t^{\#\{i\:|\: v_i<k\}}g_{v(k-1)}$
\item[(K3)] $g_{mv}=g_{v(m-1)}+qg_{vm}$.
\end{enumerate}
\end{definition}

The argument for why the $\g_{v}$ are well-defined is nearly identical to what was said about $\f_{v}$.  Again, note that
\[
g_{m^n}^{(m)} = \frac{1}{1-q} g_{m^{n-1}(m-1)}^{(m)}.
\]

\begin{remark}
Note that if $v\neq m^n$, the relations (K$i$) express $g_{v}^{(m)}$ in terms of some other sequence(s), none of which is $(m,\ldots,m)$.  Consequently, $g_{v}^{(m)}\in \N[q,t,a]$ for $v\neq m^n$.  The same is \emph{not} true of $\f_{v}$. 
\end{remark}

\begin{theorem}\label{thm:knotCase}
There is no $\Z$-torsion in $\HHH(\XB\CB_{v}^{(m)})$.  The Poincar\'e series of this homology equals $g_{v}^{(m)}$, for all $v\in \{0,1,\ldots,m\}^n$.  In particular $g^{(m)}_{m^n}$ is the KR series of the $(n,nm+1)$ torus knot and $g_{m^{n-k}(m-1)m^{k-1}}^{(m)}$ is $(1-q)$ times the KR series of the braid closure of $\FT_n^m \sigma_{n-1}\cdots \sigma_{n-k+1}\sigma_{n-k}\inv\cdots\sigma_1\inv$.
\end{theorem}

\begin{proof}
This follows along the same lines as the proof of Theorem \ref{thm:linkCase}.  All the hard work is accomplished by the relations in Propositions \ref{prop:DvMarkov1},\ref{prop:DvMarkov2}, \ref{prop:Fkv}, and \ref{prop:Fmv}.
\end{proof}

%=====================
\subsection{Some special cases}
\label{subsec:specialBraids}
%=====================
We obtain the Khovanov-Rozansky homologies of several infinite families of links as special cases of Theorems \ref{thm:linkCase} and \ref{thm:knotCase}.

Let us introduce some notation.  For any integers $1\leq i\leq j\leq n$, let $\Cox_{[i,j]}=\sigma_{j-1}\cdots\sigma_{i+1}\sigma_i$ and $\otherCox_{[i,j]}=\sigma_i\sigma_{i+1}\dots \sigma_{j-1}$.  Observe that $\Cox_n = \Cox_{[1,n]}$.  If $\b\in \Br_{n-1}$ is any braid, then
\begin{equation}\label{eq:braidSlide}
\Cox_{[1,n]}(\b\sqcup \one_1) \simeq (\one_1\sqcup \b) \Cox_{[1,n]} \hskip .5in \otherCox_{[1,n]}(\one_1\sqcup \b) \simeq (\b\sqcup \one_1) \otherCox_{[1,n]}.
\end{equation}
Diagrammatically these are obvious isotopies:
\[
\begin{minipage}{.7in}
\labellist
\small
\pinlabel $\b$ at 33 13
\endlabellist
\includegraphics[scale=.8]{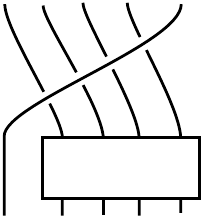}
\end{minipage}
\ \ \simeq \ \ \begin{minipage}{.7in}
\labellist
\small
\pinlabel $\b$ at 22 47
\endlabellist
\includegraphics[scale=.8]{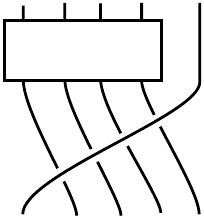}
\end{minipage}
\hskip1in
\begin{minipage}{.7in}
\labellist
\small
\pinlabel $\b$ at 22 13
\endlabellist
\includegraphics[scale=.8]{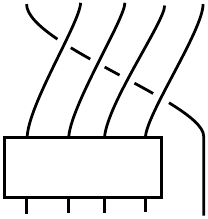}
\end{minipage}
\ \ \simeq \ \ \begin{minipage}{.7in}
\labellist
\small
\pinlabel $\b$ at 33 47
\endlabellist
\includegraphics[scale=.8]{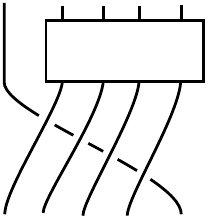}
\end{minipage}
\]

  Let $J_k=\Cox_{[1,k]}\otherCox_{[1,k]}=\sigma_{k-1}\cdots \sigma_2\sigma_1^2\sigma_2\cdots \sigma_{k-1}$ denote the \emph{Jucys-Murphy braid}.  We commonly abuse notation and identify a braid $\b \in \Br_{k}$ with its image $\b\sqcup \one_{n-k} \in \Br_n$.  We denote the Rouquier complex associated to $\b$ also by $\b$.   It is an exercise to show that $\FT_n = J_1J_2\cdots J_n$ (compare with (\ref{eq:FTrecursion})).

\begin{corollary}\label{cor:allTheLinks}
Let $\b_{i,k,m,n}$ denote the braid
\[
X_{[1,n-i]} (\one_1\sqcup \FT_{n-1}^{m-k})\FT_n^k Y_{[1,n-i]}.
\]
Then the Khovanov-Rozansky series of $\b_{i,k,m,n}$ is $\frac{1}{1-q}f_{m^{n-i},k,m^{i-1}}^{(m)}(q,t,a)$, and the Khovanov-Rozansky series of $X_n\b_{i,k,m,n}$ is $\frac{1}{1-q}g_{m^{n-i},k,m^{i-1}}^{(m)}(q,t,a)$.
\end{corollary}
\begin{proof}
From the definition, we have
\[
\CB_{m^{n-i},k,m^{i-1}}^{(m)} = X_{[1,n-i]}  (\KB_1\sqcup \FT_{n-1}^{m-k}) \FT_n^k Y_{[1,n-i]}.
\]
Then the proof follows from Proposition \ref{prop:K1superpoly}.
\end{proof}

\begin{proposition}\label{prop:hookPureBraid}
The Khovanov-Rozansky series of $J_i\inv \FT_n^m$ equals $\frac{1}{1-q}f_{m^{n-i},m-1,m^{i-1}}^{(m)}(q,t,a)$.
\end{proposition}
\begin{proof}
From the definitions,
\begin{eqnarray*}
\b_{i,m-1,m,n}& =& X_{[1,n-i]}  (\one_1\sqcup \FT_{n-1}) \FT_n^{m-1} Y_{[1,n-i]}\\
&=& X_{[n-i,n]}\inv X_{[1,n]}(\one_1\sqcup \FT_{n-1}) \FT_n^{m-1} Y_{[1,n]}Y_{[n-i,n]}\inv\\
&=& X_{[n-i,n]}\inv X_{[1,n] Y_{[1,n]} \FT_n^{m-1}}(\FT_{n-1}\sqcup\one_1)Y_{[n-i,n]}\inv\\
&=& X_{[n-i,n]}\inv \Big(X_{[1,n]} Y_{[1,n]} (\FT_{n-1}\sqcup\one_1)\Big)\FT_n^{m-1}Y_{[n-i,n]}\inv\\
&=& X_{[n-i,n]}\inv \FT_n^m Y_{[n-i,n]}\inv\\
&=& X_{[n-i,n]}\inv Y_{[n-i,n]}\inv  \FT_n^m
\end{eqnarray*}
Up to conjugation with $\HT$, this braid is $J_i\inv \FT_n^m$.  Then Corollary \ref{cor:allTheLinks} completes the proof.
\end{proof}

\begin{proposition}\label{prop:hookBraid}
The Khovanov-Rozansky series of the knot represented by the braid
\[
\FT_n^{m} \Cox_{[i,n]} \otherCox_{[1,i]}\inv \ \ = \ \  \FT_n^m \ \cdot \ \begin{minipage}{1.05in}
\labellist
\small
\pinlabel $\underbrace{\ \ \ \ \ \ \ \ \  \ \ \ \ \   }_{i-1}$ at 37 -9
\pinlabel $\underbrace{\ \ \ \ \ \ \ \ \ \ }_{n-i}$ at 81 -9
\endlabellist
\includegraphics[scale=.8]{diagrams/hookBraid2}
\end{minipage}\ \
\]\vskip7pt
equals $\frac{1}{1-q}g_{m^{n-i},m-1,m^{i-1}}^{(m)}(q,t,a)$.
\end{proposition}

\begin{proof}
By (\ref{eq:FTtimesC}) it suffices to prove this when $m=1$.  By Proposition \ref{prop:hookPureBraid} we know that
\[
\b_{i,0,1,n} \simeq X_{[n-i,n]}\inv Y_{[n-i,n]}\inv  \FT_n.
\]
Now we simplify this braid, using isotopies and conjugation.  Introduce the notation $\JB_{[i,j]}:=\Cox_{[i,j]}\otherCox_{[i,j]}$ to abbreviate, and compute
\begin{eqnarray*}
\Cox_{[1,n]}\Cox_{[n-i+1,n]}\inv \otherCox_{[n-i+1,n]}\inv \FT_n
&\sim & \Cox_{[1,n]}(\Cox_{[n-i-1,n-2]}\inv \otherCox_{[n-i-1,n-2]}\inv) \FT_n\\
&\sim & \cdots\\
&\sim & \Cox_{[1,n]} (\Cox_{[1,i]}\inv \otherCox_{[1,i]}\inv)\FT_n\\
&\simeq & \Cox_{[i,n]}\otherCox_{[1,i]}\inv \FT_n,
\end{eqnarray*}
This follows by several applications of relations of the form
\begin{eqnarray*}
\Cox_{[1,n]}(\one_j\sqcup \b\sqcup \one_{n-k-j}) \FT_n 
&\simeq & (\one_{j-1}\sqcup \b\sqcup \one_{n-k-j-1})\Cox_{[1,n]} \FT_n \\
&\sim &  \Cox_{[1,n]} \FT_n  (\one_{j-1}\sqcup \b\sqcup \one_{n-k-j-1})\\
&\simeq & \Cox_{[1,n]}  (\one_{j-1}\sqcup \b\sqcup \one_{n-k-j-1})\FT_n,
\end{eqnarray*}
which holds for any braid $\b\in \Br_k$ by (\ref{eq:braidSlide}) and the fact that $\FT_n$ is central.  This completes the proof.
\end{proof}

%%%%%%%%%%%%%%%%%%%%%%%%%%%%%%%
\section{Combinatorial formulas}
\label{sec:combinatorics}
%%%%%%%%%%%%%%%%%%%%%%%%%%%%%%%

In this section we give combinatorial formulas for $\f_v$ and $g_v^{(m)}$ and interpret the results in terms of Dyck paths.  In particular we show that $g_{m^n}^{(m)}(q,t,0)=\frac{1}{1-q}C_n^{(m)}(q,t)$.  Given Theorem \ref{thm:knotCase}, this proves a generalization of Gorsky's conjecture (Conjecture \ref{conjecture:gorsky}).

\subsection{The statistics}
\label{subsec:diagrams}
In this section we introduce the combinatorial objects and the statistics that will be involved in the sum formulas formulas for $f^{(m)}_v$ and $g^{(m)}_v$.  Let $m\geq 1$ be a fixed integer.

Let $\sigma\in \Z_{\geq 0}^n$ be given.  We will visualize $\sigma$ as bottom-justified collection of columns, as in the following:
\begin{equation}\label{eq:columnEx}
\begin{minipage}{2in}
\labellist
\small
\pinlabel $1$ at 28 15
\pinlabel $5$ at 53 15
\pinlabel $1$ at 78 15
\pinlabel $2$ at 101 15
\pinlabel $3$ at 125 15
\pinlabel $0$ at 149 15
\pinlabel $4$ at 173 15
\endlabellist
\includegraphics[scale=.7]{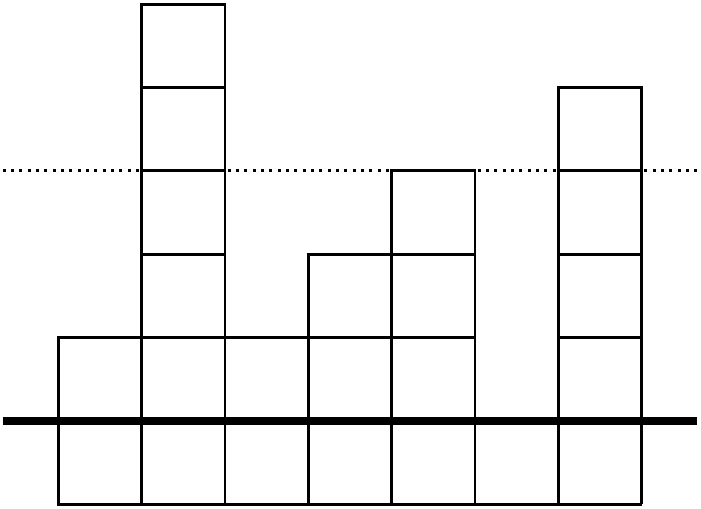}
\end{minipage}
\end{equation}
The above diagram corresponds to $\sigma=(1,5,1,2,3,0,4)$.  The bottom row of boxes has labels given by $\sigma_i$.  Above the $i$-th box we place $\sigma_i$ additional boxes.  For visual aid, we separate the bottom row from the other rows with a solid horizontal line.

Given $\sigma\in \Z_{\geq 0}^n$, let $v(\sigma)\in \{0,1,\ldots,m\}^n$ denote the sequence obtained by \emph{truncating} $\sigma$.  Precisely, $v(\sigma)_i = \sigma_i$ if $\sigma_i<m$, and $v(\sigma)_i=m$ if $\sigma_i\geq m$.

\begin{example}\label{ex:vExample}
Let $\sigma=(1,5,1,2,3,0,4)$, which is pictured via the diagram (\ref{eq:columnEx}), and let $m=3$.  We may visualize $v(\sigma)$ graphically by drawing a dashed line at height $m$ (in which the solid line corresponds to height zero).  Then $v(\sigma)$ counts the numbers of boxes between the solid and dashed lines.  In this example, $v(\sigma)=(1,3,1,2,3,0,3)$.
\end{example}

As a matter of terminology, we will say that $\sigma$ extends $v(\sigma)$.  Now we define some quantities associated with $\sigma$.  Let $\sigma\in\Z_{\geq 0}^n$ be given, and let $D=D(\sigma)$ denote the associated box diagram.  We regard $D$ as a union of \emph{cells}, where a cell is by definition a 1$\times$1 rectangle with corners in $\Z\times \Z$.  For each $c\in D$, let $l(c)$ denote the number of cells above $c$ in $D$ (in combinatorics literature this may be called the \emph{leg} of $c$), and let $l'(c)$ denote the number of cells below $c$ in $D$ (the \emph{co-leg}), so that $l(c)+l'(c)=\sigma_i$ if $c$ is the $i$-th column.  In particular if $c_1,\ldots,c_n$ are the cells in the bottom row, then $l(c_i)=\sigma_i$, by construction.

We put a total order on the set of all possible cells, as follows.  The bottom left corner is the unique minimum, and in general the successor of a cell $c$ is
\begin{itemize}
\item the next cell to the right of $c$ if this exists, or
\item the leftmost cell in the row above $c$ otherwise.
\end{itemize}
The \emph{range} of a cell $c$ is by definition $\range_c:=\{\succ(c),\ldots,\succ^{n-1}(c)\}$.  Note that $\range_c$ consists of all cells strictly to the right of $c$ and in the same row, together with all cells strictly to the left of $c$, in the row above.  We let $\range_c(D):=\range_c\cap D$.

\begin{definition}\label{def:area}
Let $|\sigma|:=\sigma_1+\cdots+\sigma_n$ and $|v|=v_1+\cdots+v_n$.
\end{definition}

\begin{definition}\label{def:dinv}
Let $\sigma\in\Z_{\geq 0}^n$ be given, and let $D=D(\sigma)$ be the associated diagram.  For each $c\in D$, let $\Dinv(D,c)$ denote the set of cells $d\in \range_c(D)$ such that
\begin{equation}\label{eq:dinvPair}
l(c)\leq l(d)\leq m-1.
\end{equation}
Let $\Dinv(D):= \bigsqcup_{c\in D} \Dinv(D,c)$, and $\dinv(D):=\#\Dinv(D)$.  We also write $\dinv(\sigma) = \dinv(D(\sigma))$.
\end{definition}

It is sometimes useful to rephrase $\Dinv(D,c)$ purely in terms of the associated sequence $\sigma$.

\begin{lemma}\label{lemma:dinvLemma}
Let $\sigma\in \Z_{\geq 0}^n$ be given, and let $D=D(\sigma)$ be the associated diagram.  Fix an index $i$, let $c\in D$ be a cell in the $i$-th column with $l(c)=k$, and let $d\in \range_c(D)$ be in the $j$-th column.  Then $(c,d)\in \Dinv(D,c)$ if and only if $0\leq k\leq \min\{m-1,\sigma_1\}$ and
\begin{equation}\label{eq:dinvPair2}
\sigma_i  \leq \bar{\sigma}_j\leq \sigma_i+m-1-k
\end{equation}
where $\bar{\sigma}_j=\sigma_j$ if $i<j$, and $\bar{\sigma}_j=\sigma_j-1$ if $j<i$.  In particular, if $c_i\in D$ is the unique cell in the top of the $i$-th column, then
\begin{equation}\label{eq:dinv_i}
\dinv(D,c_i)=\#\left\{j\in \{1,\ldots,\hat{i},\ldots,n\}\  \bigg|\  \sigma_i\leq \bar{\sigma_j}\leq \sigma_i+m-1\right\}.
\end{equation}
\end{lemma}
\begin{proof}
Let $c\in D$ be in column $i$, and set $k:=l(c)$.  Note that $k\leq \sigma_i$.  Suppose $d\in \range_c(D)$ is in column $j$. If $i<j$, then $l'(d)=l'(c)=\sigma_i-k$ since $c$ and $d$ are in the same row.  If $j<i$, then then $l'(d)=\sigma_i-k+1$.  Thus, $l(d) = \sigma_j-l'(d) = \bar{\sigma}_j -\sigma_i+k$.  We write this as
\begin{equation}\label{eq:lfjewfon;}
\bar{sigma}_j - \sigma_i = l(d)-k
\end{equation}
Now, $d\in \Dinv(D,c)$ if and only if (\ref{eq:dinvPair}) holds, which can be written
\[
k\leq l(d)\leq m-1.
\]
In particular this forces $k\leq m-1$.  Subtracting $k$ from this inequality and using (\ref{eq:lfjewfon;}), we get
\[
0\leq \bar{\sigma}_j-\sigma_i \leq m-1-k.
\]
This  completes the proof.
\end{proof}

The following statistics are defined entirely in terms of $\sigma$.

\begin{definition}\label{def:dinv_i}
For each index $i\in \{1,\ldots,n\}$, let $d(\sigma,i)$ denote the number of pairs $(i,j)$ of indices such that either
\begin{itemize}
\item[$(i)$] $i<j$ and $\sigma_i\leq \sigma_j\leq \sigma_i+m-1$, or
\item[$(ii)$] $j<i$ and $\sigma_i\leq \sigma_j-1\leq \sigma_i+m-1$.
\end{itemize}
Note that $d(\sigma,i)$ equals the right-hand side of (\ref{eq:dinv_i}).  Let $\dinv(\sigma)$ denote the number of triples $(i,j,k)$ such that either $1\leq i<j\leq n$ and $k$ solves the inequalities
\[
\sigma_i\leq \sigma_j \leq \sigma_i+m-1-k \quad \quad 0\leq k\leq \min\{m-1,\sigma_1\},
\]
or $1\leq j<i\leq n$ and $k$ solves
\[
\sigma_i\leq \sigma_j-1\leq  \sigma_i+m-1-k \quad\quad 0\leq k\leq \min\{m-1,\sigma_1\}.
\]
if $j<i$. 
\end{definition}

This definition of $\dinv(\sigma)$ agrees with that in Definition \ref{def:dinv} by Lemma \ref{lemma:dinvLemma}.
\subsection{Summation formula for $f_v^{(m)}$}
\label{subsec:stateSum1}
In this section we give an exact formula for $f^{(m)}_v$.  Let $m\geq 1$ be a fixed integer.  We will heavily use the terminology developed in the previous section.

\begin{theorem}\label{thm:linkStateSum}
For each $v\in \{0,1,\ldots,n\}^m$, we have
\begin{equation}\label{eq:fseries}
\f_v = \sum_{\sigma} q^{|\sigma|-|v|} t^{\dinv(\sigma)+e(v)}\prod_{i=1}^n(1 + at^{-d(\sigma,i)})
\end{equation}
where $e(v)$ denotes the number of pairs $(i<j)$ such that $1\leq v_i \leq v_j \leq m-1$.
\end{theorem}

We will prove this by verifying that the right-hand side satisfies the same recursion as defines $\f_v$.   We first introduce an operation that we refer to as cycling, or rotation of sequences.

\begin{definition}\label{def:rotation}
Let $\sigma\in \Z_{\geq 0}^n$ be given, and let $D=D(\sigma)$ be the associated diagram.  Recall the total order on the set of all possible cells.  If $c$ is any cell other than the bottom left corner, then let $\xi c$ denote the unique cell such that $\succ(\xi c)=c$.  Let $\xi D$ the diagram obtained by applying $\xi$ to $D\setminus\{c_0\}$.  In other words, $\xi D$ is the diagram associated to
\begin{equation}\label{eq:cycling}
\xi(\sigma): = (\sigma_2,\ldots,\sigma_n,\sigma_1-1).
\end{equation}
We say that $\xi D$ (resp.~$\xi\sigma$) is obtained from $D$ (resp.~$\sigma$) by rotation.  If $\sigma_1=0$, it is convenient to define
\[
\xi(\sigma) := (\sigma_2,\ldots,\sigma_n) \ \ \ \ \ \ (\text{ if }\sigma_1=0).
\]
\end{definition}

We make some easy observations about $\xi$.

\begin{observation}\label{obs:vRotation}
Suppose $\sigma=(\sigma_1,\ldots,\sigma_n)$ extends $v=(v_1,\ldots,v_n)$.  If $\sigma_1\leq m$, then $v_1=\sigma_1$ and $\xi\sigma$ extends $\xi v$.   However, if $\sigma_1>m$, then $\xi\sigma$ extends $(v_2,\ldots,v_n,m)$.
\end{observation}

\begin{lemma}\label{lemma:dinvRot}
Let $\sigma\in \Z_{\geq 0}^n$ be given, and set $D=D(\sigma)$.  For each cell $c\in D$ other than the bottom left corner, $d\in \Dinv(D,c)$ if and only if $\xi d\in \Dinv(\xi D,\xi c)$. Consequently
\begin{enumerate}
\item $\dinv(\sigma)=\dinv(\xi\sigma)$ if $\sigma_1\geq m$.
\item $\dinv(\sigma)=\dinv(\xi\sigma)+d(\sigma,1)$ if $0\leq \sigma_1<m$.
\end{enumerate}
\end{lemma}
\begin{proof}
If $c\in D$ is not the bottom left corner, then rotation gives a leg-preserving bijection $\range_c(D)\cong \range_{\xi c}(\xi D)$, hence $\dinv(D,c)- \dinv(\xi D,{\xi c})$ counts the number of cells $d$ such that $d\in \Dinv(D,c_0)$, where $c_0$ is the bottom left corner.  Note that $d\in \Dinv(D,c_0)$ if and only if $d$ is in the bottom row and $\sigma_1\leq \sigma_j\leq m-1$, where $j$ is the column containing $d$.  In particular if $\sigma_1\geq m$, then $\Dinv(D,c_0)$ is empty, and we obtain (1).  If $\sigma_1\leq m-1$, then the cardinality of $\Dinv(D,c_0)$ is $d(\sigma,1)$ by definition, and we obtain (2).
\end{proof}

%\begin{lemma}\label{lemma:vRotation}
%We have
%\[
%\inversions(0w)=\inversions(w)=\inversions(mw)=\inversions(wm)=\inversions(w(m-1)).
%\]
%for all $w\in\{0,1,\ldots,m\}^{n-1}$.
%\end{lemma}
%\begin{proof}
%Clear.
%\end{proof}

\begin{proof}[Proof of Therem \ref{thm:linkStateSum}]
Fix a sequence $v\in \{0,1,\ldots,m\}^n$, and let $\sigma\in \Z_{\geq 0}^n$ be a sequence extending $v$.  Let $f'_\sigma$ be the contribution of $\sigma$ to the right-hand side of (\ref{eq:fseries}).  Let $f'_v=\sum_\sigma f'_\sigma$ be the right-hand side of (\ref{eq:fseries}).  The proof splits up into three cases, depending on $v_1$.  In any case, let $w=(v_2,\ldots,v_n)$. % We will freely use Lemmas \ref{lemma:vRotation}, \ref{lemma:areaRot}, and \ref{lemma:dinvRot} without explicit reference.
\vskip10pt
\textbf{Case 1.}  Suppose first that $v_1=0$, hence $v=0w$.  Any $\sigma$ extending $v$ satisfies $\sigma_1=0$, hence $\sigma\mapsto \xi \sigma :=(\sigma_2,\ldots,\sigma_n)$ gives a bijection between sequences extending $v=0w$ and sequences extending $w$.  We have $|\sigma|=|\xi\sigma|$, $|v|=|w|$, and $\dinv(\sigma)=\dinv(\xi\sigma)+d(\sigma,1)$ by Lemma \ref{lemma:dinvLemma}.  Also, $e(v)=e(w)$.  Thus
\[
f'_{\sigma} = t^{d(\sigma,1)}(1+a t^{-d(\sigma,1)})f'_{\xi\sigma},
\]
where the extra factor is the contribution of $i=1$ to the factor $\prod_{i=1}^n (1+ at^{-d(\sigma,i)})$, which contributes to $f'_{\sigma}$ but not to $f'_{\xi \sigma}$. Since $\sigma_1=0$, we have $d(\sigma,1)=\#\{j\in\{2,\ldots,n\}\:|\: 0\leq v_j\leq m-1\}$, hence $f'$ satisfies (L1) of Definition \ref{def:linkCase}.
%\[
%f'_{0w}=(t^{\#\{j\:|\: w_j\leq m-1\}}+a) f'_w.
%\]

\vskip10pt
\textbf{Case 2.} Assume $0<v_1<m$, and let $\sigma$ extend $v$.  Then $\sigma_1=v_1$.  Rotation again gives a bijection between sequences extending $v$ and sequences extending $\xi v$.  We have $|\sigma|=\xi\sigma|+1$, $|v|=|\xi v|+1$,  $\dinv(\sigma)=\dinv(\xi\sigma)+d(\sigma,1)$ by Lemma \ref{lemma:dinvRot}.

Note that if $\sigma_1\leq \sigma_j\leq m-1$, then $\sigma_j=v_j$, since $\sigma$ is supposed to extend $v$.  We claim that
\[
e(v)=e(\xi v) - d(\sigma,1)+\{j\in \{2,\ldots,n\}\:|\: v_j<v_1\}.
\]
Indeed, if $v_1\leq v_j\leq m-1$, then $(1,j)$ contributes to $e(v)$ but not $e(\xi\sigma)$, while if $v_j<v_1$, then $(j,n)$ contibutes to $e(\xi v)$ but not $e(v)$.  Thus,
\begin{eqnarray*}
f'_{\sigma}&=&q^{\area(\sigma)}t^{\dinv(\sigma)+e(v(\sigma))}\prod_{i=1}^n (1+ at^{-d(\sigma,i)})\\
&=&q^{\area(\xi\sigma)}t^{\dinv(\xi\sigma)+e(v(\xi\sigma))+\#\{j\:|\:v_j<v_i\}} \prod_{i=1}^n (1+ at^{-d(\xi\sigma,i)})\\
&=& t^{\#\{j\:|\:v_j<v_i\}}f'_{\xi\sigma}.
\end{eqnarray*}
This shows that $f'$ satsifies (L2) of Definition \ref{def:linkCase}.

\vskip10pt
\textbf{Case 3.}  Assume $v_1=m$, and let $\sigma$ extend $v=mw$ (so $\sigma_1\geq m)$.  There are two subcases.   If $\sigma_1=m$, then $\xi \sigma =(\sigma_2,\ldots,\sigma_n,m-1)$ extends $w(m-1)$.  In this case, we have $|\sigma|=|\xi \sigma|+1$ and $|mw|=|w(m-1)|+1$, and Lemma \ref{lemma:dinvRot} implies that $\dinv(\sigma)=\dinv(\xi\sigma)$.  Further, it is clear that $e(mw)=e(w(m-1))$.  Thus, the contribution of $\sigma$ with $\sigma_1=m$ is
\[
f'_\sigma = f'_{\xi \sigma} \ \ \ \ \ \ \ \ \ \ (\sigma_1=m).
\]

On the other hand, if $\sigma_1>1$, then $\xi \sigma_1$ extends $wm$.  In this case $|\sigma|=|\xi\sigma|+1$, $|mw|=|wm|$, and Lemma \ref{lemma:dinvRot} implies that $\dinv(\sigma)=\dinv(\xi\sigma)$.  Again it is clear that $e(mw)=e(wm)$.  Thus,
\[
f'_\sigma = qf'_{\xi \sigma} \ \ \ \ \ \ \ \ \ \ (\sigma_1>m).
\]
Taking the sum over all $\sigma$ extending $v=mw$ shows that $f'$ satisfies (L3) from Definition \ref{def:linkCase}.

Thus, $f'_v$ satisfies the same recursion as defines $\f_v$, hence $f'_v= \f_v$.
\end{proof}

\subsection{Summation formula for $g_v^{(m)}$}
\label{subsec:stateSum2}
Now we adapt the above to a state sum formula for $g_v^{(m)}$.  

\begin{definition}\label{def:cyclicCondition}
Fix integers $m,n\geq 1$.  Let $v\in \{0,1,\ldots,m\}^n$ be given, and let $\sigma\in \Z_{\geq 0}$ extend $v$.  We say that $\sigma$ satisfies the \emph{cyclic $m$-Dyck condition} if
\begin{itemize}
\item $\sigma_{i+1}\leq \sigma_i+m$ for $1\leq i\leq n-1$.
\item $\sigma_1-1 \leq \sigma_n+m$.
\end{itemize}
Let $\Seq_m$ denote the the set of all sequences $\sigma\in \Z_{\geq0}^n$ satisfying the $m$-Dyck condition, and let $\Seq_m(v)\subset \Seq_m$ be the subset consisting of sequences extending $v$.
\end{definition}

\begin{definition}\label{def:peak}
Let $\sigma\in \Seq_m(v)$ be given.  An index $i$ is called a \emph{peak} of $\sigma$ if
\begin{itemize}
\item $i\neq n$ and $\sigma_{i+1}<\sigma_j+m$, or
\item $i=n$ and $\sigma_1-1 <\sigma_n+m$.
\end{itemize}
In other words, $i$ is a peak if the corresponding inequality in Definition (\ref{def:cyclicCondition}) is strict.  Let $\chi_i(\sigma)$ denote $1$ if $i$ is a peak, and zero otherwise.
\end{definition}

\begin{lemma}\label{lemma:DyckRotation}
If $\sigma$ satisfies the cyclic $m$-Dyck condition, then so does $\xi \sigma$.  Moreover, $i$ is a peak of $\sigma$ if and only if then $i-1$ (mod $n$) is a peak of $\xi \sigma$.
\end{lemma}
\begin{proof}
Straightforward.
\end{proof}
K
\begin{theorem}\label{thm:stateSumKnot}
For each $v\in\{0,1,\ldots,m\}^n$ we have
\begin{equation}\label{eq:knotStatesum}
g_v^{(m)} = \sum_{\sigma} q^{|\sigma|-|v|} t^{\dinv(\sigma)+e(v)} \prod_{i=1}^n (1+\chi_i(\sigma) at^{-d(\sigma,i)}),
\end{equation}
where the sum is over $\sigma\in \Seq_m(v)$, and $e(v)=\{(i<j)\:|\: 1\leq v_i\leq v_j\leq m-1\}$.
\end{theorem}
The proof is quite similar to the proof of Theorem \ref{thm:linkStateSum}.
\begin{proof}
Let $g'_\sigma$ denote the contribution of $\sigma$ to the right-hand side, and let $g'_v=\sum_\sigma g'_\sigma$ denote the right-hand side.  We prove that $g'_v$ satisfies the same recursion as defines $g_v^{(m)}$.  The proof splits up into several cases.  First, observe that if $\sigma_1=0$, then $\sigma_2\leq m$ is forced by the $m$-Dyck condition.  The proof splits into 4 cases, the first three of which are proven the same way as in the proof of Theorem \ref{thm:linkStateSum}.

\vskip10pt
\noindent
\textbf{Case 1.} Assume that $v_1=0$ and $v_2=k<m$, and write $v=(0kw)$.  Let $\sigma$ extend $v$.  Then $\sigma_1=0$ and $\sigma_2=k$.  In this case $i=1$ is a peak, and the argument that $g'_v$ satisfies (K1) (the second half) from Definition \ref{def:knotCase} proceeds exactly as in Case 1 of the proof of Theorem \ref{thm:linkStateSum}.
\vskip10pt
\noindent
\textbf{Case 2.} Assume that $v_1\in \{1,\ldots,m-1\}$.  In this case the proof proceeds almost verbatim as in Case 2 in the proof of Theorem \ref{thm:linkStateSum}.  The only difference here is the appearance of the $\chi_i(\sigma)\in \{0,1\}$, which indicates when $i$ is a peak of $\sigma$.  With appropriate modifications, the argument is the same.  We conclude that $g'_v$ satisfies (K2) from Definition \ref{def:knotCase}.
\vskip10pt
\noindent
\textbf{Case 3} Assume that $v_1=m$.  Once again, arguing as in Case 3 in the proof of Theorem \ref{thm:linkStateSum} shows that $g'_v$ satisfies (K3) from Definition \ref{def:knotCase}.  We omit the details, since they are straightforward.
\vskip10pt
\noindent
\textbf{Case 4.}  The final case is the most interesting.  Assume that $v_1=0$ and $v_2=m$, and write $v=0mw$.  Let $\sigma$ extend $v$.  Then $\sigma_1=0$ and $\sigma_2=m$ is forced by the $m$-Dyck condition.

We have
\[
\xi\sigma = (m,\sigma_3,\ldots,\sigma_n), \ \ \ \ \ \ \ \ \ \ \ \xi^2\sigma = (\sigma_3,\ldots,\sigma_n,m-1),
\]
and $\sigma\mapsto \xi^2\sigma$ defines a bijection $\Seq_m(0mw)\cong \Seq_m(w(m-1))$.  We have $\area(\xi^2\sigma)=\area(\sigma)$,  $e(0mw)=e(w(m-1))$, and
\begin{eqnarray*}
\dinv(\sigma)&=&\dinv(\xi^2\sigma)+\#\{j\geq 2\:|\: \sigma_1\leq \sigma_j\leq \sigma_1+ m-1\}\\
&=&\dinv(\xi^2\sigma)+\#\{j\geq 3\:|\: 0\leq v_j\leq  m-1\}\\
\end{eqnarray*}
since $\sigma_1=0$, $\sigma_2=m$, and $\sigma_j\leq m-1$ if and only if $\sigma_j=v_j$.

Now, we compare the extra factors of the form $\prod_{i=1} (1+\chi_i a t^{-d(\sigma,i)})$.  Since $1$ is not a peak, the $i=1$ factor equals $1$, and the other factors are the same as those in $g'_{\xi^2\sigma}$.  The result is
\[
g'_{\sigma}=t^{\#\{j\geq 3\:|\: 0\leq v_j\leq m-1\}} g'_{\xi^2\sigma}.
\]
Taking the sum over all $\sigma\in \Seq_m(v)$ shows that $g'_v$ satisfies (K1) from Definition \ref{def:knotCase}.  Thus $g'_v=g_v^{(m)}$, since they satisfy the same recursion.
\end{proof}

\begin{corollary}\label{cor:catalan}
The Poincar\'e polynomial $\PC_{T(n,nm+1)}(q,t,a)$ at $a=0$ equals $\frac{1}{1-q}C^{(m)}_n(q,t)$, the $m$-th higher $q,t$ Catalan polynomial.
\end{corollary}
\begin{proof}
We take equation (98) in \cite{HHLRU05} as our definition of $C^{(m)}_n(q,t)$.  The authors conjectured that $C^{(m)}_n(q,t)$ can also be expressed by equation (99) in the same paper, and they showed that this conjecture is a consequence their conjectured formula for $\nabla^m e_n$ (called the $m$-shuffle conjecture; see \S \ref{sec:shuffle}).  This conjecture is now proven \cite{CarMel-pp,MellitRational-pp}, hence (99) in \cite{HHLRU05} is valid.  But this formula is exactly equal to our sum formula for $g_{m^{n-1},m-1}^{(m)}(q,t,0)$ up to swapping $q$ with $t$, as is easily checked.  On the other hand, the shuffle theorem proves $q,t$ symmetry of $C_n(q,t)^{(m)}$, which completes the proof.
\end{proof}

This result is generalized in the next section to an interpretation of each $g_v^{(m)}$ in terms of pieces of the $m$-shuffle conjecture.

%%%%%%%%%%%%%%%%%%%%%%%%%%%%%%%%%%%
\section{Connection to the shuffle conjecture}
\label{sec:shuffle}
%%%%%%%%%%%%%%%%%%%%%%%%%%%%%%%%%%%
We interpret our results in terms of the $m$-shuffle conjecture \cite{HHLRU05}, which is now a theorem \cite{CarMel-pp,MellitRational-pp}.

\subsection{Motivation}
\label{subsec:motivation}
There is a well-known relationship between the HOMFLY-PT polynomial and symmetric functions, which follows from work of Turaev's \cite{TuraevTorus} and has been studied extensively \cite{MortonPowerSums,Lukac05}.  We summarize this below, and explain how our work fits into a conjectural categorification of this story.

We recall some basics of symmetric functions.  If $\F$ is a field of characteristic zero, we let $\L_\F=\F[x_1,x_2,\ldots]^{\Sym}$ denote the ring of symmetric functions.  The special cases $\F=\Q(q)$ and $\F=\Q(q,t)$ will be denoted by $\L_q$ and $\L_{q,t}$.  As graded algebras we have
\[
\L_\F \cong \F[e_1,e_2,\ldots]\cong \F[h_1,h_2,\ldots]\cong \F[p_1,p_2,\ldots]
\]
where $e_i$, $h_i$, and $p_i$ are the elementary, complete, and power sum symmetric functions.  We also have the Schur basis $\{s_\l\}\subset \L_\F$.  Of particular importance in the connections with Hilbert schemes, representation theory, and knot homology are the modified Macdonald polynomials $\tilde{H}_\mu\in \L_{q,t}$.  These are certain integral elements, indexed by partitions, and they form a linear basis of $\L_{q,t}$.   The Garsia-Bergeron operator  $\nabla:\L_{q,t}\rightarrow \L_{q,t}$,  introduced in \cite{BerGars96},  is defined by $\nabla \tilde{H}_\mu = t^{n(\mu)}q^{n(\mu^t)}\tilde{H}_\mu$, where $n(\mu)=\sum_i(i-1)\mu_i$ and $\mu^t$ is the transpose partition. 

Let $\ev:\L_\F\rightarrow \F[a]$ denote the unique $\F$-algebra map sending $e_n\mapsto 1+a$ for all $n\geq 1$.  Equivalently
\begin{equation}\label{eq:evaluation}
\ev(f) = \sum_{k=0}^n \ip{f,h_ke_{n-k}}a^k.
\end{equation}

In order to explain the relevance of the above to knot invariants, we now recall some basics of the HOMFLY-PT polynomial  One way defining HOMFLY-PT polynomial is in terms of the Jones-Ocneanu trace for Hecke algebras.  The Hecke algebra $\Hecke_n$ can be defined as the group algebra $\Q(q)[\Br_n]$ of the braid group, modulo the two-sided ideal generated by $(\sigma_i-1)(\sigma_i+q)$ for $1\leq i\leq n-1$.  The Jones-Ocneanu trace \cite{Jones87} is a $\Q(q)$-linear map $\Tr_{\JO}:\Hecke_n\rightarrow \Q(q)[a]$, such that $\Tr_{\JO}(\b\gamma)=\Tr_{\JO}(\gamma\b)$, and $\Tr_{\JO}(\b)$ is the HOMFLY-PT polynomial of $\b$, up to normalization.

Recall that, for any $k$-algebra $A$, the zeroth Hochschild homology is $\HH_0(A)=A/[A,A]$, where $[A,A]\subset A$ is the $k$-submodule spanned by the commutators $ab-ba$.  Since $\Tr_{\JO}(\b)$ is a trace, it factors through the \emph{universal trace} $\Hecke_n \rightarrow \HH_0(\Hecke_n)$.   Thus, to each conjugacy class of braid $\b\in \Br_n$ we have a well-defined element $[\b]\in \HH_0(\Hecke_n)$.

Juxtaposition of braids gives a bilinear maps $\sqcup:\Hecke_i\times \Hecke_j\rightarrow \Hecke_{i+j}$.  This operation makes $\bigoplus_{n\geq 0}\HH_0(\Hecke_n)$ into a (commutative) graded algebra, via $[\b][\gamma]:=[\b\sqcup \gamma]$.  In \cite{TuraevTorus}, Turaev proved that $\bigoplus_{n\geq 0}\HH_0(\Hecke_n)$ is isomorphic to the free commutative algebra on $[X_1]$, $[X_2]$, $\ldots$, where $X_n=\sigma_{n-1}\cdots\sigma_2\sigma_1$.  This is a $q$-analogue of the fact that permutations modulo conjugation are determined by cycle-type, hence $\bigoplus_{n\geq 0}\HH_0(\Q[S_n])\cong \Q[c_1,c_2,\ldots]$, where $c_n=[s_{n-1}\cdots s_2s_1]$ is the class of an $n$-cycle.

Thus, $\bigoplus_{n\geq 0}\HH_0(\Hecke_n)\cong \L_q$ as graded algebras.  There is a preferred choice of isomorphism\footnote{Our preferred isomorphism differs from the preferred choice of Morton and coathors \cite{MortonPowerSums,Lukac05}.  For the combinatorially inclined, we mention that the two are related by the plethysic substition $f\mapsto f[(1-q)\xx]$}, given by:
\[
\Phi:\bigoplus_{n\geq 0}\HH_0(\Hecke_n)\rightarrow \L_q \ \ \ \ \ \ \ \ \ \ \ [X_n]\mapsto \frac{1}{1-q} e_n.
\]
Thus, associated to each braid, we have $\Phi(\b)\in \L_q$, which depends only on the conjugacy class of $\b$.  Topologically, braids modulo conjugation can also be thought of as annular closures of braids, which is how Turaev's result was originally stated.  Thus, there is an invariant of annular braid closures taking values in $\L_q$.   The Jones-Ocneanu trace can then be recovered as
\begin{equation}\label{eq:traceAndEval}
\Tr_{\JO}(\b) = \ev(\Phi(\b)).
\end{equation}

Inserting full-twists gives an operation $[\b]\mapsto [\FT_n \b]$ on annular braid closures, which is well-defined on conjugacy classes since $\FT_n$ is central in the braid group.  The operator which corresponds to $\FT_n$ on the symmetric function side is the operator $\nabla$ of Garsia-Bergeron  at $t=q\inv$, i.e.~$\Phi(\FT_n \b) = \nabla|_{t=q\inv} \Phi(\b)$.

Work of Gorsky-Negut-Rasmussen \cite{GNR16} suggests a categorification of the above, as mentioned in \S \ref{subsec:hilb} of the introduction.   There is a combinatorial shadow of their work which is essentially a $t$-deformation of the above story.  Based on their work the following is sensible.

\begin{definition}\label{def:GNRshadow}
Let $C\in \KC^b(\SBim_n)$ be given.   Let us say that $f\in \L_{q,t}$ is a \emph{combinatorial shadow} of $C$ if
\[
\PC_{\FT_n^m\otimes C}(q,t,a) = \sum_{k=0}^n \ip{\nabla^m f, h_k e_{n-k}}a^k
\]
for all $m\geq 0$.
\end{definition}

In this section we prove that $e_n$ is a combinatorial shadow of $\XB_n\in\KC^b(\SBim_n)$.  For the comparison with $\nabla^m e_n$ we utilize the famous shuffle conjecture from combinatorics, which we now spend some time recalling.

%===================================
\subsection{$m$-Dyck paths}
\label{subsec:mDyck}
%===================================

Fix integers $n,m\geq 1$.  Let $R\subset \R^2$ denote the $n\times m$ rectangle with its bottom left corner at the origin, aligned with the $x$ and $y$ axes, with height $n$ and width $m$.  An $(n,m)$ Dyck path is a path in $R$ consisting of $n$ North steps and $m$ East steps, beginning at the origin, and terminating in at the top right corner of $R$, which stays weakly above the diagonal.

\begin{definition}\label{def:Dyck}
Let $\Dyck(n,m)$ denote the set of $(n,m)$ Dyck paths.
\end{definition}

\begin{example}\label{ex:Dyck}
The following is an example of an $(4,8)$ Dyck path.
\[
\ig{1}{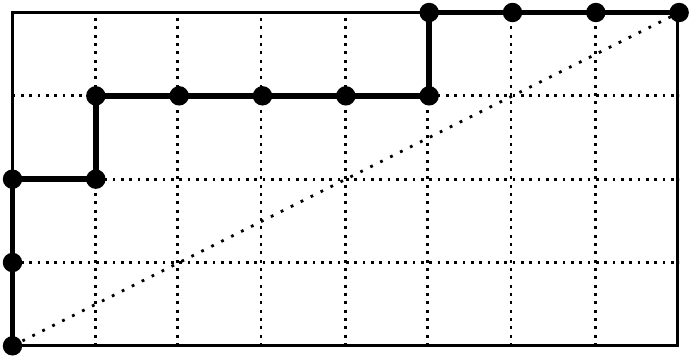}
\]
\end{example}

In this paper we only care about $(n,nm)$ Dyck paths.  These are also called $m$-Dyck paths.  The following is clear.

\begin{lemma}
There are canonical bijections relating the following sets.
\begin{itemize}
\item $(n,nm+1)$-Dyck paths and $(n,nm)$-Dyck paths.
\item  $(n,nm-1)$-Dyck paths and $(n,nm)$-Dyck paths which touch the diagonal exactly once, at the origin.
\end{itemize}\qed
\end{lemma}

\begin{definition}\label{def:diagonals}
Given an $(n,nm)$ Dyck path $D$, let $\gamma_i\in \{0,1,\ldots,m\}^n$ denote the horizontal distance between the beginning of the $i$-th north step and the diagonal $x=my$.  In other words, if $(x_i,y_i)$ are the coordinates of the beginning of the $i$-th North step in $D$, then
\[
\gamma_i=m y_i - x_i.
\]
\end{definition}
In particular $\gamma(D)_1=0$ for all $D$.

\begin{example}\label{ex:Diag}
The Dyck path $D$ from Example \ref{ex:Dyck} has $\gamma(D) = (0,2,3,1)$
\end{example}

\begin{lemma}\label{lemma:Dyckbijection}
The mapping $D\mapsto \gamma(D)$ is a bijection between $(n,mn)$ Dyck paths and sequences $\gamma\in \Z_{\geq 0}^n$ such that
\begin{itemize}
\item[(i)] $\gamma_1=0$.
\item[(ii)] $\gamma_{i+1}\leq \gamma_i+m$.
\end{itemize}
\end{lemma}
\begin{proof}
Let $D$ be an $(n,nm)$ Dyck path, let $(x_i,y_i)$ be the coordinates of the beginnings of its north steps, and let $\gamma_i=my_i-x_i$.   Since all steps have length 1, we see that $y_i=mi-1$, hence $\gamma_i=mi-m-x_i$.  The Dyck path condition is equivalent to
\[
x_1=0 \hskip.5in\text{ and } \hskip.5in x_i\leq x_{i+1}\leq mi,
\]
which is equivalent to 
\[
\gamma_1=0 \hskip.5in\text{ and } \hskip.5in mi-m-\gamma_i \leq mi-\gamma_{i+1}\leq mi,
\] 
which is clearly equivalent to (i) and (ii) from the statement.  This completes the proof.
\end{proof}

\begin{definition}\label{def:DyckSequence}
Let $\Dyck_m\subset \Z_{\geq 0}^n$ denote the subset consisting of those sequences satisfying (i) and (ii) in the statement of Lemma \ref{lemma:Dyckbijection}.  We refer to elements of $\Dyck_m$ as $m$-Dyck sequences.
\end{definition}

\subsection{Parking functions}
\label{subsec:parkingFn}
Let $D$ be an $m$-Dyck path of height $n$, and let $\gamma=\gamma(D)\in \Z_{\geq 0}^n$ be the associated sequence.

Let $\pi\in \{1,\ldots,n\}^n$ be a permutation on $n$ letters, presented in one-line notation.  If $\gamma$ is an $m$-Dyck sequence, we say that $\pi\in \PF(\gamma)$ if $\gamma_{i+1}=\gamma_i+m$ implies $\pi_i<\pi_{i+1}$.
\begin{remark}
Suppose we are given $\pi\in \PF(\gamma)$.  For each $1\leq i\leq n$, let $c_i$ denote the box immediately to the right of the $i$-th north segment of the associated Dyck path, and label this box with the integer $\pi_i$.  Then our condition on $\pi$ ensures that these labels are strictly increasing as one goes up along a repeating sequence of north steps.
\end{remark}
The pairs $(\gamma,\pi)$ with $\pi\in \PF(\gamma)$ are in bijection with $m$-parking functions.  We will visualize $(\gamma,\pi)$ as a labelled box diagram as in (\ref{subsec:diagrams}).  Consider a bottom-justified collection of boxes, arranged in $n$ columns, with heighs equal to $\gamma_{i+1}$.  The boxes in the bottom row will be labelled with the integers $\gamma_1,\ldots, \gamma_n$, and the cell above the topmost box in the $i$-th column will be labelled by $\pi_i$. 

\begin{example}\label{ex:mDyck}
If $m=2$, there is a parking function $(\gamma,\pi)$ given by $\gamma=(0,2,3,1,3,5,4)$, $\pi=(2,5,3,4,6,7,1)$.  This will be pictured as
\[
\begin{minipage}{2in}
\labellist
\small
\pinlabel $0$ at 28 15
\pinlabel $2$ at 53 15
\pinlabel $3$ at 78 15
\pinlabel $1$ at 101 15
\pinlabel $3$ at 125 15
\pinlabel $5$ at 149 15
\pinlabel $4$ at 173 15
\pinlabel $2$ at 28 35
\pinlabel $5$ at 53 83
\pinlabel $3$ at 78 106
\pinlabel $4$ at 101 59
\pinlabel $6$ at 125 106
\pinlabel $7$ at 149 154
\pinlabel $1$ at 173 130
\endlabellist
\includegraphics[scale=.7]{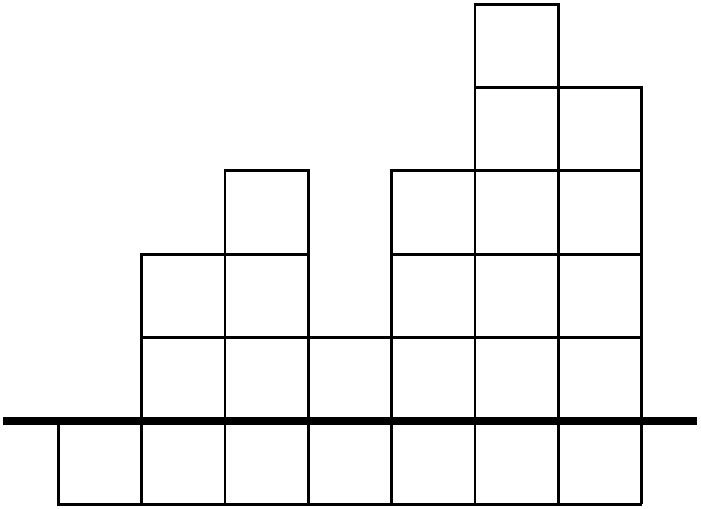}
\end{minipage}.
\]
\end{example}

The following notions are standard in the combinatorics literature \cite{HHLRU05}.  These will be used in \S \ref{subsec:quasisymmetric} in the definition of some special symmetric functions.
\begin{definition}\label{def:dinvParking}
Given $\gamma\in \Z_{\geq 0}^n$ and $\pi\in \Z_{\geq 1}^n$, let $\Dinv_1(\gamma)$ denote the set of pairs of indices $(i,j)$ such that either of the following conditions hold:
\begin{itemize}
\item $i<j$ and $\gamma_i=\gamma_j$.
\item $j<i$ and $\gamma_j-1 = \gamma_i$.
\end{itemize}
Let $\Dinv_1(\gamma,\pi)$ denote the subset of $\Dinv_1(\gamma)$ consisting of pairs $(i,j)$ with $\pi_i<\pi_j$.  Let $\dinv_1(\gamma):=\#\Dinv_1(\gamma)$ and $\dinv_1(\gamma,\pi):=\#\Dinv_1(\gamma,\pi)$.
\end{definition}
We also have a generalization of $\dinv$ to $\dinv$ for $m\geq 1$. 

\begin{definition}
Given $\gamma\in \Z_{\geq 0}$, let $m\ast \gamma\in \Z_{\geq 0}^{mn}$ denote the \emph{$m$-expansion} of $\gamma$, defined as follows:
\[
m\ast \gamma := (\gamma_1, \gamma_1+1,\ldots,\gamma_1+m-1, \ldots,\ldots, \gamma_n, \gamma_n+1,\ldots, \gamma_n+m-1).
\]
Define $m\ast \pi\in \Z_{\geq 1}^{mn}$ by
\[
m\ast \pi = (\:\underbrace{\pi_1,\ldots,\pi_1}_m,\ldots, \underbrace{\pi_n,\ldots,\pi_n}_m\:).
\]
Now, set $\dinv(\gamma,\pi):=\dinv_1(m\ast \gamma, m\ast \pi)$.
\end{definition}

\begin{example}
Let us continue Example \ref{ex:mDyck}.  Expanding gives
\[
\left(\begin{tabular}{c}
$2\ast \gamma$\\
$2\ast \pi$
\end{tabular}
\right) \ \ = \ \ \left( \begin{tabular}{cccccccccccccc}
0&1&2&3&3&4&1&2&3&4&5&6&4&5 \\  2 &2& 5& 5& 3& 3& 4& 4& 6& 6& 7& 7& 1& 1\end{tabular}\right).
\]
\end{example}

It will be useful to have an alternate characterization of $\dinv(\gamma,\pi)$.  Note that $m\ast \gamma$ is given by
\[
(m\ast \gamma)_i = \gamma_j+k,
\]
where $i=1+k+(j-1)m$ with $0\leq k\leq m-1$ and $1\leq j\leq n$.  From the definitions, if $\pi\in \PF(\gamma)$, then $\dinv(\gamma,\pi)$ counts quadruples $(j,k,j',k')$ such that
\begin{enumerate}
\item $0\leq k,k'\leq m-1$ and $1\leq j,j'\leq n$.
\item $\pi_j<\pi_{j'}$.
\item Either of the following conditions hold:
\begin{enumerate}
\item $j<j'$ and $\gamma_j+k = \gamma_{j'}+k'$, or
\item $j'<j$ and $\gamma_j+k = \gamma_{j'}+k'-1$.
\end{enumerate}
\end{enumerate}

\begin{definition}\label{def:attacks}
Say index $i$ \emph{attacks} $j$ if either $i<j$ and $\gamma_i\leq \gamma_j$, or $\gamma_i\leq \gamma_j-1$.  We may put a partial order $<_a$ on the set of indices, by declaring that $i<_a j$ if $i$ attacks $j$.  Given a pair of indices $(i,j)$ with $i<_a j$ let us define $\gamma\bar{\gamma_j}=\gamma_j$ if $i<j$ and $\bar{\gamma}_j = \gamma_j-1$ if $j<i$.
\end{definition}
The notation $\bar{gamma}_j$ is convenient shorthand, but we must be careful to clearly indicate the distinguished index $i$, which we are suppressing from the notation.  The following is straightforward.
\begin{lemma}\label{lemma:parkingDinv}
If $i$ attacks $j$, then the contribution of $(i,j)$ to $\dinv(\gamma,\pi)$ is the number of elements $k\in\{0,1,\ldots,m-1\}$  such that
\begin{itemize}\setlength{\itemsep}{3pt}
\item[$(i)$] if $\pi_i<\pi_j$, then $\gamma_i\leq \bar{\gamma}_j+k\leq \gamma_i+m-1$.
\item[$(ii)$] if $\pi_i>\pi_j$, then $\gamma_i\leq \bar{\gamma}_j+k+1\leq \gamma_i+m-1$.
\end{itemize} \qed
\end{lemma}

Given $\pi\in \PF(\gamma)$, the \emph{reading word} of $\pi$ is the permutation $w(\pi)$ given by reading the entries of $\pi$ right-to-left, row by row, starting at the top row.  In other words, utilizing the total order on cells discussed in \S \ref{subsec:diagrams}, the reading word is the word $(w_1,\ldots,w_n)$ uniquely characterized by:
\begin{itemize}
\item $w$ is a permutation of $\pi$.
\item the box with label $w_{i+1}$ comes before the box with label $w_i$ in the total order on boxes.
\end{itemize}

Note that if $i$ attacks $j$, then $\pi_j$ comes before $\pi_i$ in the reading word.  An important special case of Lemma \ref{lemma:parkingDinv} is the following.

\begin{corollary}\label{cor:parkingDinvTrivial}
Let $\gamma\in \Dyck_m$ be given, and let $\pi\in \PF(\gamma)$ be the unique parking function with reading word $(n,n-1,\ldots,1)$.  Then $\dinv(\gamma):=\dinv(\gamma,\pi)$ counts the number of triples $(i,j,l)$ such that $i<_a j\in \{1,\ldots,n\}$, $l\in \{0,1,\ldots,m-1\}$, and 
\[
\gamma_i-l \leq \bar{\gamma}_j\leq \gamma_i-l+m-1.
\]
\end{corollary}
\begin{proof}
If $i$ attacks $j$, then $\pi_i<\pi_j$, since $\pi_j$ occurs before $\pi_i$ in the reading word $(n,\ldots,2,1)$.  The contribution of $(i,j)$ to $\dinv(\gamma,\pi)$ is then calculated from Lemma \ref{lemma:parkingDinv}. 
\end{proof}

\subsection{From parking functions to symmetric functions}
\label{subsec:quasisymmetric}
Recall that composition of $n$ is a sequence $(\a_1,\ldots,\a_r)$ with $\a_i\in \{1,\ldots,n\}$ and $\sum_i\a_i=n$.  We indicate this by writing $\a\vDash n$.

Let $S\subset \{1,\ldots,n-1\}$ be given.  Associated to $S$ there is a unique composition  $\a=(\a_1,\ldots,\a_r)$, characterized by
\[
S = \{\a_1, \a_1+\a_2,\ldots,\a_1+\cdots +\a_{r-1}\}.
\]
By definition, the \emph{monomial quasi-symmetric function} is the formal infinite sum
\[
M_\a(\xx) := \sum_{i_1<i_2<\ldots<i_r}x_{i_1}^{\a_1}\cdots x_{i_r}^{\a_r}.
\]
Here, $\xx = (x_1,x_2,\ldots)$ denotes a countable list of formal variables.  Given a composition $\b\vDash n$ and an index $1\leq i\leq r-1$, we obtain a coarser composition of the form $\a=(\b_1,\ldots,(\b_i+\b_{i+1}),\ldots,\b_r)$.  We say that $\b$ \emph{refines $\a$} if $\a$ is obtained from $\b$ by a sequence of such operations.  The \emph{Gessel quasi-symmetric function} \cite{Gessel84}, is by definition
\[
Q_\a(\xx):=\sum_{\b\text{ refines } \a} M_{\b}(\xx).
\]
If $\a$ is associated with a subset $S\subset \{1,2,\ldots,n-1\}$, then we also write $Q_S$ for $Q_\a$

Now, fix an $m$-parking function $(\gamma,\pi)$, and fix $a\in \{1,\ldots,n-1\}$.  Let $i$ and $j$ be the indices such that $\pi_i=a$ and $\pi_j=a+1$.  We say $a$ is a \emph{descent} of $(\gamma,\pi)$ if either of the following conditions holds:
\begin{itemize}
\item $\gamma_i<\gamma_j$, or
\item $\gamma_i=\gamma_j$, and $i<j$.
\end{itemize}
Let $\des(\pi,\gamma)\subset \{1,\ldots,n-1\}$ denote the set of descents (elsewhere this may also be denoted by $\ides(\gamma,\pi)$, since it is actually the descent set of the \emph{inverse} of the reading word of $(\gamma,\pi)$).

Finally we have the following functions.
\begin{definition}\label{def:DgammaSymFun}
For each $m$-Dyck path $\gamma$ of height $n$, let $D_\gamma(\xx;q,t)$ denote the following formal sum
\[
D_\gamma(\xx;q,t) := \sum_{\pi\in \PF(\gamma)} q^{\area(\gamma)}t^{\dinv(\gamma,\pi)} Q_{\des(\gamma,\pi)}(\xx).
\]
\end{definition}

These functions originally appeared in \cite{HHLRU05}, where it was shown that $D_\gamma(\xx;q,t)$ is a symmetric function for all $\gamma$ (in fact Schur positive), and it was conjectured that $\nabla^m e_n = \sum_{\gamma} D_\gamma(\xx;q,t)$.  This conjecture is now a theorem due to Carlson-Mellit \cite{CarMel-pp} and Mellit \cite{MellitRational-pp}. 
\begin{theorem}[Shuffle conjecture \cite{CarMel-pp,MellitRational-pp}]\label{thm:shuffleThm}
We have
\[
\nabla^m e_n = \sum_{\gamma} D_\gamma(\xx;q,t),
\]
where the sum is over $m$-Dyck paths of height $n$.
\end{theorem}

\subsection{Relation to our complexes}
\label{subsec:connectionToShuffle}

In this section we state and prove a connection between our complexes of Soergel bimodules $\XB \CB_v^{(m)}$ and the pieces of the shuffle conjecture.

Let $\sigma\in \Seq_m$ be given.  Then there is a unique $\gamma\in \Dyck_m$ of the form $\xi^h \sigma$, where $\xi$ is the operation of rotation (Definition \ref{def:rotation}).  In other words
\[
\gamma = (\sigma_i - r, \ldots, \sigma_n-r, \sigma_1-r-1,\ldots,\sigma_{i-1}-r-1),
\]
where $r=\min\{\sigma_1,\ldots,\sigma_n\}$ and $i=\min\{j\:|\:\sigma_j=r\}$.

We obtain a many-to-one surjective map $p:\Seq_m\rightarrow \Dyck_m$.  There is a right inverse $\iota:\Dyck_m\rightarrow \Seq_m$ defined by
\[
\iota(\gamma)=(m+\gamma_2,\ldots,m+\gamma_n, m-1).
\]
Suppose $\sigma$ extends $v$, and assume that $v\neq m^n$.  Let $r=r_\min$ denote $\min\{v_1,\ldots,v_n\}$, and let $i=i_\min$ denote the smallest index such that $v_i=r$.  Observe that
\[
\iota\circ p(\sigma) = (\sigma_{i+1}+(m-r),\ldots,\sigma_n+(m-r),\sigma_1+(m-r-1),\ldots,\sigma_{i-1}+(m-k-1),m-1),
\]
and this sequence extends $(m^{n-1},m-1)$.  We can also describe $\sigma^\ast:=\iota\circ p(\sigma)$ as the unique sequence extending $(m^{n-1},m-1)$ such that $p(\sigma^\ast)=p(\sigma)$.  We also have $\sigma= \xi^{n(m-r)-i}\sigma^\ast$.

\begin{definition}\label{def:corrections}
Retain notation as above.  Define the following numbers:
\[
\corr_q(v) := \sum_{j>i} (v_j-r)+\sum_{j<i} (v_j-r-1)
\]
and
\[
\corr_t(v) = \sum_{j<j', i<j' \atop
0\leq v_{i}\leq v_j\leq m-2} (m-1-v_j) +  \sum_{j<j', i<j' \atop
0\leq v_{i}+1\leq v_j\leq m-2} (m-1-v_j). 
\]
\end{definition}

\begin{lemma}\label{lemma:gsigma}
Let $\sigma\in \Seq_m$ be a sequence extending $v$, and let $g_\sigma(q,t,a)$ be the contribution of $\sigma$ to the sum on the right-hand side of (\ref{eq:knotStatesum}).  In other words,
\[
g_\sigma(q,t,a) = q^{|\sigma|-|v|}t^{d(\sigma)+e(v)}\prod_{i=1}^n(1+\chi_i(\sigma) a t^{-d(\sigma,i)}).
\]
Then
\[
g_{\sigma}(q,t,a) = q^{-\corr_q(v)}t^{-\corr_t(v)} g_{\iota\circ p(\sigma)}(q,t,a).
\]
\end{lemma}
\begin{proof}
Retain the notation as above.  Since $\sigma $ is obtained from $\iota\circ p(\sigma)$ by rotation, the proof of Theorem \ref{thm:stateSumKnot} shows that $g_\sigma$ and $g_{\iota\circ p(\sigma)}$ are related by an explicit monomial in $q$ and $t$.  It is routine to show that the monomial is as claimed in the statement.
\end{proof}

\begin{theorem}\label{thm:connectionWithShuffle}
Let $v\in \{0,1,\ldots,m\}^n$ be given, let $r=\min\{v_1,\ldots,v_n\}$, and let $i$ be the smallest index with $v_i=r$.  Then
\begin{equation}\label{eq:shufflePoincareSeries}
g_v^{(m)}(q,t,a) = q^{-\corr_q(v)}t^{-\corr_t(v)}\sum_{k=0}^n  \left\langle \sum_{\gamma}\nabla^m D_\gamma(\xx;q,t), h_ke_{n-k} \right\rangle a^k
\end{equation}
where the sum is over $m$-Dyck sequences $\gamma\in \Dyck_m$ such that
\begin{itemize}
\item if $v_j\leq m-1$ and $i+1\leq j\leq n$  then $\gamma_{j+1-i}=v_{j}-r$.
\item if $v_j\leq m-1$ and  $1\leq j\leq i-1$, then $\gamma_{j+1-i+n}=v_{j}-r-1$.
\item if $v_j=m$ and $i+1\leq j\leq n$  then $\gamma_{j+1-i}\geq m-r$.
\item if $v_j=m$  and  $1\leq j\leq i-1$, then $\gamma_{j+1-i+n}\geq m-r-1$.
\end{itemize}
Here $g_v^{(m)}(q,t,a)$ is the Poincar\'e series of $\HHH(\XB\CB_v^{(m)})$, as computed in Theorem \ref{thm:knotCase}.
\end{theorem}
\begin{proof}
Theorem \ref{thm:stateSumKnot} combined with Lemma \ref{lemma:gsigma} yields
\[
g_v^{(m)}(q,t,a)=q^{-\corr_q(v)}t^{-\corr_t(v)} \sum_{\sigma\in \Seq_m(v)} g_{\iota\circ p(\sigma)}^{(m)}(q,t,a).
\]
Each sequence $\iota\circ p(\sigma)$ on the right-hand side corresponds uniquely to a Dyck path $\gamma$ as in the statement, as is easily verified.  We will prove by direct computation that the coefficient of $a^k$ in $g_{\iota\circ p(\sigma)}^{(m)}(q,t,a)$ is equal to $\ip{D_\gamma(\xx;q,t),h_ke_{n-k}}$.  Note that this  computation is completely independent of $v$.

Thus, let $\sigma\in \Seq_m$ be a sequence extending $m^{n-1}(m-1)$.  Let
\[
\gamma=(0,\sigma_1-m,\ldots,\sigma_{n-1}-m)\in \Dyck_m
\]
be the associated $m$-Dyck sequence, hence $\iota(\gamma)=\sigma$.   Theorem 6.2.4 in \cite{HHLRU05} computes the right hand side of (\ref{eq:shufflePoincareSeries}), modulo the shuffle conjecture.  Thus,
\[
\ip{D_\gamma(\xx;q,t),h_ke_{n-k}} =  \sum_{\pi\in \PF_k(\gamma)} q^{\area(\gamma)}t^{\dinv(\gamma,\pi)},
\]
where $\PF_k(\gamma)\subset \PF(\gamma)$ is the set of parking functions whose reading word $(w_1,\ldots,w_n)$ contains $(n-k,\ldots,2,1)$ and $(n-k+1,\ldots,n-1,n)$ as subsequences.

In case $k=0$, a direct comparison of Lemma \ref{lemma:dinvLemma} and Corollary \ref{cor:parkingDinvTrivial} shows that the $a$-degree zero part of (\ref{eq:shufflePoincareSeries}) holds.

For general $k$, let $\pi\in \PF_k(\gamma)$, with reading word $w$.  Let $S(\pi)\subset \{1,\ldots,n\}$ denote the set of indices $i$ such that $\pi_i\in\{n-k+1,\ldots,n\}$.  The parking function $\pi$ is completely determined by $S(\pi)$.

\vskip10pt
\noindent
\textbf{Claim: }if $i\in S$, then $\gamma_{i+1}<\gamma_i+m$ (so that $\chi_i(\gamma)=1$ in the language of Theorem \ref{thm:stateSumKnot}).  Indeed, if $\sigma_{i+1}=\sigma_i+m$, then $\pi_i<\pi_{i+1}$ by definition of parking functions.  On the other hand, $\pi_{i+1}$ occurs earlier in the reading word for $\pi$.  Hence the pair $(\pi_{i+1}>\pi_i)$ is a subword of $w$, which is only possible if either $\pi_{i}<\pi_{i+1}\leq n-k$ or $\pi_i\leq n-k<\pi_{i+1}$ (otherwise $\pi_i$ and $\pi_{i+1}$ would a part of an increasing subsequence of $w$, which we know is impossible).   This proves the claim.
\vskip10pt

Thus, $\ip{D_\gamma(\xx;q,t),h_k e_{n-k}}$ can be expressed as a sum over subsets $S\subset \{i\:|\: \chi_i(\sigma)=1\}$ of size $l$.  Fix such a subset $S$, and let $(i,j)$ be a pair of indices so that $i$ attacks $j$.  Then $\pi_j$ comes before $\pi_i$ in the reading word.  Let us compute the contribution of $(i,j)$ to $\dinv(\gamma,S)$, working case by case.

First, assume that $i\not\in S$.  Then if $j\not\in S$, then $\pi_j,\pi_i\leq n-k$, which forces $\pi_j>\pi_i$ since $(\pi_j,\pi_i)$ occurs as a subword of $(n-k,\ldots,1)$.  If $j\in S$, then $\pi_j>\pi_i$ since $\pi_i\leq n-k$ and $\pi_j>n-k$.  Either way, $\pi_i<\pi_j$, and the contribution of $(i,j)$ to $\dinv(\gamma,\pi)$ is the same as the contribution of $(i,j)$ to $\dinv(\gamma)$ by Lemma \ref{lemma:parkingDinv}.

Next assume that $i\in S$.  We claim that $\pi_i>\pi_j$ is forced on us.  Indeed, if $j\in S$ then $\pi_j<\pi_i$ since in this case $(\pi_j,\pi_i)$ is a subword of $(n-k+1,\ldots,n)$.  If $j\not \in S$ then $\pi_j\leq n-k<\pi_i$.  Either way, $\pi_i>\pi_j$.  Using Lemma \ref{lemma:parkingDinv}, we now compute the contribution of $(i,j)$ to $\dinv(\gamma,\pi)$.

The contribution of $(i,j)$ to $\dinv(\gamma)$ is the number of solutions to the inequalities
\begin{equation}\label{eq:someIneq1}
\gamma_i\leq \bar{\gamma}_j+k\leq \gamma_i+m-1 \quad\quad 0\leq k\leq m-1,
\end{equation}
while the contribution of $(i,j)$ to $\dinv(\gamma,\pi)$ equals the number of solutions to the inequalities
\begin{equation}\label{eq:someIneq2}
\gamma_i\leq \bar{\gamma}_j+k+1\leq \gamma_i+m-1 \quad\quad 0\leq k\leq m-1
\end{equation}
where $\bar{\gamma}_j$ is as in Definition \ref{def:attacks}.  
These inequalities have the same number of solutions, unless $\gamma_i\leq \bar{\gamma}_j\leq \gamma_i+m-1$, in which case $k=m-(\bar{\gamma}_j-\gamma_i)-1$ is a solution of (\ref{eq:someIneq1}) but not (\ref{eq:someIneq2}).

We conclude that
\[
\dinv(\gamma,\pi) = \dinv(\gamma) - d(\gamma,i),
\]
where $d(\gamma,i)$ is the number of indices $j$ such that $\gamma_i\leq \bar{\gamma}_j\leq \gamma_i+m-1$ (note that this is exactly as in Definition \ref{def:dinv_i}).   Summing over all $k$ and all subsets $S$ gives
\[
\sum_{k=0}^n \ip{D_\gamma(\xx;q,t), h_k e_{n-k}}a^k= q^{\area(\gamma)} t^{\dinv(\gamma)}\prod_{i=1}^n (1+\chi_i(\gamma) at^{-d(\gamma,i)}),
\]
which equals $g_{\iota(\gamma)}(q,t,a)$, as claimed.
\end{proof}

\subsection{Special cases}
\label{subsec:specialCasesShuffle}
As a corollary of Theorem \ref{thm:connectionWithShuffle} and the $m$-shuffle conjecture (which we recorded in Theorem \ref{thm:shuffleThm}), we have the following.

\begin{corollary}
The symmetric function $\frac{1}{1-q}e_n$ is a combinatorial shadow of $\XB_n\FT_n^m$ in the sense of Definition \ref{def:GNRshadow}.  In particular, the Poincar\'e series of the $(n,nm+1)$ torus knot equals
\[
g_{m^n}^{(m)}(q,t,a) = \frac{1}{1-q} g_{m^{n-1},m-1}^{(m)} = \frac{1}{1-q}\left\langle \nabla^m  e_n, h_ke_{n-k} \right\rangle a^k
\]
\end{corollary}

We believe that the complexes $\CB_{m^{n-i},k,m^{i-1}}^{(m)}$ should also have interesting combinatorial shadows.  If $k=m-1$ then we expect that the Schur function associated to a hook shape provides such a combinatorial shadow (see Conjecture \ref{conj:hookShadow}), but for general $k$ we do not formulate precise conjecture in this direction.  In any case observe the following, which follows immediately from Theorem \ref{thm:connectionWithShuffle}.

\begin{corollary}\label{cor:generalizedHooks}
We have
\[
g_{m^{n-i},k,m^{i-1}}^{(m)}(q,t,a) = q^{-n(m-k)+(m-k)+(n-i)}\sum_{k=0}^n  \left\langle \sum_{\gamma} D_\gamma(\xx;q,t), h_ke_{n-k} \right\rangle a^k
\]
where the sum is over $\gamma\in\Dyck_m$ with 
\begin{itemize}
\item $\gamma_2,\ldots,\gamma_{i}\geq m-k$.
\item $\gamma_{i+1},\ldots, \gamma_n\geq m-k-1$.
\end{itemize}
\end{corollary} 

Taking $k=m-1$ we obtain
\begin{corollary}\label{cor:hookBraids}
We have
\[
 g_{1^{n-i},0,1^{i-1}}^{(m)}(q,t,a) = q^{1-i}\sum_{k=0}^n  \left\langle \sum_{\gamma}\nabla D_\gamma(\xx;q,t), h_ke_{n-k} \right\rangle a^k
\]
in which the sum is over Dyck sequences $\gamma\in \Dyck_m$ with $\gamma_2,\ldots,\gamma_{i}\geq 1$, where this condition is vacuous if $i=1$.
\end{corollary}

Finally, it is worth stating explicitly our result for $m=1$.
\begin{corollary}\label{cor:compositional}
Let $1^n\neq v\in \{0,1\}^n$ be a sequence, and let $i$ be the smallest index such that $v_i=0$.  Then
\[
g_v^{(1)} = q^{1-i} \left\langle\sum_{\gamma} D_\gamma(\xx;q,t)\right\rangle,
\]
where the sum is over Dyck sequences $\gamma\in \Dyck_1$ such that if  $i+1\leq j\leq n$  then $v_j=0$ and iff $\gamma_{j+1-i}=0$.
\end{corollary}

\subsection{Conjectures}
\label{subsec:moreconnections}

Our work suggests the following combinatorial conjectures.  We first introduce some terminology.  Given $v\in \{0,1\}^n$ with $v_1=0$, there is an associated composition $\a=\vDash n$ in a standard way.  To be precise, let$1:=i_0<i_1<\cdots<i_{r-1}\leq n$ denote the indices $i_j$ such that $v_{i_j}=0$, set $i_r:=n+1$, and define $\a_j:=i_{j}-i_{j-1}$ for $1\leq j\leq r$.  This describes a bijection between compositions of $n$ and binary sequences $v\in\{0,1\}^n$ with $v_1=0$.

 In \cite{HMZ12} the authors consider Hall-LIttlewood symmetric functions, indexed by binary seqeunces  $C_v(\xx;t)\in \L_t$, indexed by sequences $v\in \{0,1,\}^n$ with $v_1=0$, such that
\[
e_n = \sum_{\a \vDash n} C_v(\xx,t).
\]
See Proposition 5.2 in \cite{HMZ12}.  Note that we are swapping the roles of $t$ and $q$, relative to the preferred conventions in \cite{HMZ12}.  The authors go on to conjecture that $\nabla C_\a=\sum_\gamma D_\gamma(\xx;q,t)$, where the sum is over Dyck sequences $\gamma$ such that $\gamma_i=0$ iff $v_i=0$, where $v\in \{0,1\}^n$ is the binary sequence associated to $\a$.  Comparison with Corollary \ref{cor:compositional} suggests the following.

\begin{conjecture}\label{conj:HallLittlewoodShadow}
Let $v\in \{0,1\}^n$ be a sequence such that $v_1=0$, and let $\a\vDash n$ be the associated composition.  Then $\nabla_\a C_\a(\xx;t)$ is a combinatorial shadow of $\CB_v^{(1)}$, in the sense of Definition \ref{def:GNRshadow}.
\end{conjecture}

Proposition 5.3 in \cite{HMZ12} states that
\[
s_{i,1^{n-i}} = (-t)^{i-1} \sum_{\a\vDash n\atop \a_1\geq i} C_\a(\xx,t),
\]
which combined with the shuffle theorem gives
\[
\nabla s_{i,1^{n-i}} = (-t)^{i-1} \sum_{\gamma} D_\gamma(\xx,t),
\]
where the sum is over Dyck sequences $\gamma$ with $\gamma_2,\ldots,\gamma_{i}\geq 1$.    Comparing with Corollary \ref{cor:hookBraids} suggests the following.

\begin{conjecture}\label{conj:hookShadow}
The symmetric function $(-qt)^{1-i} s_{i,1^{n-i}}$ is a combinatorial shadow of the complex $\XB\CB_{1^{n-1}01^{i-1}}^{(1)}\in \KC^b(\SBim_n)$.  In particular
\begin{equation}\label{eq:KRpoly_hookConj}
g_{m^{n-i},m-1,m^{i-1}}^{(m)}(q,t,a) = (-qt)^{1-i}\sum_{k=0}^n  \left\langle \nabla^{m} s_{i,1^{n-i}}, h_ke_{n-k} \right\rangle a^k.
\end{equation}
Furthermore,
\begin{equation}\label{eq:mShuffle_hook}
\nabla^m s_{i,1^{n-k}} =  (-qt)^{i-1}\sum_{\gamma} D_\gamma(\xx,q,t)
\end{equation}
where the sum is over $\gamma\in \Dyck_m$ with  with $\gamma_2,\ldots,\gamma_{i}\geq 1$.   
\end{conjecture}

Equations (\ref{eq:KRpoly_hookConj}) and (\ref{eq:mShuffle_hook}) are consitent with one another, given Corollary \ref{cor:generalizedHooks}.

As is well-known, the power sum symmetric function $p_n\in \L_{q,t}$ equals $\sum_{i=1}^n (-1)^{n-i}s_{i,1^{n-i}}$.  Combining this with the previous conjectures yields:
\begin{conjecture}
We have
\[
(-1)^{n-1}\nabla^m p_n = \sum_{i=1}^n(qt)^{i-1}\sum_{\gamma} D_\gamma(\xx;q,t),
\]
where in the second summation, $\gamma\in\Dyck_m$ ranges over all $m$-Dyck sequences with $\gamma_2,\ldots,\gamma_{i}\geq 1$.  In particular
\[
\sum_{i=1}^n (-qt)^{i-1} g_{m^{n-i},(m-1),m^{i-1}}^{(m)}(q,t,a) = (-1)^{n-1}\sum_{k=0}^n \ip{\nabla^m p_n, h_k e_{n-k}}.
\]
\end{conjecture}

The above conjectures have focused on the knot case, i.e.~the polynomials $g_v$ rather than the \emph{rational functions} $f_v$.  For the link case, the trivial braid $\one_1$ will correspond to the symmetric function $\frac{1}{1-q}p_1$ under the GNR correspondence, hence one expects $\one_n = \one_1\sqcup \cdots \sqcup \one_1$ to correspond to $\left(\frac{1}{1-q}\right)^n p_1^n$.  This gives rise to the following.

\begin{conjecture}
$\FT_n = \CB^{(1)}_{1^n}$ has a combinatorial shadow given by $p_1^n$. In particular
\[
f_{m^n}^{(m)} = \left(\frac{1}{1-q}\right)^n\ip{\nabla^m p_1^n, h_ke_{n-k}}a^k.
\]
\end{conjecture}

In case $m=1$ this reduces to a conjecture of Andy Wilson's \cite{Wilson16-pp}.

\printbibliography

\end{document}